\documentclass[11pt]{amsart}
\usepackage{amsmath,amsthm, amscd, amssymb, amsfonts}
\usepackage{textcomp}
\usepackage[all]{xy}
\usepackage{enumitem}

\usepackage[latin1]{inputenc}
\usepackage[dvips, dvipsnames, usenames]{color}
\numberwithin{equation}{section}\theoremstyle{plain}
\usepackage{color}

\newtheorem*{stepi}{Step i}
\newtheorem{theorem}{Theorem}[section]

\newtheorem{lem}[theorem]{Lemma}
\newtheorem{cor}[theorem]{Corollary}

\newtheorem{pro}[theorem]{Proposition}
\newtheorem{claim}{Claim}[section]

\theoremstyle{definition}

\newtheorem{exa}[theorem]{Example}

\theoremstyle{remark}
\newtheorem{rem}[theorem]{Remark}

\renewcommand{\_}[1]{_{\left( #1 \right)}}

\newcommand{\uvi}{\mathtt u}

\newcommand{\ydh}{{}^{H}_{H}\mathcal{YD}}
\newcommand{\ydhp}{{}^{H'}_{H'}\mathcal{YD}}

\newcommand{\hyd}{\mathcal{YD}^H_H}

\newcommand{\ydg}{{}^{\Gamma}_{\Gamma}\mathcal{YD}}
\newcommand{\ydkg}{{}^{\k\Gamma}_{\k\Gamma}\mathcal{YD}}

\def\zt{\Z^{\theta}}

\def\G{\mathbb{G}}
\def\Gb{\mathbf{G}}
\def\Lb{\mathbf{L}}

\newcommand\id{\operatorname{id}}
\newcommand\Id{\operatorname{Id}}

\newcommand\ord{\operatorname{ord}}
\newcommand\Hom{\operatorname{Hom}}
\newcommand\Alg{\operatorname{Alg}}
\newcommand\GL{\operatorname{GL}}

\newcommand\End{\operatorname{End}}

\newcommand\simple{\mathfrak S}
\def\qb{\mathfrak{q}}

\newcommand\gr{\operatorname{gr}}
\newcommand\co{\operatorname{co}}

\newcommand\ad{\operatorname{ad}}

\newcommand\Irr{\operatorname{Irr}}

\def\Sym{\mathbb{S}}

\newcommand{\Ss}{{\mathcal S}}
\def\k{\Bbbk}
\newcommand\Isom{\operatorname{Isom}}

\def\ot{\otimes}

\def\s{\mathbb{S}}
\def\C{\mathbb{C}}
\def\Z{\mathbb{Z}}
\def\N{\mathbb{N}}

\def\B{\mathfrak{B}}
\def\R{\mathfrak{R}}

\def\eps{\epsilon}
\def\mT{\mathcal{T}}

\def\mB{\mathcal{B}}
\def\mI{\mathcal{I}}

\def\mL{\mathcal{L}}
\def\mQ{\mathcal{Q}}
\def\Oc{\mathcal{O}}
\def\mR{\mathcal{R}}

\def\mE{\mathcal{E}}
\def\mA{\mathcal{A}}
\def\mH{\mathcal{H}}

\def\bH{\mathtt{H(p)}}
\newcommand{\ug}{\mathfrak{u}} 
\def\bs{\boldsymbol}
\def\tt{\texttt{t}}

\newcommand{\toba}{{\mathcal B}}

\newcommand{\tobat}{\widetilde{\toba}}
\newcommand{\tobah}{\widehat{\toba}}
\newcommand{\pit}{\widetilde{\pi}}

\def\mA{\mathcal{A}}
\def\mE{\mathcal{E}}
\def\mH{\mathcal{H}}
\def\mL{\mathcal{L}}
\def\mT{\mathcal{T}}
\def\mtA{\widetilde{\mA}}
\def\mhA{\widehat{\mA}}
\def\mtE{\widetilde{\mE}}
\def\mhE{\widehat{\mE}}
\def\mtH{\widetilde{\mH}}
\def\mhH{\widehat{\mH}}
\def\mtL{\widetilde{\mL}}
\def\mhL{\widehat{\mL}}
\newcommand\I{\mathbb I}

\def\tz{\zeta}
\newcommand{\J}{{\mathcal J}}
\newcommand{\Gc}{{\mathcal G}}

\newcommand{\up}{\delta}

\def\pf{\begin{proof}}
\def\epf{\end{proof}}

\def\q{{\bf q}}

\def\lg{\langle}
\def\rg{\rangle}

\newcommand{\prov}{{\mathcal H}}

\newcommand{\prova}{{\mathcal A}}
\newcommand{\cA}{{\mathcal A}}

\newcommand\w{\widetilde}

\allowdisplaybreaks

\newcommand\Cleft{\operatorname{Cleft}}

\begin{document}
 \title[Pointed Hopf algebras of type A]{\small  Liftings of Nichols algebras of diagonal type \\ I. Cartan type A}
\author[Andruskiewitsch; Angiono; Garc\'ia Iglesias]{Nicol\'as Andruskiewitsch;
Iv\'an Angiono; Agust\'in Garc\'ia Iglesias}

\address{FaMAF-CIEM (CONICET), Universidad Nacional de C\'ordoba,
Medina A\-llen\-de s/n, Ciudad Universitaria, 5000 C\' ordoba, Rep\'
ublica Argentina.} \email{(andrus|angiono|aigarcia)@famaf.unc.edu.ar}

\thanks{The work was partially supported by CONICET,
FONCyT-ANPCyT, Secyt (UNC), the MathAmSud project
GR2HOPF}

\begin{abstract} 
After the classification of the finite-dimensional Nichols algebras of diagonal type \cite{H-inv,H-classif}, the determination of its defining relations \cite{A-jems,Ang-crelle}, and the verification of the \emph{generation in degree one} conjecture \cite{Ang-crelle},
there is still one step missing in the classification of complex finite-dimensional Hopf algebras with abelian group, without restrictions on the order of the latter:
the computation of all deformations or liftings. A technique towards solving this question was developed in \cite{AAnGMV}, built on cocycle deformations.
In this paper, we elaborate further and present an explicit algorithm to compute liftings. In our main result we classify  all liftings of finite-dimensi\-onal Nichols algebras of Cartan type $A$, over a cosemisimple Hopf algebra $H$. This extends \cite{AS2}, where it was assumed that the parameter is a root of unity of order $>3$ and that $H$ is a commmutative group algebra. When the parameter is a root of unity of order 2 or 3, new phenomena appear:  
the quantum Serre relations can be deformed;  this allows in turn the power root vectors to be deformed to elements in  lower terms of the coradical filtration, but not necessarily in the group algebra. 
These phenomena are already present in the calculation of the liftings in type $A_2$ at a parameter of order 2 or 3 over an abelian group \cite{BDR,helbig},
done by a different method using a computer program. As a by-product of our calculations, we present new infinite families of finite-dimensional pointed Hopf algebras.
\end{abstract}

\maketitle

\section{Introduction}

\subsection{The general context}\label{subsec:gral-context}
This is the first article of a series  intended to determine all liftings of 
finite-dimensi\-onal Nichols algebras of diagonal type over an algebraically closed field of characteristic zero $\k$.
The end of this series will also conclude
the classification of the finite-dimensi\-onal pointed Hopf algebras with abelian group of group-likes, without
restrictions on the order of the group.
The setting, slightly different than in \cite{AAnGMV}, is the following. We fix:
\begin{itemize} [leftmargin=*]\renewcommand{\labelitemi}{$\circ$}
\item A cosemisimple Hopf algebra $H$.
\item A braided vector space  of diagonal type $(V,c)$, with a principal realization
 in $\ydh$, such that the Nichols algebra $\B(V)$ is finite-dimensional.
\end{itemize}
We place ourselves in this more general context in order to contribute to the classification of
Hopf algebras with finite Gelfand-Kirillov dimension, and more precisely to those that are co-Frobenius; see
\cite{AAH}.

A \emph{lifting} of  $V \in\ydh$ is a Hopf
algebra $L$ such that $\gr L= \B(V) \# H$, where $\gr L$ is the graded
Hopf algebra associated to the coradical filtration. In other words \cite[2.4]{AV}, $L$ is a lifting of $V$ iff
there is an epimorphism of Hopf algebras  $\phi: \mT(V) := T(V)\# H
\rightarrow L$  such that $\phi_{|H}=\id_H$ and 
\begin{align} \label{eq:lifting-map}
  \phi_{|H \oplus V\# H}: H \oplus V\#H \rightarrow L_1 \text{ is an isomorphism of Hopf bimodules}.
\end{align}
Such $\phi$ is called a \emph{lifting map}.
If emphasis on $H$ is needed, then we say that $L$ is a lifting of $V$ over $H$; if $H = \k G$
is the group algebra of the group $G$, then we also say that $L$ is a lifting of $V$ over $G$.

The aim of the series is to compute all liftings of every $V$ as above. 
It seems very hard, and probably not feasible, to give a uniform answer to this problem, {\it i.e.} compact formulae valid for all $V$.
We proceed then by a case-by-case analysis of the list in the classification of \cite{H-classif}.
Let $r_1, \dots, r_M$ be the defining (homogeneous) relations of $\B(V)$,  computed in \cite{Ang-crelle}; 
let $n_j = \deg r_j$. If $\phi$ is a lifting map as above, then there exists $p_j \in \bigoplus\limits_{0\le i < n_j} 
T^i(V)\# H$
such that $\phi(r_j) = \phi(p_j)$, for all $j$. Our approach is:

\begin{itemize} [leftmargin=*]\renewcommand{\labelitemi}{$\Diamond$}
	\item To establish the general form of the $p_j$'s, in terms of the $r_j$'s and some parameters.
	\item  To define a Hopf algebra $L = \mT(V) /\langle r_1 - p_1, \dots, r_M - p_M \rangle$ for each choice of the parameters alluded above
	and to prove that $\gr L \simeq \B(V) \# H$.
	\item To show that every lifting can be obtained in this way.
\end{itemize}

In situations considered in previous work \cite{AS2, AS3} the $\phi(r_j)$'s belong to $H$ and were computed recursively, while 
the remaining points were dealt with by ad-hoc manners. In this series we proceed recursively again but following the strategy in \cite{AAnGMV}, inspired by \cite{GM1, M0}; namely, we compute a sequence of quotients of $\mT(V)$ as cocycle deformations 
of a parallel sequence of quotients of the form $\B \# H$,  describing eventually  $L = \mT(V) /\langle r_1 - p_1, \dots, r_M - p_M \rangle$
as a cocycle deformation of  $\B(V) \# H = \mT(V) /\langle r_1, \dots, r_M \rangle$. See Section \ref{sec:strategy}.

\smallbreak
In this paper we compute all liftings for $V$ of Cartan type $A$, over a root of unit $\xi$ of order 2 or 3. The case when $\xi$ has order $>3$ is known for group algebras of finite abelian groups \cite[\S 6]{AS2};
we  extend this to a general cosemisimple Hopf algebra. See Theorems \ref{thm:main-N=2}, \ref{thm:main-N=3}  and \ref{thm:main-N>3}.
 There are three reasons to start with Cartan type $A$. First, it is the Dynkin diagram of \emph{Her all-embracing Majesty} \cite{W}. Second, formulae for the Nichols algebras of this type are much more explicit than for other types. Third, the experience and results for this type would help to understand and solve the other types.

\subsection{The main result} Let $\theta \in \N$ and $\I = \{1, \dots, \theta\}$. 
Let $V$ be a braided vector space of diagonal type with basis $(x_i)_{i\in \I}$ and braiding matrix  $\qb = (q_{ij})_{i,j\in \I}$.
Let  $(\alpha_i)_{i\in \I}$ be the canonical basis of $\zt$. 
The braided Hopf algebra $T(V)$ is $\zt$-graded by $\vert x_i\vert = \alpha_i$, $i\in \I$.
Let $\chi: \zt \times \zt \to \k$ be the bilinear form defined by $\chi(\alpha_i, \alpha_j) = q_{ij}$, $i,j\in \I$;
set $q_{\alpha \beta} = \chi(\alpha, \beta)$, $\alpha, \beta \in \zt$.
The braided commutator is defined on $\zt$-homogeneous elements $u, v \in T(V)$  by
\begin{align*}
[u, v]_c &= uv - q_{\vert u\vert \vert v\vert} vu.
\end{align*}
Then $\ad_c x_i(v) := [x_i, v]_c$. 
Let $\xi$ be a primitive $N$-th root of unity, $N \ge 2$.
We fix a  braiding matrix  $(q_{ij})_{i,j\in \I}$ such that
\begin{align}\label{eq:braiding-A}
q_{ii}&=\xi, & q_{ij}q_{ji}&=\begin{cases}
\xi^{-1}, & |i- j|=1,\\
1, & |i- j|>1, \end{cases} & i,j &\in \I.
\end{align}
This is a braided vector space of Cartan type $A_\theta$ and the corresponding generalized Dynkin diagram, cf. \cite{H-classif}, is $\xymatrix{ \overset{\xi}{\circ}\ar  @{-}[r]^{\xi^{-1}}  &
	\overset{\xi}{\circ}\ar
	@{-}[r]^{\xi^{-1}} &  \overset{\xi}{\circ}\ar@{.}[r] & \overset{\xi}{\circ} \ar
	@{-}[r]^{\xi^{-1}}  & \overset{\xi}{\circ}}$.
The corresponding Nichols algebra is indeed the multiparametric version of the positive part of the small quantum group, or Frobenius-Lusztig kernel, of type  $A_\theta$.
For  $i \le j\in \I$, we denote $(i\, j ) = \sum_{i \le k\le j} \alpha_k \in \zt$; clearly $\{(i\, j ):i \le j\in \I \}$ is the set of positive roots of the root system $A_\theta$. The associated Lyndon words are defined recursively by
\begin{align*}
& x_{(i\, j )}=\begin{cases} x_{j}, & i=j, \\
[x_i, x_{(i+1 \, j)}]_c & i+1 \leq j,
\end{cases}& 
\end{align*}
in $T(V)$ or any quotient thereof. We also need the notation $x_{ij} = [x_i, x_j]_c$, $i < j\in \I$.
We now state the presentation of $\toba(V)$ by generators and relations.
Part \eqref{item:presentation-N>2}  was proved in \cite{AS2}, inspired by 
\cite{T}; \eqref{item:presentation-N=2} is from \cite{AD}.
\begin{pro}\label{prop:presentation-nichols}
\begin{enumerate}[leftmargin=*]
\item\label{item:presentation-N>2}  Assume that $N > 2$. Then $\B(V)$ is  generated by
$(x_i)_{i\in \I}$ with relations
\begin{align}
\label{eq:rels-Atheta-a}
\ x_{ij} &= 0,&  & i < j - 1;\\
\label{eq:rels-Atheta-b3}
(\ad_c x_i)^{2}(x_{j}) &= 0,& &\vert j- i\vert = 1;
\\
\label{eq:rels-Atheta-c}\ x_{(i\,j)}^N &=0,& & i \leq j.
\end{align}
The distinguished  pre-Nichols algebra $\tobat(V)$ \cite[Definition 1]{Ang}  is generated
by $(x_i)_{i\in \I}$ with relations \eqref{eq:rels-Atheta-a} and
\eqref{eq:rels-Atheta-b3}; this  is denoted $\widehat{\toba}(V)$ in  \cite[\S 6.3]{AS2}.

\medbreak
\item\label{item:presentation-N=2}  Assume that $N = 2$. Then $\B(V)$ is  generated by
$(x_i)_{i\in \I}$ with relations \eqref{eq:rels-Atheta-a}, \eqref{eq:rels-Atheta-c} and 
\begin{align}\label{eq:rels-Atheta-b2}
[x_{(i-1\,i+1)}, x_i]_c &= 0, & & 1<i<\theta.
\end{align}
The distinguished pre-Nichols algebra $\tobat(V)$ \cite[Definition 1]{Ang}  is generated
by $(x_i)_{i\in \I}$ with relations \eqref{eq:rels-Atheta-a},
\eqref{eq:rels-Atheta-b3} and \eqref{eq:rels-Atheta-b2}. 
\end{enumerate}	
\end{pro}

\begin{rem}\label{rem:distinguished}
Relations \eqref{eq:rels-Atheta-b2}  hold for $N> 2$, by \eqref{eq:rels-Atheta-a} and \eqref{eq:rels-Atheta-b3}.

When $N=2$, \eqref{eq:rels-Atheta-b3}  becomes 
\begin{align*}
x_i^2x_j+q_{ij}^2x_jx_i^2&=0 ,& &\vert j- i\vert = 1.
\end{align*}
Since $x_i^2=0$ by \eqref{eq:rels-Atheta-c},  \eqref{eq:rels-Atheta-b3} holds in $\B(V)$; thus $\B(V)$ is a quotient of $\tobat(V)$.
\end{rem}

\begin{rem}
The distinguished pre-Nichols algebra $\tobat(V)$ is meant to have the same set of PBW generators, hence the same root system, as $\toba(V)$.
By this reason, the choice of the defining relations is performed so as to guarantee this property. 
In particular, one needs relations to reorder any pair of PBW generators.

Assume $V$ is of Cartan type $A$. If $N>3$, then this is automatically attained provided that the quantum Serre relations \eqref{eq:rels-Atheta-a} and \eqref{eq:rels-Atheta-b3} hold, see \cite[Lemmas 6.4 \& 6.7]{AS2}.
When $N=2$, then quantum Serre relations \eqref{eq:rels-Atheta-a} and \eqref{eq:rels-Atheta-b3} are not enough, as we can not reorder the PBW generators $x_{(i-1\,i+1)}$ and $x_i$, $1<i<\theta$; hence the need of \eqref{eq:rels-Atheta-b2}. Now, this enlarged set of relations suffices, as it is shown in Lemma \ref{lem:rels preNichols}.
\end{rem}

Assume that  $N>3$. Then all  liftings of  $V$ (over a finite abelian group) are classified in \cite[Theorem 6.25]{AS2}.
In this paper we classify all liftings of $V$ when $N = 2$  or $3$. 
To present our main results, we need more notation. Let $\left(g_i, \chi_i\right)_{i \in \I}$ be a principal realization of $V$ over $H$, see \S  \ref{subsec:ppal-realization}; let 
\begin{align*}
\Gamma &= \lg g_1,\dots,g_\theta\rg.
\end{align*}
For $i_1,\dots,i_k\in\I$ distinct, $k\in\N$, set
\begin{align*}
g_{i_1,\dots, i_k}&:=g_{i_1}\dots g_{i_k}, &  \chi_{i_1,\dots, i_k}&:=\chi_{i_1}\dots \chi_{i_k}, &x_{i_1,\dots, i_k}&:=[x_{i_1},[x_{i_2\dots, i_k}]_c]_c.; \\
g_{(i\,j)}&:=g_{i,i+1,\dots, j}, &  \chi_{(i\,j)}&:=\chi_{i,i+1,\dots, j}, &&  i 
\leq j\in\I.
\end{align*}
Also, if $i < j\in\I$, then let us fix $g_{(j\,i)}:=1$, $\chi_{(j\,i)}:=\eps$.

\subsubsection{Component in $\Gamma$}\label{sec:component} Here $N\geq 2$ is arbitrary. For $i\leq j\in \I$ we set
\begin{align}\label{eqn:Cp}
C_p=C^{j}_{ip} = (1-q^{-1})^N\chi_{(i\,p)}(g_{(p+1\,j)})^{N(N-1)/2}.
\end{align}
If the quantum Serre relations \eqref{eq:rels-Atheta-a} and \eqref{eq:rels-Atheta-b3} 
are not deformed, then the lifting problem is equivalent to the following question, which amounts to solving an equation in $\k\Gamma$, see 
\cite[(6-36)]{AS2}, \cite[\S 3]{AD}:
\begin{itemize}[leftmargin=*]\renewcommand{\labelitemi}{$\circ$}
\item Find all families $(u_{(i\,j)})_{i\leq j\in\I}$ of elements in $\k\Gamma$, such that
\begin{align}\label{eqn:sol-in-Gamma}
\Delta(u_{(i\,j)})=u_{(i\,j)}\ot1+g_{(i\,j)}^N\ot u_{(i\,j)}+\sum_{i\leq p<j} C^{j}_{ip} u_{(i\,p)}g_{(p+1\,j)}^N\ot u_{(p+1\,j)}.
\end{align}
\end{itemize} 
The solutions to \eqref{eqn:sol-in-Gamma} are given in \cite[Theorem 6.18]{AS2}. 
These are defined recursively on $j-i\geq 0$  \cite[6-40]{AS2} as elements $u_{(i\,j)}(\boldsymbol{\gamma})$, for each family\footnote{The parameters $\gamma_{ij}$ are called $\mu_{ij}$ (more precisely $\mu_\alpha$, $\alpha$ a root) in \cite{AS3}. This is the notation we shall adopt in this article.} of scalars $
\boldsymbol{\gamma} = (\gamma_{ij})_{i\leq j\in\I}$,   by 
\begin{align}\label{eqn:part-in-Gamma}
u_{(i\,j)}(\boldsymbol{\gamma}) = \gamma_{ij}(1-g_{(i\,j)}^N)+\sum_{i\leq p<j}C^{j}_{ip}\gamma_{ip} \, u_{(p+1\,j)}(\boldsymbol{\gamma}).
\end{align}
If $(u_{(i\,j)}(\bs\gamma))_{i\leq j\in\I}$ is a solution, then the quotient of $T(V)\#H$ by the ideal generated by
\begin{align}\label{eqn:lift-sin-qs}
\begin{split}
r&=0, \qquad  r \text{ (generalized) quantum Serre relation};\\
a_{(i\,j)}^N&=u_{(i\,j)}(\bs\gamma), \qquad  i\leq j\in\I,
\end{split}
\end{align}
is a lifting of $V$, by Theorems \ref{thm:main-N=2}, \ref{thm:main-N=3} and \ref{thm:main-N>3}.
It was shown in \cite[6.25]{AS2} that all liftings arise like this if $H$ is a commutative group algebra and $N>3$. In Theorem \ref{thm:main-N>3} we extend this to any cosemisimple $H$. We also compute all liftings when $N\geq 3$.  

\smallbreak A key difference in the case $N\leq 3$ is that 
solutions to \eqref{eqn:sol-in-Gamma} are a {\it part} of the general solution, see \eqref{eqn:lift-con-qs-3}  below. 
In particular, we show that the deformations do not necessarily restrict to the coradical. See for instance the concrete Examples \ref{ej:not-in-H-N2} and \ref{ej:not-in-H-N3}.

\begin{rem}\label{lem:Gamma}
In the present article we use an equivalent version of \eqref{eqn:part-in-Gamma}. Namely, we consider families of scalars 
$\boldsymbol\mu = (\mu_{(k\, l)})_{k\le l\in\I}$
subject to 
\begin{align}\label{eqn:cond-mu-gral}
\begin{aligned}
\mu_{(k\, l)}&=0,&  &\mbox{ if }\chi_{(kl)}^N\neq\eps,&
&\mbox{ or }g_{(kl)}^N = 1. 
 \end{aligned}
\end{align}
We define recursively $\uvi_{(j\,k)}=\uvi_{(j\,k)}(\bs\mu) \in \k \Gamma$, $j \leq  k\in \I$, by $\uvi_{(jj)}=0$, and  
\begin{align}\label{eqn:def-gral}
\uvi_{(j\,k)}&=-\sum_{j\leq p<k} C_p\mu_{(p+1\,k)}
\Big(\uvi_{(j\,p)}+\mu_{(j\,p)}(1-g_{(j\,p)}^N)\Big)g_{(p+1\,k)}^N.
\end{align}
The comparison with the previous solution is as follows: define
 $\bs\gamma=\bs\gamma(\bs\mu)$, $\bs\gamma=(\gamma_{ij})_{i\leq j\in I}$, by
$\gamma_{ij}=\mu_{(i\,j)}-\sum_{i\leq p<j} C_p\gamma_{ip}\mu_{(p+1\,j)}$,  $i\leq j$.
Then
\begin{align*}
u_{(jk)}(\bs\gamma) = \uvi_{(j\,k)}(\bs\mu) + \mu_{(j\,k)}(1-g_{(j,k)}^N).
\end{align*}
\end{rem}

\subsubsection{The shape of the liftings}
In the general case $N\geq 2$, we show that the lifting problem  is equivalent to solving an algorithm, described synthetically in \S \ref{sec:algorithm}. An  equation similar to \eqref{eqn:part-in-Gamma} must be solved recursively, this time with solutions in the previous term of the coradical filtration. 
We show for type $A_\theta$ in Theorems \ref{thm:main-N=2}, \ref{thm:main-N=3} and \ref{thm:main-N>3} that any lifting of $V$ is given by
\begin{itemize}
\item[i.] a solution $(u_{(i\,j)}(\bs\mu))_{i\leq j\in\I}$.
\item[ii.] elements $v_r(\bs\lambda)\in\k\Gamma$, one for each (generalized) quantum Serre relation and associated to  scalars $\bs\lambda=(\lambda_r)_r$. See \eqref{eq:rels-liftings prenichols-N2-1-intro},  \eqref{eq:rels-liftings prenichols-N2-3-intro}, \eqref{eq:rels-liftings prenichols-N3-2-intro}.
\item[iii.] elements $\sigma_{(i\,j)}(\bs\lambda,\bs\mu)\in T(V)\#H$,  computed algorithmically.
\end{itemize}
The corresponding lifting is the quotient of $T(V)\#H$ by
\begin{align}\label{eqn:lift-con-qs-3}
\begin{split}
r&=v_r(\bs\lambda), \quad  r \text{ (generalized) quantum Serre relation};\\
a_{(i\,j)}^N&=u_{(i\,j)}(\bs\mu)+\sigma_{(i\,j)}(\bs\lambda,\bs\mu), \quad i\leq j\in\I.
\end{split}
\end{align}
Compare with \eqref{eqn:lift-sin-qs}. 
When $N>3$, $\lambda_r= 0$ for all $r$ and thus $v_r(\bs\lambda)=0$, also $\sigma_{(i\,j)}(\bs\lambda,\bs\mu)=0$.
When $N=3$, $\lambda_r\neq 0$ only for $r$ of type \eqref{eq:rels-Atheta-b3} as relations \eqref{eq:rels-Atheta-a} remain unchanged. 

The case $N=2$ is actually a bit more involved, as the deformation of the generalized quantum Serre relations \eqref{eq:rels-Atheta-b2} depends on the deformation  of the powers of  the simple root vector relations, see \eqref{eq:rels-liftings prenichols-N2-3-intro}. 
Also, in Theorem \ref{thm:main-N=2}, the last line of \eqref{eqn:lift-con-qs-3} is expressed as
$\zeta_{(i\,j)}^2 =u_{(i\,j)}$, $i\leq j\in\I$,
as $\zeta_{(i\,j)}^2=a_{(i\,j)}^2+$ terms $\sigma(\bs\lambda,\bs\mu)$, see Remark \ref{rem:zeta}. 
The family $\bs\nu = (\nu_i)_{i\in\I}$ controls the deformations of the generalized quantum Serre relations.

\subsubsection{The main result, $N=2$}\label{subsubsec:N=2} Here $\xi =-1$.
We fix a family of scalars 
$\bs\mu=(\mu_{(k\,l)})_{k\leq l\in\I}$
subject to the constraints and normalizations \eqref{eqn:cond-mu-gral}.
We consider two more families of scalars 
\begin{align*}
\boldsymbol\lambda &=(\lambda_{ij})_{i<j-1\in \I},& 
\boldsymbol\nu &=(\nu_{i})_{1<i<\theta}
\end{align*}
subject to the constraints and normalizations
\begin{align}\label{eq:condiciones escalares N=2-intro}
\begin{aligned}
\lambda_{ij}&=0, & &\mbox{ if }\chi_{ij} \neq\eps,&
 &\mbox{ or }g_{ij} = 1; \\
 \nu_{i}&=0, & &\mbox{ if }\chi_{i-1,i,i,i+1} \neq\eps&
 &\mbox{ or }g_{i-1,i,i,i+1} = 1.
 \end{aligned}
\end{align}
We define families of elements in $\mT(V)$ attached to these parameters in the following way. To distinguish from the sequence of pre-Nichols algebras, we denote now by $(a_{i})_{i\in \I}$ the generators of $T(V)$; correspondingly, we denote $a_{i j}$, $a_{(i j)}$, $a_{i_1,\dots, i_k}$, instead of $x_{i j}$, $x_{(i j)}$, $x_{i_1,\dots, i_k}$.

\smallbreak
Let $i, j \in \I$, $|i-j|\ge2$. We define recursively scalars $d_{i\,j}(s)$, $b_{i\,j}(s)$, $s\in \N_0$, as follows:
$d_{i\,j}(0)=2\lambda_{i\,j}$, $b_{i\,j}(0)=-2\chi_j(g_{(i\,j)})\lambda_{i\,j}$, and for $s>0$,
\begin{align}\label{eq:corchete yi yjk - formula coef}
d_{i\,j}(s)&= q_{ij} \sum_{0 \le t <s} d_{i\,j+1}(t)d_{j\,j+2t+2}(s-t-1),\\
b_{i\,j}(s)&= \sum_{0 \le t <s} b_{i+1\,j}(t)d_{i\,i+2t+2}(s-t-1).
\end{align}

We define recursively $\tz_{(j\,k)}\in\mT(V)$ as follows: $\tz_{(j\,j)}=a_j$ and for $j<k$
\begin{multline}\label{eq:def ztilde-intro}
\tz_{(j\,k)} = [a_{j},\tz_{(j+1\,k)}]_c+ d_{jk}(0)\chi_{(j\, k)}(g_j) \, \tz_{(j+1\,k-1)} g_{jk}\\ + 2 \sum_{1 \le t \le (k-j-1)/2 } d_{j k-2t}(t) 
\chi_{(j+1\, k-2t-1)}(g_{j}) \tz_{(j+1\,k-2t-1)} g_{j}g_{(k-2t\,k)}.
\end{multline}

Let $\ug(\boldsymbol\lambda,\bs\mu,\bs\nu)$ be the quotient of $\mT(V)$
by the relations
\begin{align}
a_{ij} &=\lambda_{ij}(1-g_ig_j); \label{eq:rels-liftings prenichols-N2-1-intro}\\
[a_{(i-1\,i+1)},a_i]_c &= \nu_{i}(1-g_i^2g_{i-1}g_{i+1}) \label{eq:rels-liftings prenichols-N2-3-intro} \\
& \qquad -4\chi_i(g_{i-1})\mu_{(i)}\lambda_{i-1\,i+1}g_{i-1}g_{i+1}(1-g_i^2);	\notag
\\
\label{eq:rels-liftings-N2-intro}
\tz_{(j\,k)}^2&=\mu_{(j\,k)}(1-g_{(j\,k)}^2)+ \uvi_{(j\,k)},
\end{align}
for $\uvi_{(j\,k)}=\uvi_{(j\,k)}(\bs\mu)$ as in \eqref{eqn:def-gral}.

The relations \eqref{eq:rels-liftings prenichols-N2-1-intro}
	are deformations of  \eqref{eq:rels-Atheta-a}, while \eqref{eq:rels-liftings-N2-intro} are deformations of \eqref{eq:rels-Atheta-c}, and \eqref{eq:rels-liftings prenichols-N2-3-intro} are deformations of \eqref{eq:rels-Atheta-b2}.

\begin{rem} The quotient $\widetilde\ug(\boldsymbol\lambda,\bs\mu,\bs\nu)$  of $\mT(V)$ by the relations \eqref{eq:rels-liftings prenichols-N2-1-intro},   \eqref{eq:rels-liftings prenichols-N2-3-intro} and \eqref{eq:rels-liftings-N2-intro} for $j=k$ is a  cocycle deformation of $\tobat(V) \# H$.
\end{rem}

Recall that that $V$ is of type $A_\theta$ at $\xi =-1$.

\begin{theorem}\label{thm:main-N=2}
The algebra $\ug(\boldsymbol\lambda,\bs\mu,\bs\nu)$ is a Hopf algebra quotient of $\mT(V)$ and is a lifting of $V$. 
Reciprocally every lifting of $V$ over $H$ is isomorphic to 
$\ug(\boldsymbol\lambda,\bs\mu,\bs\nu)$ for some family of scalars $\boldsymbol\lambda$, $\bs\mu$, $\bs\nu$ as 
in \eqref{eq:condiciones escalares N=2-intro}.
In particular, every lifting is a cocycle deformation of $\toba(V) \# H$.
\end{theorem}
\pf
We follow the strategy in \S \ref{sec:strategy}: If  $\mH=\toba(V) \# H$, then 
$\ug=\ug(\boldsymbol\lambda,\bs\mu,\bs\nu)$ arises as $L(\mA,\mH)$ for a given 
$\mA=\mA(\boldsymbol\lambda,\bs\mu,\bs\nu)\in\Cleft\mH$ such that $\gr\ug\simeq \mH$. The corresponding {\it stratification} cf. \S \ref{sec:cleft} of the set of generators of the ideal defining $\toba(V)$ is given by 
 $\Gc_0=\{\eqref{eq:rels-Atheta-a}, \eqref{eq:rels-Atheta-b2}, x_i^2, i \in \I\}$, $\Gc_1=\{\eqref{eq:rels-Atheta-c}\}$. The converse follows by Theorem \ref{thm:lift-as-cocycles}.

The cleft objects $\mA$ are obtained in Theorem
\ref{thm:A-N=2}, while the algebras $\ug$ are described in Theorem \ref{thm:L(A,H)-N=2}.
\epf

\subsubsection{The main result, $N=3$}\label{subsubsec:N=3} 
Here $\xi^3 = 1$, $\xi \neq 1$. 
We fix a family of scalars 
$\bs\mu=(\mu_{(k\,l)})_{k\leq l\in\I}$
subject to the constraints and normalizations \eqref{eqn:cond-mu-gral}. Pick an extra family of scalars $\boldsymbol\lambda=(\lambda_{iij})_{i,j\in\I, |i-j|=1}$
subject to the constraints and normalizations
\begin{align}\label{eq:condiciones escalares N=3-intro}
\begin{aligned}
 \lambda_{iij}&=0 \qquad \text{ if } \chi_{iij}\neq\eps,& &\text{or } g_{iij}=1.
\end{aligned}
\end{align}
 
We define families of elements in $\mT(V)$ attached to these parameters. As in \S \ref{subsubsec:N=2}, we denote now by $(a_{i})_{i\in \I}$ the generators of $T(V)$; 
and correspondingly $a_{i j}$, $a_{(i j)}$, $a_{i_1,\dots, i_k}$. 
Let us fix $i\leq p<l\in\I$ and set $q:=p+1$, $r:=p+2$.

First, we  define 
$h_{il}(\bs\lambda)\in\Gamma$ via
\begin{align}
\label{eqn:ul-2-intro} h_{il}(\bs\lambda)&=-9\mu_{(i+2\,l)}\lambda_{ii+1i+1}\lambda_{iii+1}(1-g_{iii+1})g_{ii+1i+1}g_{(i+2\,l)}^3.
\end{align}
Next, we consider the following elements in $T(V)\# H$:
\begin{multline*}
\varsigma^{p}(\bs\lambda,\bs\mu)=\lambda_{qrr}\Big(
\xi^2a_{(i\,p)}a_{(i\,q)}a_{(i\,r)}
+\chi_{p+2}(g_{(1\,p)})
a_{(i\,p)}a_{(i\,r)}a_{(i\,q)}\\
+a_{(i\,r)}a_{(i\,p)}a_{(i\,q)}\Big).
\end{multline*}
Now, we fix $s_p=-3(1-\xi^2)$, $p<l-2$, $s_{l-2}=1$, and set
$$
d_{i\,l}(p)=\chi_{(i\,q)}(g_{(q\,l)}g_{(r+1\,l)})
\chi_{(i\,p)}(g_{(r+1\,l)}).
$$
We set 
\begin{align}\label{eqn:monomials-intro}
\varsigma_{il}(\bs\lambda,\bs\mu)&=-3\xi^2\sum_{i\leq p<l}\mu_{(p+3\,l)}
\chi_{r}(g_{(p+3\,l)})d_{il}(p)\varsigma^p(\bs\lambda,\bs\mu)g_{qrr}g_{(p+3\,l)}^3,
\end{align}
cf. Remark \ref{rem:explicit} below for a more complete description. Finally, we set 
\begin{align}\label{uil-intro}
\sigma_{(i\,l)}(\bs\lambda,\bs\mu)=h_{il}(\bs\lambda)+\varsigma_{il}(\bs\lambda,\bs\mu).
\end{align}
Let $\ug(\boldsymbol\lambda,\bs\mu)$ be the quotient of $\mT(V)$
by the relations
\begin{align}
\label{eq:rels-liftings prenichols-N3-1-intro} a_{ij} &= 0, \qquad i < j - 1; \\
\label{eq:rels-liftings prenichols-N3-2-intro} a_{iij} &= \lambda_{iij}(1-g_{iij}), \qquad \vert j - i\vert = 1;\\
\label{eq:rels-liftings-N3-intro}  a_{(i\,l)}^3&=\mu_{(i\,l)}(1-g_{(i\,l)}^3)+\uvi_{(i\,l)}+\sigma_{(i\,l)}, \qquad i\leq l\in\I.
\end{align}
for $\uvi_{(i\,l)}=\uvi_{(i\,l)}(\bs\mu)$ as in \eqref{eqn:def-gral} and
$\sigma_{(i\,l)}=\sigma_{(i\,l)}(\bs\lambda,\bs\mu)$ as in \eqref{uil-intro}.

The relations \eqref{eq:rels-liftings prenichols-N3-1-intro}
	are deformations of  \eqref{eq:rels-Atheta-a}, while 
\eqref{eq:rels-liftings prenichols-N3-2-intro} are deformations of \eqref{eq:rels-Atheta-b3} and \eqref{eq:rels-liftings-N3-intro} are deformations of \eqref{eq:rels-Atheta-c}.

\begin{rem} The quotient $\widetilde\ug(\boldsymbol\lambda)$  of $\mT(V)$ by the relations \eqref{eq:rels-liftings prenichols-N3-1-intro} and \eqref{eq:rels-liftings prenichols-N3-2-intro} is a  cocycle deformation of $\tobat(V) \# H$.
\end{rem}

\begin{theorem}\label{thm:main-N=3}
The algebra $\ug(\boldsymbol\lambda,\bs\mu)$ is a Hopf algebra quotient of $\mT(V)$ and is a lifting of $V$. 
Reciprocally every lifting of $V$ is isomorphic to 
$\ug(\boldsymbol\lambda,\bs\mu)$ for some families  $\boldsymbol\lambda$ and $\bs\mu$ as in 
\eqref{eq:condiciones escalares N=3-intro}.
In particular, every lifting is a cocycle deformation of $\toba(V) \# H$.
\end{theorem}
\pf
Similar to the proof of Theorem \ref{thm:main-N=2}, following the strategy in \S \ref{sec:strategy}. The corresponding stratification of the set of defining relations for $\toba(V)$ is given by $\Gc_0=\{\eqref{eq:rels-Atheta-a}, \eqref{eq:rels-Atheta-b3}\}$, $\Gc_1=\{\eqref{eq:rels-Atheta-c}\}$. The converse follows by Theorem \ref{thm:lift-as-cocycles}.
In this case, cleft objects $\mA = \mA(\boldsymbol\lambda,\bs\mu)$ are obtained in Theorem \ref{thm:A-N=3}, while the algebras $\ug(\boldsymbol\lambda,\bs\mu)=L(\mA,\toba(V) \# H)$ are described in Theorem \ref{thm:L(A,H)-N=3}.
\epf

\begin{rem}\label{rem:explicit}
We give an explicit description of $\varsigma_{(i\,l)}$ in terms of the PBW basis. 
See Corollary \ref{cor:suma-a}. To ease up the notation, we fix 
\begin{align*}
&j:=i+1, && k:=i+2, && q:=p+1, && r:=p+2.
 \end{align*}
Let the symmetric group $\mathbb{S}_3$ act on $\{r,q,p\}$ via $(12)(r)=q$, $(23)(q)=p$. 
If $p=i,j$, then $\varsigma^{p}_i(\bs\lambda,\bs\mu)=0$. When $p>i+2$, 
\begin{align}\label{eqn:ul-3-intro}
\begin{split}
\varsigma^{p}_{i}(\bs\lambda,\bs\mu)=-3\lambda_{qrr}&\lambda_{qqr}\chi_{(i\,p)}(g_{q})a_{(i\,p)}^3g_{qqr}\\
&-3\lambda_{qrr}\lambda_{iij}\sum_{\sigma\in\Sym_3}(-1)^{|\sigma|}h_{\sigma,i}a_{(k\,\sigma(p))}a_{(j\,\sigma(q))}a_{(i\,\sigma(r))}.
\end{split}
\end{align}
for $h_{\sigma,i}\in\k$, $\sigma\in\Sym_3$, given by:
\begin{align*}
h_{\id,i}&=\xi \chi_{qqr}(g_{(i\,p)})\chi_{(i\,r)}(g_{(j\,q)}), &h_{(12),i}&=(\xi^2-1)\chi_{qqr}(g_{(i\,p)})\chi_{i}(g_{(k\,q)}), \\
h_{(23),i}&=\xi \chi_{r}(g_{i})\chi_{i}(g_{(j\,p)}), & h_{(13),i}&=\xi (\xi-2)\chi_{(k\,p)}(g_{ij}), \\
h_{(123),i}&=2\chi_{r}(g_{(i\,p)})\chi_{i}(g_{(k\,p)}), &h_{(132),i}&=\xi^2\chi_{(k\,q)}(g_{(i\,r)})\chi_{(j\,p)}(g_{r}). 
\end{align*} 
\end{rem}

\subsubsection{The main result, $N>3$}\label{subsubsec:N>3} 
Here $\xi$ is a root of unity of order $N>3$. We fix a family of scalars 
$\bs\mu=(\mu_{(k\,l)})_{k\leq l\in\I}$
subject to the constraints and normalizations \eqref{eqn:cond-mu-gral}.

Let $\ug(\bs\mu)$ be the quotient of $\mT(V)$
by the relations
\begin{align*}
 a_{ij} &= 0, \qquad  i < j - 1; \\
 a_{iij} &= 0, \qquad \vert j - i\vert = 1;\\
 a_{(i\,l)}^N&=\mu_{(i\,l)}(1-g_{(i\,l)}^N)+\uvi_{(i\,l)}(\bs\mu), \qquad i\leq l\in\I,
\end{align*}
for $\uvi_{(i\,l)}=\uvi_{(i\,l)}(\bs\mu)$ as in \eqref{eqn:def-gral}.
\begin{theorem}\label{thm:main-N>3}
The algebra $\ug(\bs\mu)$ is a Hopf algebra quotient of $\mT(V)$ and is a lifting of $V$. 
Reciprocally every lifting of $V$ is isomorphic to 
$\ug(\bs\mu)$ for some  $\bs\mu$ as in 
\eqref{eqn:cond-mu-gral}.
Hence, every lifting is a cocycle deformation of $\toba(V) \# H$.
\end{theorem}
\pf
As in the case $N=3$, the stratification of the set of defining relations for $\toba(V)$ is given by $\Gc_0=\{\eqref{eq:rels-Atheta-a}, \eqref{eq:rels-Atheta-b3}\}$, $\Gc_1=\{\eqref{eq:rels-Atheta-c}\}$. We set $\mH=\B(V)\#H$, $\mtH=\mT(V)/\lg\Gc_0\rg$. 

In this case, the relations in $\Gc_0$ cannot be deformed cf. \cite[Theorem 5.6]{AS2}. As a result, $\Cleft'\mtH=\{\mtH\}$ and thus the corresponding deformation $\mL_1=L(\cdot,\mtH)\simeq \mtH$. This shows \eqref{eqn:condition-recursive} for $j=0$  trivially. Let us denote by $\mtA=\mtH$ as $(\mL_1,\mtH)$-bicleft object.

Pick $\bs\mu$ as in \eqref{eqn:cond-mu-gral}; set $\mA(\bs\mu)$ the quotient of $\mtH$ by the ideal generated by
\begin{align*}
 y_{(i\,l)}^N=\mu_{(i\,l)}, \qquad i\leq l.
\end{align*}
It follows from \cite{Ang} that $\,^{\co\mH}\mtH\leq\mtH$ is a normal coideal subalgebra and thus \cite[Theorem 4]{G}, see also \cite[Theorem 3.1]{AAnGMV}, yields:
\begin{align*}
\Cleft'(\mH)=\{\mA(\bs\mu)|\bs\mu\text{ as in } \eqref{eqn:cond-mu-gral}\}.
\end{align*}
In particular, \eqref{eqn:condition-recursive} holds for $j=1$. Now, we use \eqref{eqn:u-tilde} recursively, as in pp. \pageref{recursive-1} and \pageref{recursive-2}. More precisely, let $\up:\mtA\to \mL_1\ot \mtA$ denote the left coaction. Assume $i=1$ to simplify the notation and set 
\begin{align*}
A=a_{(1\,l)}\ot 1, && B=g_{(1\,l)}\ot y_{(1\,l)}, && X_p=a_{(1\,p)}g_{(p+1\,l)}\ot
y_{(p+1\,l)}, \quad
1\leq p<l,
\end{align*}
so $\up (y_{(1\,l)})=A+B+(1-\xi^{-1})\sum_{1\leq p<l} X_p$. Set also $C_p$ as in \eqref{eqn:Cp} and $r=y_{(1\,l)}^N$. By \cite[Remark 6.10]{AS2} we have:
\begin{align*}
\up (r)=a_{(1\,l)}^N\ot 1+ g_{(1\,l)}^N\ot y_{(1\,l)}^N 
+\sum_{1\leq p<l}C_{p}a_{(1\,p)}^Ng_{(p+1\,l)}^N\ot y_{(p+1\,l)}^N.
\end{align*}
We apply the deformation procedure following
\cite[Corollary 5.12]{AAnGMV}, {\it i.e.} we assume recursively $y_{(p+1\,l)}^N=\mu_{(p+1\,l)}$, and thus we get cf. \eqref{eqn:u-tilde}:
\begin{align}\label{eqn:ul-1-N>3}
\begin{split}
&\tilde r=-\sum_{1\leq p<l}C_{p}\mu_{(p+1\,l)}\Big(\uvi_{p}+\mu_{(1\,p)}(1-g_{(1\,p)}^N)\Big)g_{(p+1\,l)}^N.
\end{split}	
\end{align}
Hence, $L(\bs\mu)=L(\mA(\bs\mu),\mH)\simeq \ug(\bs\mu)$, by Proposition \ref{pro:summary} {\it (c)}.

The converse follows from Theorem \ref{thm:lift-as-cocycles}.
\epf

\subsubsection{Applications}\label{subsubsec:classif}
The classification of all finite-dimensional pointed Hopf algebras over a group algebra $H=\k G$ whose infinitesimal braiding $V$ 
is a principal realization of a braided vector space with braiding matrix \eqref{eq:braiding-A}
follows from our main results because such Hopf algebras are generated in degree one \cite{AGI}. 
When $\ord \xi > 3$, the classification was obtained in \cite[Theorem 6.25]{AS2} assuming that $G$ is abelian; the methods in \S \ref{sec:strategy} show that this hypothesis is not necessary. We extend this classification to the case in which $H$ is any cosemisimple Hopf algebra.

As a byproduct, new examples of Hopf algebras are defined, as deformations of 
intermediate  pre-Nichols algebras, see Propositions \ref{pro:lift-pre-N=2} 
and \ref{pro:lift-pre-N=3}. Also, new examples of co-Frobenius Hopf algebras arise, see \S \ref{subsubsec:cofrob} next.

\subsubsection{New examples of co-Frobenius Hopf algebras}\label{subsubsec:cofrob}

Let $\Gb$ be an algebraic group and let $H=\Oc(\Gb)$ be its function algebra; thus $\Alg(\Oc(\Gb),\k) \simeq \Gb$. 
A YD-pair for $H$, cf. \S \ref{subsec:ppal-realization}, is  $(g,x)$, where $g\in \Hom_{\text{alg gp}}(\Gb,\k^\times)$, $x\in Z(\Gb)$.

Let  $\Gb=\GL_{n}(\k)$. As
$Z(\Gb)=\k^{\times} \Id$ and $\Hom_{\text{alg gp}}(\Gb,\k^\times) =\lg\det\rg$, a YD-pair $(g,x)$ as above identifies with $(h,t)\in \Z\times \k^\times$ via $g=\det^h$, $x=t\,\Id$.

Let $V$ a braided vector space of type $A_2$, with parameter $\xi$. 
Then there is a principal YD-realization $V\in\ydh$ if and only if there are $(h_1,h_2)\in\Z^2$ and $(t_1,t_2)\in\C^2$ such that, if $u_i:=t_i^{n}$, then
\begin{align*}
&\xi =u_1^{h_1}=u_2^{h_2}; && \xi^{-1}=u_1^{h_2}u_2^{h_1}. 
\end{align*}
Each solution yields a realization $V\in\ydh$ and as a consequence of our main results Theorems \ref{thm:main-N=2}, \ref{thm:main-N=3} or \ref{thm:main-N>3}, we obtain new families of co-Frobenius Hopf algebras over $\Oc(\Gb)$. Examples of solutions are given by:
\begin{itemize}
\item $N=3$, $(u_1,u_2)=(\xi,\xi)$ and $(h_1,h_2)=(1,1)$. 
\item $N=7$, $(u_1,u_2)=(\xi,\xi^4)$ and $(h_1,h_2)=(1,2)$. 
\end{itemize} 
More examples arise considering  $\Gb= \GL_{n_1}(\k) \times \GL_{n_2}(\k) \times \dots \times \GL_{n_s}(\k)$.
\section{Preliminaries}
\subsection{Conventions}
If $n\in\N$, we set $\I_n=\{1,\dots,n\}$; we omit the subscript when it is clear from the context. We denote by $\Sym_n$ the symmetric group in $n$ letters. 
Also, $\G_n$ denotes the group of $n$th roots of 1, and $\G'_n$ is the subset of primitive  $n$th roots.

Let $H$ be a Hopf algebra; we always assume that its antipode  is bijective. We  use the Heynemann-Sweedler notation  for the comultiplication and coaction.  We 
denote by $G(H)$ the group of group-like elements of $H$ and by $\ydh$, respectively 
$\hyd$ the category of left, respectively right, Yetter-Drinfeld modules over $H$. 
If $A$ is an algebra and $S \subset A$, then $\lg S \rg$ denotes the two-sided ideal generated by $S$.

If $H'$ is a Hopf algebra, we denote by $\Isom(H,H')$ 
the set of Hopf algebra 
isomorphisms $\varphi:H\to H'$. If $A, A'$ are right $H$-comodule algebras, 
then $\Alg^H(A,A')$ is the set of comodule
algebra morphisms between them.  We shall denote by $\Alg_H^H(A,B)$ the 
set of algebra morphisms between  
two algebras $A,B\in\hyd$. When $H=\k$, we omit any reference as 
$\Alg(A,B)=\Alg_H^H(A,B)$ becomes the set of $\k$-algebra maps $A\to B$.

\subsection{Principal realizations}\label{subsec:ppal-realization}

Let $H$ be a Hopf algebra.
Let $(g, \chi)$ be a \emph{YD-pair} \cite{AAnGMV}, that is $g\in G(H)$ and $\chi\in \Alg(H, \k)$ satisfy $$\chi(h)\,g = \chi(h\_{2}) h\_{1}\, g\, \Ss(h\_{3})$$ for all $h\in H$; this implies that $g\in Z(G(H))$.
Then   $\k_g^{\chi} := \k$ with coaction given by $g$ and action given by $\chi$ is an object in $\ydh$.

Let $V$ be a braided vector space of \emph{diagonal type}, that is,
there are a basis $\left(x_{i}\right)_{i \in \I}$ of $V$ and a matrix $\q=\left(q_{ij}\right)_{i, j \in \I}$  such that $c(x_i\ot x_j)=q_{ij}x_j\ot x_i$.
A \emph{principal realization} of $V$ over $H$ is a family 
$\left((g_i, \chi_i)\right)_{i \in \I}$ of YD-pairs such that
$\chi_j(g_i)=q_{ij}$, $i, j \in \I$; so that $V\in\ydh$ 
up to identifying $\k x_i \simeq \k_{g_i}^{\chi_i}$, and the braiding $c$ is the categorical one from $\ydh$. Clearly
\begin{align}\label{eqn:Gamma}
	\Gamma=\lg g_1,\dots,g_\theta\rg \leq Z(G(H))
\end{align}
and we can realize $V$ as an object in $\ydg:=\ydkg$.

\begin{exa} There are $V \in \ydh$ with diagonal braiding but not from a
	principal realization. Let $H = \k \bH$ where $\bH$ is the finite Heisenberg group of upper triangular matrices $\left(\begin{matrix}
	1 & a & c \\ & 1 & b \\ & & 1 \end{matrix}\right)$ with coefficients
	in the finite field $\mathbb F_p$, $p$ a prime. The  conjugacy classes in 
	$\bH$ are:
	\begin{align*}
		\Oc_c &= \left\{\left(\begin{matrix}
			1 & 0 & c \\ & 1 & 0\\ & & 1 \end{matrix}\right) \right\}, &
		\Oc_{(a,b)} &= \left\{\left(\begin{matrix}
			1 & a & c \\ & 1 & b\\ & & 1 \end{matrix}\right): c\in \mathbb F_p \right\},
	\end{align*}
	for all $ c\in \mathbb F_p$, $(a,b) \in \mathbb F_p^2 - 0$. Then 
	
	\begin{itemize}[leftmargin=*]\renewcommand{\labelitemi}{$\diamond$}
		
		\item If $\rho\in \Irr \bH$, then the $M(\Oc_c, \rho) \in \ydh$ is of diagonal type, but does not arise from a principal realization unless $\dim \rho = 1$.
		
		\item If $(a,b) \in \mathbb F_p^2 - 0$ and $x \in \Oc_{(a,b)}$, then the isotropy group $\bH^x \simeq \Z_p \times \Z_p$ and $\Oc_{(a,b)}$ is an abelian rack. Hence $M(\Oc_{(a,b)}, \rho) \in \ydh$ is of diagonal type, but does not arise from a principal realization, for all $\rho\in \Irr \bH^{x}$.
	\end{itemize}
	\end{exa}

\subsection{Nichols and pre-Nichols algebras}\label{sec:prenichols}
Let $H$ and $V$ be as in \S \ref{subsec:gral-context}.
As usual, we denote by $\B(V)$ the \emph{Nichols algebra} of $V$
and by $\J(V)\subset T(V)$ its defining  ideal: $\B(V)=T(V)/\J(V)$, see \cite{AS2}. 
A \emph{pre-Nichols algebra} is a  Hopf algebra $\R=T(V)/\J\in\ydh$ with $\J\subset \J(V)$  a 
graded Hopf ideal. Every pre-Nichols algebra $\R$ is a $\zt$-graded semisimple object in $\ydh$.
The following identities are well-known. If  $x,y,z\in \R$
are $\zt$-homogeneous, then
\begin{align}\label{eq:qjacob}
[[x,y]_c,z]_c &=[x,[y,z]_c]_c+ q_{\vert y\vert \vert z\vert} [x,z]_cy-  q_{ \vert x\vert\vert y\vert} y[x,z]_c, \quad \text{(q-Jacobi)},
\\\label{eq:deriv}
\begin{split}
[xy,z]_c&=x[y,z]_c+ q_{\vert y\vert \vert z\vert} [x,z]_cy,\\
[x,yz]_c&=[x,y]_cz+ q_{ \vert x\vert\vert y\vert} y[x,z]_c.
\end{split}
\end{align}

Assume that the generalized Dynkin diagram of $V$ is connected.
The generators of the ideal $\J(V)$ were computed theoretically in \cite{A-jems} and concretely, case-by-case in the list of \cite{H-classif}, in \cite{Ang-crelle}. 
These relations can be informally organized into two types:

\begin{itemize} [leftmargin=*]\renewcommand{\labelitemi}{\tiny$\diamondsuit$}
	\item Quantum Serre relations and generalizations--that sometimes involve more than two simple roots, see \emph{e.g.} \eqref{eq:rels-Atheta-b2}.
	
	\item Powers of root vectors.
\end{itemize}

Now there are some special roots called \emph{Cartan roots} \cite[(20)]{Ang}.
There is a distinguished pre-Nichols algebra of $V$ with favourable properties,
denoted by $\tobat(V)$, cf. \cite{Ang}. The defining ideal $\mI(V)$
of $\tobat(V)$ is generated by the same relations as for $\toba(V)$, but excluding the powers of Cartan root vectors, and possibly adding some quantum Serre relations redundant for $\J(V)$.
We set 
\begin{align*}
\mT(V) = T(V)\# H, && \mH=\toba(V)\# H, && \mtH=\tobat(V)\# H,
\end{align*}
and $\pi: \mT(V)\to \mH$, $\pit:\mT(V)\to\mtH$  the natural projections.

\subsection{Cleft objects and 2-cocycles}\label{subsec:cleft}
A (normalized) Hopf 2-cocycle  is a convolution invertible linear map $\sigma:H\ot H\to \k$ such that, for  $x,y,z\in H$:
\begin{align*}
\sigma(x,1)&=\sigma(1,x)=\eps(x),\\
\sigma(x_{(1)},y_{(1)}) \sigma(x_{(2)}y_{(2)},z)&=\sigma(y_{(1)},z_{(1)}) \sigma(x,y_{(2)}z_{(2)}). 
\end{align*}
If $\sigma$ is a Hopf 2-cocycle, then it is possible to perturb the multiplication $m(x\ot y)=xy$ on $H$ on several ways, obtaining new associative products on the vector space $H$. First, we may consider $m_{(\sigma)}, m_{(\sigma^{-1})}:H\ot H\to H$ as:
\begin{align*}
m_{(\sigma)}(x\ot y)&=\sigma(x_{(1)},y_{(1)})x_{(2)}y_{(2)}, \qquad \text{respectively }\\ m_{(\sigma^{-1})}(x\ot y)&=\sigma^{-1}(x_{(2)},y_{(2)})x_{(1)}y_{(1)}.
\end{align*}
The corresponding algebras will be denoted by $H_{(\sigma)}$, respectively, $H_{(\sigma^{-1})}$. The comultiplications $\Delta:H_{(\sigma)}\to H_{(\sigma)}\ot H$,   $\Delta:H_{(\sigma^{-1})}\to H\ot H_{(\sigma^{-1})}$, remain  algebra maps and hence $H_{(\sigma)}$, respectively $H_{(\sigma^{-1})}$, is a right, resp. left, $H$-comodule algebra. Yet another associative multiplication $m_\sigma$ is defined:
\begin{align*}
m_\sigma(x\ot y)&=\sigma(x_{(1)},y_{(1)})x_{(2)}y_{(2)}\sigma^{-1}(x_{(3)},y_{(3)}), \quad x,y\in H.
\end{align*}
The corresponding algebra, denoted by $H_{\sigma}$ is actually a Hopf algebra with comultiplication $\Delta$--see  \cite{DT} for the explicit form of the antipode $\Ss_\sigma$. This Hopf algebra $H_{\sigma}$ is referred to as a {\it cocycle deformation} of $H$.

\subsubsection{Cleft objects}\label{subsubsec:cleft}

A (right) $H$-comodule algebra $A$ with trivial coinvariants, {\it i.e.} $A^{\co H}=\k$, is a cleft object of $H$ 
when there exists an $H$-colinear convolution-invertible map
$\gamma:H\to A$. This map  can be assumed to satisfy $\gamma(1)=1$, in which 
case it is called a {\it section}. Left, respectively bi-,cleft objects are defined 
accordingly. 

We shall denote by $\Cleft H$ the set of (isomorphism classes of) right cleft 
objects of $H$. If $A\in \Cleft H$, then $A$ is an algebra in $\hyd$ via
the {\it Miyashita-Ulbrich action} \cite{DT1}.

For every cleft object $A$ there is a Hopf 2-cocycle $\sigma:H\ot H\to\k$  such that
$A\simeq H_{(\sigma)}$. Indeed, a section $\gamma:H\to A$ determines $\sigma$ by
\begin{align*}
\sigma(x,y)=\gamma(x_{(1)})\gamma(y_{(1)})\gamma^{-1}(x_{(2)}y_{(2)}), \quad x,y\in H.
\end{align*}

If $A\in\Cleft H$, then there is an associated Hopf algebra $L=L(A,H)$ \cite{S} 
in such a way that $A$ becomes $(H,L)$-bicleft. Moreover, if $A=H_{(\sigma)}$, then $L\simeq H_{\sigma}$. Hence, $L$ is a cocycle deformation of $H$ and every cocycle deformation can be obtained in this way, see {\it loc.cit.} 

\section{The strategy}\label{sec:strategy}

Let $H$ be a cosemisimple Hopf algebra and $V$ as in \S \ref{subsec:gral-context}.
We recall and expand here the strategy developed in \cite{AAnGMV} to compute
the cocycle deformations of  $\B(V)\# H$. 
Accordingly, let $\Gamma$ be the abelian group as in \eqref{eqn:Gamma}.

\begin{rem}
 In \cite[\S 1.1]{AAnGMV} $H$ is assumed to be finite-dimensional. This assumption, however, can be omitted. Indeed, it
 is only used in \cite[Lemma 5.7]{AAnGMV} and  in \cite[\S 5.9, Question]{AAnGMV}. These two instances are 
 independent of the strategy and both of them deal with the evidence of an 
 ``intermediate Gunther's Theorem'' to simplify the recursive step.

On the other hand, a dimension argument is used to prove exhaustion in the 
examples, see \cite[Theorem 5.20]{AAnGMV}. We provide an alternative argument in Theorem \ref{thm:lift-as-cocycles}, valid in general.
\end{rem}

\subsection{The main idea}\label{sec:cleft}

We explain how to compute all Hopf algebras $L$ which are cocycle deformations of $\mH:=\B(V)\# H$ and
satisfy $\gr L\simeq \mH$. 
These Hopf algebras arise as $ L(\mA,\mH)$ for suitable $\mA\in\Cleft\mH$,
cf. \S \ref{subsec:cleft};
in turn, we compute the cleft extensions $\mA$ recursively by a method
from \cite{G}.

Let $\Gc$ be the  set of generators of the ideal $\J(V)$ described in \cite{Ang-crelle} for each connected component, union the $q$-commutators of vertices in different components. 
Notice that every $r \in \Gc$ belongs to $T(V)_{g_r}^{\chi_r}$ for some
$g_r\in \Gamma$, $\chi_r \in \Alg(H, \k)$.
We decompose $\Gc$
as a disjoint union $\Gc=\Gc_0\sqcup\dots \sqcup\Gc_\ell$.
Let \begin{align}\label{eqn:Bi}
\B_{i} &:=\begin{cases}
T(V), & i=0;\\
T(V)/\langle \Gc_0\cup\dots\cup
\Gc_{i-1}\rangle, & i>0;
\end{cases}& \mH_i &= \B_i\# H.
\end{align}
In particular, $\mH_{\ell+1}=\mH$.
We choose this decomposition in such a way that 
\begin{flalign}\label{eq:stratif1}
&	\text{the elements in	(the
	image of) $\Gc_{i}$, $i<\ell$, are primitive in $\B_{i}$;}
\\ \label{eq:stratif2}
&\text{$\Gc_{\ell}$ consists of powers of Cartan root vectors. }
\end{flalign}
In plain words, the strategy is to deform the relations in $\Gc$ step by step,
{\it i.e.} first those in $\Gc_0$, then those in $\Gc_1$ and so on. By \eqref{eq:stratif1},  the form of the deformed relations is particularly simple in the steps 0 to $\ell$ and depends on a suitable parameter. To check that the proposed deformation has the right properties, we proceed indirectly by defining  first a cleft extension for each possible parameter; then the proposed deformation appears as the corresponding cocycle deformation, cf. \S \ref{subsubsec:cleft}.
In the last step, the deformations of the powers of Cartan root vectors requires a delicate combinatorial analysis; but the definition of the cleft extensions is facilitated because  the algebra of coinvariants $\,^{\mH_{\ell+1}}\mH_{\ell}$ is a $q$-polynomial algebra \cite[Theorem 4.10]{Ang}. To organize the information we pack all the cleft extensions arising in the $i$-th step in a subset $\Cleft'\mH_i$ of $\Cleft\mH_i$. 

Concretely, the inductive procedure starts with 
\begin{itemize}[leftmargin=*]\renewcommand{\labelitemi}{$\circ$}
 \item the Hopf algebra $\mH_0$;
 \item the trivial $\mA_0=\mH_0\in \Cleft (\mH_0)$, where the section $\gamma:\mH_0\to \mA_0$ is
the identity map;
 \item and the corresponding Hopf algebra
$\mL_0=L(\mA_0,\mH_0)\simeq \mH_0$. 
\end{itemize}
 We now define recursively a subset $\Cleft'\mH_i$ of $\Cleft\mH_i$, 
$0\leq i\leq \ell+1$, see \cite[\S 5.2]{AAnGMV} for more details.
First, we clearly have $$\Cleft'\mH_0:=\{\mA_0\}.$$
Given $i\geq 0$, $\Cleft'\mH_{i+1}$ consists of quotients of each $\mA\in \Cleft'\mH_i$. To explain this, we fix $\mA\in \Cleft'\mH_i$; it comes equipped with 
\begin{itemize}[leftmargin=*]\renewcommand{\labelitemi}{$\diamond$}
	\item a section  $\gamma:\mH_{i}\to \mA$ such that the restriction $\gamma_{|H}:H\to \mA$ is an algebra map-- see \cite[Proposition 6.2 (b)]{AAnGMV};
	\item an algebra $\mE\in\ydh$
	such that $\mA=\mE\# H$ \cite[Proposition 5.8 (d)]{AAnGMV};
	actually $\mE$ is the image of $T(V)$ under the projection $\cA_0 = T(V) \# H \twoheadrightarrow \mA$. 
\end{itemize}

\smallbreak
Then we collect in $\Cleft'\mH_{i+1}$ all $\prova'$ given either as
\begin{align}
\label{eqn:A1}\prova' & =\cA/\cA \psi(X_i^+),& &\text{where}& X_i &:= 
\,^{\co\prov_{i+1}}\prov_i,
& \psi &\in\Alg_{\prov_i}^{\prov_i}(X_i,\prova);
\\
\notag &\text{or else as}& &&& &&
\\
\label{eqn:A2}\prova' &=\cA/\langle  \varphi(Y_i^+)\rangle,& &\text{where}& Y_i 
&:= \k\langle \Ss(\Gc_i)\rangle,& \varphi &\in \w\Alg^{\prov_i}(Y_i,\prova);
\end{align} 
here  $\w\Alg^{\prov_i}(Y_i,\prova) := \{\varphi\in \Alg^{\prov_i}(Y_i,\prova)| 
\langle \varphi (Y_i^+)\rangle \neq \cA\}$.

\begin{rem}\label{rem:psi-vphi} The subalgebra $X_i$ is the normalizer of 
$Y_i$  \cite[Remark 5.4]{AAnGMV}. If $\psi 
\in\Alg_{\prov_i}^{\prov_i}(X_i,\prova)$, then $\psi_{|Y_i}=:\varphi \in 
\w\Alg^{\prov_i}(Y_i,\prova)$ and $\langle \varphi (Y_i^+)\rangle=\cA 
\psi(X_i^+)$.
\end{rem}

\pf On one hand, $\langle \varphi (Y_i^+)\rangle \subseteq \langle \psi (X_i^+)\rangle = \cA \psi(X_i^+)$, the last equality by \cite[Theorem 4]{G}. 
Hence $\langle \varphi (Y_i^+)\rangle \neq \cA$, cf. \emph{loc.cit.}
The other inclusion follows because $X_i = N(Y_i)$.
\epf

More explicitly, given a family of scalars $\Lambda_i:= (\lambda_r)_{r\in \Gc_i}$  we  define
\begin{align}\label{eqn:E}
\mE(\Lambda_i)= \mE/\langle \gamma(r)-\lambda_r : r\in \Gc_i\rangle, &&
\mA' = \mA(\Lambda_i)=\mE(\Lambda_i) \# H.
\end{align}

Set $\mL=L(\mA,\prov_i)$. Recall that for  $r \in \Gc_i$, there are
$g_r\in \Gamma$, $\chi_r \in \Alg(H, \k)$ such that $r \in T(V)_{g_r}^{\chi_r}$.
By \cite[Corollary 5.12]{AAnGMV}, 
\begin{align}\label{eqn:u-tilde}
\nabla(r) &:=\gamma(r)_{(-1)}\otimes \gamma(r)_{(0)}
 - g_r \otimes \gamma(r) \in \mL \otimes 1, & \text{for all } 
r& \in \Gc_i.
\end{align}
Thus  $\nabla (r) = \tilde r \otimes 1$ and by {\it loc.cit.} $\tilde r$ is $(g_r,1)$-primitive in $\mL$. Set
\begin{align}\label{eqn:L-prima}
\mL' = \mL(\Lambda_i):= \mL/\langle \tilde 
r-\lambda_r(1-g_r): r\in \Gc_i\rangle.
\end{align}

The following proposition is a summary of \cite[\S 5.6]{AAnGMV}; we
add a short proof since in \emph{loc.cit.} this is stated for a single 
element in $\Gc_i$.

\begin{pro}\label{pro:summary}
Let $\mA'= \mA(\Lambda_i)$, $\Lambda_i\in \k^{\Gc_i}$.
\begin{enumerate}
\item[(a)] If $\mA'\neq 0$, then $\mA'\in \Cleft'\mH_{i+1}$.
\item[(b)] If $\chi_r\neq\eps$ and $\lambda_r\neq 0$ for some $r \in \Gc_i$, then $\mA'=0$.
\item[(c)] $L(\mA',\mH_{i+1})\simeq \mL(\Lambda_i)$.
\item[(d)] If $i= \ell$, then $\gr \mL(\Lambda_\ell)\simeq \toba(V) \# H$, {\it i.e.}  $\mL(\Lambda_\ell)$ is a lifting of $V$.
\end{enumerate}
\end{pro}
\pf
{\it (a)} Assume that $i<\ell$. Let us fix a numeration $r_{1},\dots, 
r_{s}$ of $\Gc_i$. Let $\mB_i^{(0)}:=\mB_i$, 
$\mB_i^{(t)}:=\mB_i/\langle
r_1,\dots, r_t\rangle$, $t\in \I_s$, so  $\mB_i^{(s)}=\mB_{i+1}$. By abuse of notation, the image of $r_j$ is denoted 
by $r_j$ throughout.
Set, as well, 
\begin{align*}
\mE^{(t)} &:=\mE/\langle
\gamma(r_j)-\lambda_{r_j}: j\in \I_t\rangle,& \mA^{(t)}&:=\mE^{(t)}\# 
H, & \pi^{(t)}&:\mA\to \mA^{(t)}
\end{align*}
the natural projection.  Notice that
$\mA^{(1)}\neq 0$ since it  projects onto $\mA'$ and thus $\mA^{(1)}\in\Cleft'\mB_i^{(1)}\#H$ by \cite[Remark 
5.11]{AAnGMV}. Let $\gamma^{(1)}:\mB_i^{(1)}\# H\to \mA^{(1)}$ be the
section. Observe that $\gamma^{(1)}(r_2)=\pi^{(1)}\left(\gamma(r_2)\right)$. 
Indeed, the coaction $\rho:\mA^{(1)}\to \mA^{(1)}\ot \mB_i^{(1)}\# H$ satisfies
$$
\pi^{(1)}\left(\gamma^{(1)}(r_2)\right) \overset{\rho}{\longmapsto} \pi^{(1)}(\gamma^{(1)}(r_2))\ot 1 
+ g_{r_2}\ot r_2$$
as $\pi^{(1)}$ is a comodule algebra projection that preserves $H$ cf.
the {\it
snapshot} in \cite[p. 696]{AAnGMV}. Hence we may iterate the argument and conclude that $\mA^{(t)}\in
\Cleft'\mB_i^{(t)}\# H$, $t \in \I_s$, and $\mA'=\mA^{(s)}\in
\Cleft'\mH_{i+1}$.

Next we consider the case $i=\ell$. In this step we allow the subset $\Gc_\ell$ to contain 
non-primitive elements. However, the previous analysis extends to 
this case. To see this, we decompose, in turn, 
$$\Gc_\ell=\Gc_\ell^{(0)}\sqcup\dots\sqcup \Gc_\ell^{(r)}$$
as a disjoint union of sets satisfying that $\Gc_\ell^{(0)}$ contains
primitive elements in $\mB_\ell$ and  that (the image of)  $\Gc_\ell^{(i)}$, 
$i>1$ is
composed of primitive elements in $$\mB^{(i)}_\ell=\mB_\ell/\langle 
\Gc_\ell^{(0)}\cup\dots\cup
\Gc_\ell^{(i-1)} \rangle.$$
We decompose, accordingly,
$\Lambda_\ell=\Lambda^{(0)}_\ell\times \dots\times \Lambda^{(r)}_\ell$ and 
proceed as before.

{\it (b)} follows by conjugating $\gamma(r)=\lambda_r$ by $g\in G(H)$
with $\chi_r(g)\neq 1$.

{\it (c)} follows by a iterative application of \cite[Corollary 
5.12]{AAnGMV}, as we proceed element-by-element as in {\it (a)}.

{\it (d)} For each $\mA'\in \Cleft'\mH_{i+1}$, the
section  $\gamma':\mH_{i+1}\to \mA'$ is 
such that the restriction $\gamma'_{|H}:H\to
\mA'$ is an algebra map \cite[Proposition 6.2 (b)]{AAnGMV}. Hence
$\gr L(\mA',\mH_{i+1})\simeq \mH_{i+1}$ by \cite[Proposition 4.14 
(c)]{AAnGMV}.
\epf
If $\Lambda = (\lambda_r)_{r\in \Gc}\in \k^{\Gc}$ and $0 \le i \le \ell$, then we set 
	$\Lambda_i=(\lambda_r)_{r\in \Gc_i}\in \k^{\Gc_i}$.
Set $\mA_0=T(V)\# H$ and define--using the assignment $\mA_i \rightsquigarrow \mA_{i+1}:=\mA_i(\Lambda_i)$ cf. \eqref{eqn:E}--the set of  deformation parameters
\begin{align}\label{eqn:R} \begin{split}
\mR &= \{\Lambda = (\lambda_r)_{r\in \Gc}\in \k^{\Gc} \,|\,  \mA_i(\Lambda_i)\neq 0, \ \forall i \text{ and }  \lambda_r=0 
\text{ if } g_r=1\}.
\end{split}
\end{align}
By Proposition \ref{pro:summary} we have:
\begin{cor}\label{cor:chain}
For each $\Lambda\in\mR$, we obtain a 
chain of Hopf algebra quotients
\begin{align}\label{eqn:chain}
\mL_0:=T(V)\# H\twoheadrightarrow 
\mL_1:=\mL_0(\Lambda_0)\twoheadrightarrow\dots \twoheadrightarrow
\mL_{\ell+1}:=\mL_\ell(\Lambda_\ell)
\end{align}
such that $\mL_i$ is a cocycle deformation of $\toba_i \# H$. \qed 
\end{cor}

For $\Lambda\in\mR$, we set $\mL(\Lambda):=\mL_{\ell+1}$. In 
this way, we obtain a family $\mL(\Lambda)$, $\Lambda\in\mR$, of 
cocycle deformations of $\B(V)\# H$ that are liftings of $V$.
Next we check when this family is exhaustive.
We consider the following condition on $V \in \ydh$: for $\Lambda= (\lambda_r)_{r\in \Gc} \in \k^{\Gc}$,
\begin{align}\label{eqn:condition} 
\Lambda\in \mR\ \text{ if and only if }  \lambda_r=0, \text{ when } \chi_r \neq \eps.
\end{align}
Observe that the ``only if'' implication always holds, by Proposition \ref{pro:summary} {\it (b)}.
Actually, we need a recursive version of \eqref{eqn:condition}:  

\noindent Suppose we are given  $0 \le j \le \ell$, and families 
$\Lambda_i=(\lambda_r)_{r\in \Gc_i}\in \k^{\Gc_i}$ for $i \le j$ such that $\lambda_r=0$, when $\chi_r \neq \eps$. Define recursively
$\mA_0 = \mT(V)$, $\mA_1 = \mA_0(\Lambda_0)$, $\mA_i = \mA_{i-1}(\Lambda_{i-1})$. The recursive version of \eqref{eqn:condition} is
\begin{align}\label{eqn:condition-recursive} 
\mA_j \neq 0.
\end{align}

\begin{theorem}\label{thm:lift-as-cocycles}
Assume that \eqref{eqn:condition-recursive} holds for all $j \ge 0$. If $L$ is a lifting of $V$, then there is $\Lambda\in \mR$ such that $L\simeq \mL(\Lambda)$.
In particular, $L$ is a cocycle deformation of
$\B(V)\#  H$.
\end{theorem}
\pf
Let $\phi: \mL_0 = T(V)\# H \to L$ be a lifting map. We shall attach to $\phi$ a family $(\lambda_r)_{r\in \Gc} \in \k^{\Gc}$ such that $\lambda_r=0$ if either $g_r=1$ or else
$\chi_r \neq \eps$.  Let $\simple$ be the set of simple subcoalgebras of $H$.
A direct computation shows that $V \# H \subset \sum\limits_{i\in \I, C \in \simple} g_iC \wedge C$.
Since $\phi$ is a lifting map, 
\begin{align*}
L_1 \overset{\eqref{eq:lifting-map}}{=} \phi(H \oplus V \# H) = \sum\limits_{C \in \simple}  C + \sum\limits_{i\in \I, C \in \simple} g_iC \wedge C.
\end{align*}
If $r\in\Gc_0$, then $r$ is $(g_r,1)$-primitive in $\mL_0$, hence so is $\phi(r)\in L$. That is, $\phi(r)\in \k g_r \wedge \k \subset L_1$. 
Then either $\phi(r)\in H$ or $\phi(r) \in g_iC \wedge C$ for some $i\in \I$, $C \in \simple$. In the former case, $\phi(r) = \lambda_r(1 - g_r)$ for some $\lambda_r\in \k$. 
As $g_r\in \Gamma < Z(H) \cap G(H)$, conjugation by $h\in H$ determines that  $\lambda_r = 0$ whenever $\chi_r\neq\eps$.
In the latter, $g_iC = \k g_r$ and $C = \k$, thus $g_r = g_i$ and 
\begin{equation*}
\phi(r) = \lambda_r(1 - g_r) + \sum_{j\in \I: g_j = g_r} \mu_j x_j,
\end{equation*}
for some $\lambda_r, \mu_j\in \k$. Conjugation by $h\in H$ shows that
\begin{align*}
\lambda_r\chi_r &= \lambda_r\eps,  &  \mu_j \chi_r &= \mu_j \chi_j, & \text{for } g_j &= g_r.
\end{align*}
Now, by \cite[Proposition 6.2]{AKM} the pair $(\chi_r,g_r)$ is different from $(\chi_i,g_i)$, $i \in \I$.
Thus $\mu_j = 0$ for all such $j$, hence $\phi(r) = \lambda_r(1 - g_r)$ and $\lambda_r = 0$ whenever $\chi_r\neq\eps$.
In either case, we can normalize $\lambda_r=0$ when $g_r=1$.
Set $\Lambda_0 = (\lambda_r)_{r\in \Gc_0}\in \k^{\Gc_0}$; by \eqref{eqn:condition-recursive} for $j = 0$, $\mL_1:=\mL'_0(\Lambda_0)$ is a well-defined cocycle deformation of $\mH_1$, 
and clearly  $\phi$ factorizes through $\mL_1$. 

We proceed
inductively: let $i>0$ and assume that  $\phi$ factorizes through
$\mL_{i}:=\mL'_{i-1}(\Lambda_{i-1})$, $\Lambda_{i-1}\in\k^{\Gc_{i-1}}$. Observe 
that for each $r\in \Gc_i$, the corresponding image $\tilde r\in \mL_i$ is $(g_r,1)$-primitive; cf. \eqref{eqn:u-tilde}. 
Arguing as in the previous paragraph, we conclude that
$\phi(\tilde r) = \lambda_r(1 - g_r)$ and $\lambda_r = 0$ whenever $\chi_r\neq\eps$ or $g_r=1$.
Hence there is $\Lambda_{i}$ such that $\phi$ factorizes
through $\mL_{i+1}=\mL_i'(\Lambda_{i})$, which is a well-defined cocycle deformation of $\mH_{i+1}$ by \eqref{eqn:condition-recursive} for $j = i$.
In the final step $\ell$ we proceed in the same way, splitting $\Gc_\ell$ as in the proof of Proposition \ref{pro:summary}. 
We conclude that there exists $\Lambda\in\mR$ such that
$\phi$ factorizes through $\mL(\Lambda)$.

Now, 
the lifting map $\phi$ is injective when restricted 
to 
$V\# H$ by definition, and so is the factorization 
$\phi:\mL(\Lambda)\twoheadrightarrow L$, {\it i.e.}  $\phi$ is injective 
when restricted to $\mL(\Lambda)_1$. 
Then $\phi$ is injective \cite[Theorem 5.3.1]{Mo} and thus $L\simeq 
\mL(\Lambda)$.
\epf

\subsection{Isomorphism classes}

Let $(H,V)$ be as in \S \ref{subsec:gral-context}, with braiding matrix  $\qb = (q_{ij})_{i,j\in \I}$, and $((g_i,\chi_i))_{i\in\I}$ a principal realization. 
We assume that the generalized Dynkin diagram of $V$ is connected.
Let $\Lambda\in\mR$ and $\mL(\Lambda)$ be as in Theorem \ref{thm:lift-as-cocycles}. 

\subsubsection{The block group}\label{subsec:exceptions}
Let
$$
\I(i)=\{j\in\I | g_j=g_i\text{ and }\chi_j=\chi_i\}\subseteq \I, \quad i\in \I.
$$

\begin{rem}\label{rem:exceptions} Either of the following holds:

\noindent (1) $|\I(i)|=1$ for all  $i \in \I$.

\noindent (2) There exists $i \neq j$ such that $j\in\I(i)$ is not adjacent to $i$. Then the generalized Dynkin diagram is one of the following:

\begin{enumerate} [leftmargin=*]\renewcommand{\theenumi}{\alph{enumi}}\renewcommand{\labelenumi}{(\theenumi)}
	\item Type $A_3$ with $q=-1$ \cite[Table 2, Row 1]{H-classif} and matrix 
	$ \left( \begin{smallmatrix} -1& x &-1\\ -x^{-1}& -1 &-x^{-1}\\ -1& x &-1 \end{smallmatrix} \right)$, where $x \in \k^{\times}$.
	\item \cite[Table 2, Row 8, diagrams 3 (\& 4)]{H-classif} with matrix 
	$ \left( \begin{smallmatrix} -1& x &-1\\ (qx)^{-1}& q &(qx)^{-1}\\ -1& x &-1 \end{smallmatrix} \right)$ (one diagram is obtained from the other by $q\mapsto q^{-1}$),
	 where $x \in \k^{\times}$ and $q\in \G_n$ for some $n \in \N$.
	\item \cite[Table, 2, Row 15, diagrams 2, resp. 3]{H-classif} with matrix 
	$ \left( \begin{smallmatrix} -1& x &-1\\ \xi x^{-1}& -1 &\xi x^{-1}\\ -1& x &-1 \end{smallmatrix} \right)$ resp. 
	$ \left( \begin{smallmatrix} -1& x &-1\\ \xi x^{-1}& -\xi^2 &\xi x^{-1}\\ -1& x &-1 \end{smallmatrix} \right)$, where $\xi \in \G_3'$ and  $x \in \k^{\times}$.
	\item Type $D_\theta$, $\theta\geq 4$, with $q=-1$ \cite[Table 3, Row 5 \& Table 4, Row 8]{H-classif}.
	\item \cite[Table 3, Row 18, diagrams 5 \& 6]{H-classif} (rank 4).
\end{enumerate}

\smallbreak 
\noindent (3) There exists $\xi \in \G_3'$ such that the generalized Dynkin diagram is one of the following:
\begin{enumerate} [leftmargin=*]\renewcommand{\theenumi}{\alph{enumi}}\renewcommand{\labelenumi}{(\theenumi)}
	\item\label{item-A2} Type $A_2$ with $q=\xi$, \cite[Table 1, Row 1]{H-classif}, and matrix 
	$ \left( \begin{smallmatrix} \xi&\xi\\ \xi&\xi \end{smallmatrix} \right)$.
	\item \cite[Table 2, Row 15, diagram 4]{H-classif} and matrix 
	$ \left( \begin{smallmatrix} \xi& \xi & x\\ \xi & \xi & x\\ \xi^2 x^{-1}& \xi^2 x^{-1} &-1 \end{smallmatrix} \right)$, where $x \in \k^{\times}$.
\end{enumerate}
\end{rem}

\pf First, observe that if $j\in\I(i)$, then every vertex not adjacent to $i$ cannot be adjacent to $j$ as
$1=\chi_k(g_i)\chi_i(g_k)=\chi_k(g_j)\chi_j(g_k)$.
Next, if $j\in\I(i)$, $j \neq i$,  is not adjacent to $i$, then $\chi_j(g_j)=\chi_i(g_i)=-1$ as
$1=\chi_j(g_i)\chi_i(g_j)=\chi_j(g_j)^2$.
Then (2) and (3) follow by inspection in \cite{H-classif}. \epf

Let us denote by 
\begin{align*} 
\Lb:=\{s\in \GL_\theta(\k)|s_{ij}=0 \text{ if } j\notin\I(i)\}.
\end{align*}
Observe that $\Lb\simeq \{(s_i)_{i\in\I}\in\k^{\times\theta}\}$ if 
 the generalized Dynkin diagram is as in Remark \ref{rem:exceptions} (1).
If the diagram is as in (3)(a), then $\Lb=\GL_{2}(\k)$.

\subsubsection{Isomorphisms}\label{subsec:isom}
We fix a new pair $(H',V')$ as in \S \ref{subsec:gral-context}. Set $\theta'=\dim V'$, $\I'=\I_{\theta'}$. Fix a principal realization $((g'_i,\chi'_i))_{i\in\I'}$ of $V'$ in
$\ydhp$ and let $\Gamma'=\lg g_i'\mid i\in\I'\rg\leq H'$ be as in \eqref{eqn:Gamma}. 

Let $\Gc'$ be 
the set of generators of the ideal defining $\toba(V')$ and 
$\mR'\subseteq\k^{\Gc'}$ as in 
\eqref{eqn:R}. Pick $\Lambda'\in\mR'$ and consider the Hopf algebra 
$\mL(\Lambda')$. 
 Let
$$
\Sym_{\qb} =\{\sigma\in\Sym_\theta| q_{ij} = q_{\sigma(i) \sigma(j)} \,\forall i, j \in \I \}.
$$
 
\begin{lem}\label{lem:iso-coradical}
Let $\psi:\mL(\Lambda)\to \mL(\Lambda')$ be a Hopf algebra isomorphism.
 Then
$\varphi=\psi_{|H}:H\to H'$ is a Hopf algebra isomorphism and $T=\psi_{|V}:V\to V'$ is an isomorphism of braided vector spaces. In particular $\theta=\theta'$. Moreover,
\begin{itemize}
\item[(i)] There is $\sigma\in\Sym_{\qb}$  such that $\varphi(g_{j})=g'_{\sigma(j)}$ and $\chi'_{\sigma(j)}\circ \varphi=\chi_j$, $j\in\I$.
\item[(ii)] There is $s = (s_{ij})\in\Lb$ such that $T(a_i)=\sum_{j\in\I(\sigma(i))} s_{ij}a'_j$, $i\in\I$.
\end{itemize}
\end{lem}
\pf
Follows since the map $\psi$ preserves both the comultiplication and the coradical filtration, as well as the adjoint action.
\epf

\begin{rem}\label{rem:subsets}
When $|\I(i)|=1$, $i\in\I$, Lemma \ref{lem:iso-coradical} (ii) reads
\begin{itemize}
\item[{\it (ii')}] There are scalars $\{s_{i}\}_{i\in\I}$ such that $
T(a_i)=s_{i}a'_{\sigma(i)}$. 
\end{itemize}
\end{rem}

Assume that $\theta=\theta'$, $H'\simeq H$. We fix $\varphi\in\Isom(H,H')$, $\sigma\in\Sym_\theta$ and  
$s\in\Lb$. 
We say that a triple $(\varphi,\sigma,s):(H,V,\Lambda)\to(H',V',\Lambda')$ is a 
{\it lifting data isomorphism} if 
\begin{itemize}
\item $\sigma\in\Sym_{\qb}$. 
\item $g'_{i}=\varphi(g_{\sigma(i)})$ and $\chi'_{i}=\chi_{\sigma(i)}\circ\varphi$, $i\in \I$.
\item $\Lambda'=s\cdot\Lambda^\sigma$, cf. Lemmas \ref{lem:s-Lambda} and \ref{lem:sigma-Lambda}.
\end{itemize}
Set 
$\Isom(\Lambda,\Lambda')=
\{\text{lifting data isomorphisms}:(H,V,\Lambda)\to(H',V',\Lambda')\}$.

\begin{theorem}\label{teo:iso}
 $\Isom(\mL(\Lambda),\mL(\Lambda'))\simeq \Isom (\Lambda,\Lambda')$.
\end{theorem}
\pf
By Lemma 
\ref{lem:iso-coradical}, any $\psi\in \Isom(\mL(\Lambda),\mL(\Lambda'))$ univocally determines a triple 
$(\varphi,\sigma,s)\in \Isom (\Lambda,\Lambda')$.

Conversely, let $(\varphi,\sigma,s)\in \Isom (\Lambda,\Lambda')$. In particular, $\mL(\Lambda')$ is an $H$-module via $\varphi$. Consider the linear map $T=T_s^\sigma:V\to \mL(\Lambda')$ given by $T^\sigma_s(a_i)=\sum_{j\in\I(\sigma(i))} s_{ij}a'_j$, $i\in\I$. By assumption, $T$ is $H$-linear and hence it defines an algebra epimorphism $F:T(V)\#H\to \mL(\Lambda')$ with $F_{|H}=\varphi$ and $F(a_i)=T(a_i)$, $i\in\I$. 
By a combination of Lemmas \ref{lem:s-Lambda} and \ref{lem:sigma-Lambda}, the map $F$ induces an isomorphism
 $\widetilde F\in \Isom(\mL(\Lambda),\mL(\Lambda'))$.
The assignment $(\varphi,\sigma,s)\mapsto \widetilde F$ is injective, as 
each triple determines a Hopf algebra map in the first term of the 
coradical filtration, hence in the whole algebra.

These constructions are inverse to each other and define a bijective correspondence  $\Isom(\mL(\Lambda),\mL(\Lambda'))\simeq \Isom (\Lambda,\Lambda')$.
\epf

We set $\mH_i=\toba_i(V)\# H$, $i\geq 0$, see \eqref{eqn:Bi}. If $\Lambda\in\mR$, then we set, cf. \eqref{eqn:E}: $\mA_0(\Lambda):=\mT(V)\in\Cleft(\mH_0)$, $\mA_{i+1}(\Lambda):=\mA_i(\Lambda_i)\in\Cleft(\mH_{i+1})$ . Let 
$$\rho_i:\mA_i\to\mA_i\ot\mH_i, \qquad \gamma_i:\mH_i\to \mA_i.$$
denote the coaction and section. Also we set $\mL_i(\Lambda):=\mL_{i}$ as in \eqref{eqn:chain}. 

\begin{lem}\label{lem:s-Lambda}
There is a well-defined action $\Lb\times\mR\to 
\mR$ so that if $s\in\Lb$, $\Lambda\in\mR$, then $\mL_i(s\cdot \Lambda)\simeq \mL_i(\Lambda)$ as Hopf algebras.
\end{lem}
\pf
We fix $\Lambda\in\mR$, $s\in\Lb$.
We shall assume for simplicity that each stratum $\Gc_i$ of $\Gc$ cf. \eqref{eqn:Bi} contains {\it all} primitive elements of $\toba_i(V)$. The general case follows analogously.

We define $s\cdot \Lambda\in\mR$. That is, we define for each $i\geq 0$ a family of scalars $s\cdot \Lambda_i\in\k^{\Gc_i}$ such that the algebras 
defined recursively as $\mA^{(s)}_{0}=\mA_{0}$ and $\mA^{(s)}_{i+1}=\mA^{(s)}_{i}(s\cdot \Lambda_i)$, cf. \eqref{eqn:E}, are nonzero. Hence $s\cdot \Lambda:=(s\cdot \Lambda_i)_{i\geq 0}\in\mR$. Moreover, we show that 
$\mA^{(s)}_{i}(s\cdot \Lambda_i)\simeq \mA_i(\Lambda_i)$ as cleft objects, all $i$. As a result, $\mL_i(s\cdot \Lambda)\simeq \mL_i(\Lambda)$ as Hopf algebras.

Let $V_s$ be the vector space with basis $\{F_s(x_k)\}_{k\in\I}$. Then $V_s$ is braided, with the braiding from $V$ by assumption on $s\in\Lb$. Set $\mH_{s\cdot i}=\B_i(V_s)\#H$. Let $F_0:\mH_0\to\mH_{s\cdot 0}$ be the unique algebra automorphism with
\begin{align*}
F_{0|H}=\id \quad \text{ and } \quad F_0(x_k)=F_s(x_k), \qquad k\in\I.
\end{align*}
By assumption, $F_0(\k\Gc_0)=\k\Gc_0\subset T(V_s)$ and thus it induces an algebra automorphism $F_1:\mH_1\to\mH_{s\cdot 1}$. Similarly, $F_1(\k\Gc_1)=\k\Gc_1$ and, in general, there is an induced automorphism $F_i:\mH_i\to\mH_{s\cdot i}$, $i\geq 0$.

\begin{claim}
There is $\mA_{s\cdot i}\in\Cleft(\mH_{s\cdot i})$ together with an  algebra automorphism $f_i:\mA_i\to  \mA_{s\cdot i}$ such that 
\begin{align}\label{eqn:claim-iso}
&\rho_{s\cdot i}\circ f_i=(f_i\ot F_i)\circ \rho_i, \qquad f_i\circ \gamma_i(r)=\gamma_{s\cdot i}(F_i(r)), \ r\in\Gc_i.
\end{align}
\end{claim}
This is clear when $i=0$, for $f_0=F_0$, $\mA_{s\cdot 0}:=\mH_{s\cdot 0}$, 
$\rho_{s\cdot 0}=\Delta$, $\gamma_{s\cdot 0}=\id$.

Assume that, for a given $i\geq 0$, we have defined $\mA_{s\cdot i}$ so that \eqref{eqn:claim-iso} holds. If $r\in\Gc_i$, then $x=\gamma_{s\cdot i}(F_s(r))\in \mA_{s\cdot i}$ is unique such that
$$
\rho_{s\cdot i}(x)=x\ot 1+g_r\ot F_s(r)\in\mA_{s\cdot i}\ot\mH_{s\cdot i}.$$ This is satisfied by $x=f_i\circ \gamma_i(r)$ and hence 
$f_i$ descends as to an isomorphism 
\begin{align*}
f_{i+1}:\mA_{i+1}\to \mA_{s\cdot (i+1)}:=\mA_{s\cdot i}/\lg \gamma_{s\cdot i}(F_i(r))-\lambda_r : r\in\Gc_i\rg,
\end{align*}
and \eqref{eqn:claim-iso} defines a structure $\mA_{s\cdot (i+1)}\in\Cleft(\mH_{s\cdot (i+1)})$.

Now, the composition $\mA_0\twoheadrightarrow \mA_i(\Lambda_i)\overset{f_i}\to \mA_{s\cdot i}$ defines a family of scalars $s\cdot \Lambda_i\in\k^{\Gc_i}$ with $\mA^{(s)}_{i}(s\cdot \Lambda_i)\simeq \mA_i(\Lambda_i)$, $i\geq 0$. Hence $s\cdot \Lambda\in\mR$.
\epf

We consider the action of $\Sym_{\qb}$ on $T(V)$ by permutations of
the generators. If $\Lambda=(\lambda_r)_{r\in\Gc}\in\mR$, then we set 
$\Lambda^\sigma:=(\lambda_{\sigma\cdot r})_{r\in\Gc}\in\k^\Gc$, $\sigma\in\Sym_{\qb}$.

\begin{lem}\label{lem:sigma-Lambda}
There is a well-defined action $\Sym_{\qb}\times\mR\to 
\mR$ so that if $\sigma\in\Sym_{\qb}$, $\Lambda\in\mR$, then $\mL_i(\Lambda^\sigma)\simeq \mL_i(\Lambda)$ as Hopf algebras.
\end{lem}
\pf
Proceed as in Lemma \ref{lem:s-Lambda}, {\it mutatis mutandis}.
\epf

We give  examples of the action $\Lb\times\mR\to 
\mR$ from Lemma \ref{lem:s-Lambda}.

\begin{exa}\label{rem:action}
(1) Assume $|\I(i)|=1$, $i\in \I$; hence $\Lb\simeq \k^{\times\theta}$. If $s=(s_i)_{i\in\I}\in\Lb$
and $r\in T(V)$ is a 
$\Z^\theta$-homogeneous element with $\deg r=(d_1,\dots,d_\theta)$, then we set $s_r:=s_1^{d_1}\cdots s_\theta^{d_\theta}\in\k^\times$. 
If $\Lambda=(\lambda_r)_{r\in\Gc}\in\mR$ and $s\in\Lb$, then 
$s\cdot \Lambda:=(s_r\lambda_r)_{r\in\Gc}$.

\smallbreak
(2) Assume $V$ is as in Remark \ref{rem:exceptions} (3)(a), so $\Lb=\GL_{2}(\k)$. In this case $\Lambda=(\lambda_{112},\lambda_{122},\mu_1,\mu_2,\mu_{12})\in\k^5$ by Theorem \ref{thm:main-N=3}. Let $s= \left( \begin{smallmatrix} s_{11}&s_{12}\\ s_{21}&s_{22} \end{smallmatrix} \right)\in\Lb$ and denote $s\cdot \Lambda:=(\lambda_{112}^s,\lambda_{122}^s,\mu_1^s,\mu_2^s,\mu_{12}^s)$. Then 
\begin{align*}
\mu_1^s&=s_{11}^3\mu_1 + s_{12}^3\mu_2 + s_{11}^2s_{12}\lambda_{112} 
+ s_{11}s_{12}^2\lambda_{122},\\
\mu_2^s&=s_{21}^3\mu_{1} + s_{22}^3\mu_{2} + s_{21}^2s_{22}\lambda_{112} + s_{21}s_{22}^2\lambda_{122},  \\
\lambda_{112}^s&=3s_{11}^2s_{21}\mu_{1} + 3s_{12}^2s_{22}\mu_{2} + 
(s_{11}^2s_{22}+2s_{11}s_{12}s_{21})\lambda_{112}\\
& \hspace*{4cm}+ (2s_{11}s_{12}s_{22}+
s_{12}^2s_{21})\lambda_{122}, \\
\lambda_{122}^s&=3s_{11}s_{21}^2\mu_{1} + 3s_{12}s_{22}^2\mu_2+ 
(2s_{11}s_{21}s_{22}+s_{12}s_{21}^2)\lambda_{112} \\
& \hspace*{4cm}+ (s_{11}s_{22}^2+2s_{12}s_{21}s_{22})\lambda_{122},\\
\mu_{12}^s&=(s_{11}s_{22}-s_{12}s_{21})^3\mu_{12}.
\end{align*}
\end{exa}

\subsection{The algorithm}\label{sec:algorithm}
Our strategy reduces the lifting problem to an algorithm, that we describe next.

Let $H$, $V$ be as in \S \ref{subsec:gral-context}, $\Gamma$ as in \eqref{eqn:Gamma}. Let $\Gc$ be the  set of generators of the ideal $\J(V)$ defining $\B(V)$ as described in \cite{Ang-crelle}  for each connected component, union the $q$-commutators of vertices in different components. Decompose it as $\Gc=\Gc_0\sqcup\dots \sqcup\Gc_\ell$ so that \eqref{eq:stratif1} and \eqref{eq:stratif2} hold.

The algorithm involves $\ell+1$ recursive steps.
At each Step {\bf i}, the {\it input} are two Hopf algebras $\mH_i$ and $\mL_i$, a ($\mL_i,\mH_i$)-bicleft object $\mA_i$, with coactions $\rho_i,\up_i$ and a choice of scalars $\Lambda_i=(\lambda_r)_{r\in \Gc_i}\in \k^{\Gc_i}$ such that
\begin{align*}
\begin{aligned}
\lambda_r&=0,&  &\mbox{ if }\chi_r \neq \eps,&
&\mbox{ or }g_r=1. 
 \end{aligned}
\end{align*}
The {\it output} is a new triple $(\mH_{i+1},\mA_{i+1},\mL_{i+1})$, as quotient of the input data. Step {\bf 0} starts with $\mH_0=\mL_0=T(V)\# H$ and $\mA_0=\mH_0$, with $\rho_0=\delta_0=\Delta$.

The {\it final outcome} of the algorithm is a list of liftings of $V$ in terms of families $\Lambda\in\k^{\Gc}$. All of them are cocycle deformations of $\B(V)\#H$.  If no step produces a zero object, then this list is exhaustive.

The recursive step is the following:
\begin{stepi} \quad
\begin{enumerate}
 \item Compute $r'\in\mA_i$, $r\in\Gc_i$. These elements are defined by the equation: 
\begin{align*}
 \rho_i(r')=r'\ot 1+g_r\ot r, \quad r\in\Gc_i.
\end{align*}
 \item Set 
 $
 \mA_{i+1}:=\mA_i(\Lambda_i)=\mA_i/\langle  r'-\lambda_r : r\in \Gc_i  \rangle
 $ and check $\mA_{i+1}\neq 0$.
 \item Compute $\tilde r\in\mL_i$, $r\in\Gc_i$. These elements are defined by the equation:
$$
\up_i(r')=\tilde r\ot 1+g_r\ot r', \quad r\in\Gc_i.
$$ 
\item Set $\mH_{i+1}:=\mH_i/\lg \Gc_i\rg$,
  $\mL_{i+1}:=\mL_i/\langle \tilde r-\lambda_r(1-g_r) : r\in \Gc_i \rangle$.
\end{enumerate}
\end{stepi}

\begin{rem}
We make some comments regarding the recursive step.

1. At Step {\bf 0}, $r'=r$, for each $r\in\Gc_0$.

2. At Step $\bs\ell$, $\mA_{\ell+1}\neq 0$ automatically.

3. At Step {\bf i}, $1\leq i\leq\ell$, the verification of (2) is facilitated by the fact that $\mA_{i+1}=\mE_{i+1}\# H$, for 
$\mE_{i+1}\in \ydh$, $i\geq 0$, the algebra defined recursively by
\begin{align*}
\mE_0=T(V), \qquad 
\mE_{i+1}=\mE_i/\langle  r'-\lambda_r : r\in \Gc_i  \rangle.
\end{align*}
\end{rem}

\section{The case $N=2$}\label{sec:N=2}

Let $H$, $V$ as in \S \ref{subsec:gral-context}, $\Gamma$ as in \eqref{eqn:Gamma}. Assume moreover that $V$ is of type $A_\theta$, $\theta\in\N$, associated to $\xi=-1$. Let $\B(V)$ be the corresponding 
Nichols algebra.
In this section we compute the liftings of $V$. We show that all of them arise as cocycle deformations of $\B(V)\# H$.

Recall the definition of the distinguished pre-Nichols algebra $\tobat(V)$, see Proposition \ref{prop:presentation-nichols} (2). Set $\mtH=\tobat(V)\# H$.

\begin{lem}\label{lem:rels preNichols} Let $i\le j\le k\le l$.
The following relations holds in $\mtH$:
\begin{align}
\label{eqn:lem4.1-1} [x_{(i\,j)}, x_{(i\,k)}]_c&= 0,  & [x_{(i\,k)}, x_{(j\,k)}]_c&=0,   \\
\label{eqn:lem4.1-2} [x_{(i\,l)}, x_{(j\,k)}]_c&=0,  & [x_{(i\,k)}, x_{(j\,l)}]_c&= 2\chi_{(j\,k)}(g_{(i\,k)}) x_{(j\,k)} x_{(i\,l)}.
\end{align}
The coproduct of $\mtH$ satisfies
\begin{align*}
\Delta(x_{(i\,j)})&=x_{(i\,j)}\ot 1+ g_{(i\,j)}\ot x_{(i\,j)}+ 2 \sum_{k=i}^{j-1} x_{(i\,k)}g_{(k+1\,j)}\ot x_{(k+1\,j)}, \\
\Delta(x_{(i\,j)}^2)&=x_{(i\,j)}^2\ot 1+ g_{(i\,j)}^2\ot x_{(i\,j)}^2 \\
&\qquad + 4 \sum_{k=i}^{j-1} \chi_{(i\,k)}(g_{(k+1\,j)}) x_{(i\,k)}^2g_{(k+1\,j)}^2\ot x_{(k+1\,j)}^2.
\end{align*}
\end{lem} 
\pf
It follows as in \cite[Section 6]{AS2}, see also \cite[Section 3]{AD}, by induction. The key point to show \eqref{eqn:lem4.1-2} is to use  \eqref{eq:rels-Atheta-b2} as the initial step. On the other hand, relations \eqref{eqn:lem4.1-1} follow from \eqref{eq:rels-Atheta-a} and \eqref{eq:rels-Atheta-b3}. The formula for the coproduct now follows. 
\epf

In the remaining part of this section, we deal with a quotient of $\tobat(V)$, namely we fix the 
the algebra $\tobah(V)$ generated by $x_1, \dots, x_{\theta}$ with
relations
\begin{align}\label{eq:rels-Atheta-N2-prenichols}
x_{ij} &= 0, \, i< j - 1; & [x_{(i-1i+1)},x_i]_c &= 0, \, 2\leq i<\theta; & x_{k}^2 &=0, \, k\in\I_\theta.
\end{align}

This algebra is an intermediate quotient between $\tobat(V)$ and $\toba(V)$, see Proposition \ref{prop:presentation-nichols} and Remark \ref{rem:distinguished}. We prefer the quotient \eqref{eq:rels-Atheta-N2-prenichols} as it is more suitable for our computations. We set $\mhH=\tobah(V)\# H$. Observe that Lemma \ref{lem:rels preNichols} holds for $\mhH$.

Recall also that the Nichols algebra $\toba(V)$ is generated by $x_1, \dots, x_{\theta}$ with the previous defining
relations and also $x_{(i\,j)}^2=0$ for $i<j$. We set $\mH=\toba(V)\# H$. Let $\pi:\mhH\twoheadrightarrow\mH$ be
the canonical Hopf algebra map. Recall that $\mhH^{\co \pi}$ is the subalgebra generated by $x_{(i\,j)}^2$, $i<j$,
which is a polynomial algebra with these elements as generators.

\subsection{Cleft objects}\label{subsec:cleft N2}

Let $\boldsymbol\lambda=(\lambda_{ij})_{1\leq i<j-1< \theta}$, $\boldsymbol\mu=(\mu_{(k\, l)})_{1\le k\le l\le\theta}$, 
$\boldsymbol\nu=(\nu_{i})_{1<i<\theta}$ be families of scalars such that
\begin{align}\label{eq:condiciones escalares N=2}
\begin{split}
\lambda_{ij}=0 & \mbox{ if }\chi_i\chi_j\neq\eps, \qquad  \mu_{(k\, l)}=0 \mbox{ if }\chi_{(kl)}^2\neq\eps,
\\ \nu_{i}=0 & \mbox{ if }\chi_i^2\chi_{i-1}\chi_{i+1}\neq\eps.
\end{split}
\end{align}
Let us set, following Proposition \ref{pro:summary}, $\mhA=\mhA(\boldsymbol\lambda)$
the quotient of $T(V)\# H$ by the
relations
\begin{align}\label{eq:rels-Atheta-N2-galois}
\begin{split}
y_{ij} &=\lambda_{ij}, \ i < j - 1; \qquad [y_{(i-1i+1)},y_i]_c = \nu_{i}, \ 2\leq i<\theta; \\
y_{k}^2&=\mu_{(k)}, \ 1\leq k\leq \theta.
\end{split}
\end{align}
Here, we have renamed the basis $\{x_1,\dots,x_\theta\}$ of $V$ by $\{y_1,\dots,y_\theta\}$.

\begin{pro}\label{pro:cleft-tilde-N2}
The algebras $\mhA(\boldsymbol\lambda,\bs\mu,\bs\nu)$ are cleft objects for $\mhH$. Hence
\begin{align*}
\Cleft'\mhH=\{\mhA(\boldsymbol\lambda,\boldsymbol\mu,\boldsymbol\nu)|\boldsymbol\lambda,\boldsymbol\mu,\boldsymbol\nu\text{ as in }\eqref{eq:rels-Atheta-N2-galois}\}.
\end{align*}
\end{pro}
In particular, this shows that \eqref{eqn:condition-recursive} holds for $j=0$.
\pf
Set $\mhA=\mhA(\boldsymbol\lambda)$ and $\mhE$ the quotient of $T(V)$ by the ideal $I$ generated by \eqref{eq:rels-Atheta-N2-galois}.
Observe that $\mhA\simeq \mhE\# H$, as $I$ is an object in $\ydh$.
Hence we need to show that $\mhE\neq0$.

For this we use Diamond Lemma \cite[Theorem 1.2]{B - diamond}. We introduce a notation close to the one in \emph{loc. cit.}
Let $\Xi_{ij}=(w_{ij},f_{ij})$ be the pair associated to the relation $y_{ij}-\lambda_{ij}$; we choose $w_{ij}=y_iy_j$, so
$f_{ij}=q_{ij}y_jy_i+\lambda_{ij}$. Similarly we set $\Xi_i=(y_i^2,\mu_{(i)})$ for $1\le i\le\theta$, and $\Xi_i'=(y_{i-1}y_iy_{i+1}y_i,f_i')$,
$2\le i\le\theta-1$, for the relation $[y_{(i-1i+1)},y_i]_c- \nu_{i}$, where for $i=2$,
\begin{align*}
f_2'&=q_{12}^2q_{13} y_2y_3y_2y_1-q_{12}q_{13}q_{23}y_3y_2y_1y_2-q_{12}q_{32} y_2y_1y_2y_3 \\
& \qquad +2q_{23}\lambda_{13}\mu_{(2)}+\nu_2.
\end{align*}
There are no \emph{inclusion ambiguities}. There are eight \emph{overlap ambiguities}:
\begin{enumerate}
  \item $(\Xi_{ij},\Xi_{jk},y_i,y_j,y_k)$. Both $y_if_{jk}$ and $f_{ij}y_k$ reduce to
$$ q_{ij}q_{ik}q_{jk} y_ky_jy_i+\lambda_{ij}y_k+q_{jk}\lambda_{ik}y_j+\lambda_{jk}y_i, $$
since $\chi_k(g_ig_j)\lambda_{ij}=\chi_{ij}^{-1}(g_k)\lambda_{ij}=\lambda_{ij}$.
  \item $(\Xi_{ij},\Xi_{j},y_i,y_j,y_j)$. As $q_{ij}^2\mu_j=\mu_j$ and $\lambda_{ij}(1+q_{ij})=0$, both $y_if_{j}$
and $f_{ij}y_j$ reduce to $\mu_jy_i$.
  \item $(\Xi_{i},\Xi_{ij},y_i,y_i,y_j)$. Analogous to the previous case.
  \item $(\Xi_{i},\Xi_{i+1}',y_i,y_i,y_{i+1}y_{i+2}y_{i+1})$. For simplicity set $i=1$. To prove that $y_1f_2'$ reduces
to $\mu_{(1)} y_2y_3y_2$ we use the identities $[y_1,y_{123}]_c=0$ obtained from $\Xi_1$, and $[y_{12},y_{123}]_c=0$, which is obtained from
$\Xi_2'$ and the previous relation.
  \item $(\Xi_{i+1}',\Xi_{i+1},y_iy_{i+1}y_{i+2},y_{i+1},y_{i+1})$. Again set $i=1$. Then $f_2'y_2$ reduces to $\mu_{(2)} y_1y_2y_3$ up to reduce by 
$\Xi_2'$.
  \item $(\Xi_{i+1}',\Xi_{i+1\,j},y_iy_{i+1}y_{i+2},y_{i+1},y_j)$. Again set $i=1$. If $j>4$, then both $f_2'y_j$ and $y_1y_2y_3f_{2j}$ reduce to
\begin{align*}
& q_{1j}q_{2j}^2q_{3j} \Big( q_{12}^2q_{13} y_jy_2y_3y_2y_1-q_{12}q_{13}q_{23}y_jy_3y_2y_1y_2-q_{12}q_{32} y_jy_2y_1y_2y_3
\\ & \quad +2q_{23}\lambda_{13}\mu_{(2)}y_j+\nu_2 y_j \Big) + \lambda_{1j}q_{2j}^2q_{3j} y_2y_3y_2 + \lambda_{2j}\lambda_{13} q_{2j}q_{3j} y_2
\\ & \quad + \lambda_{2j}q_{13}q_{3j} y_3y_1y_2 + \lambda_{3j}\mu_{(2)} q_{2j}y_1+\lambda_{2j} y_1y_2y_3
\end{align*}
by direct computation. If $j=4$, then use the relation $[y_{(14)},y_2]_c=0$ obtained from $f_2'$ and $f_2$ to reduce the word
$y_1y_2y_3y_4y_2$ and obtain the same reduction for both $f_2'y_4$ and $y_1y_2y_3f_{24}$.
  \item $(\Xi_{ij},\Xi_{j+1}',y_i,y_j,y_{j+1}y_{j+2}y_{j+1})$. Fix $j=1$, $j=3$ to simplify the notation. By direct computation we reduce both 
expressions to
\begin{align*}
\lambda_{13} &y_4y_5y_4 +\lambda_{14} q_{13}q_{35} y_5 y_3 y_4  +\lambda_{14} q_{13}q_{14}q_{15} y_3y_4y_5+\lambda_{15}\mu_{(4)} q_{13}q_{14} y_3 \\
&+q_{13}q_{14}^2q_{15} \Big( q_{34}^4q_{35} y_4y_5y_4y_3-q_{34}q_{35}q_{45}y_5y_4y_3y_4 -q_{34}q_{54} y_4y_3y_4y_5\Big)y_1 \\
&+\lambda_{14} \lambda_{35} q_{13} y_4+ 4q_{45}\lambda_{35}\mu_{(4)}y_1+\nu_4y_1.
\end{align*}
  \item $(\Xi_{i+1}',\Xi_{i+2}',y_iy_{i+1}y_{i+2},y_{i+1},y_{i+2}y_{i+3}y_{i+2})$. Again assume that $i=1$. Both $y_1y_2y_3 f_3'$ and
$f_2'y_3y_4y_3$ reduce to
\begin{align*}
q_{12}q_{43} & \lambda_{13}\mu_{(2)}\mu_{(3)} y_4 -2q_{12}q_{23}q_{24} \lambda_{14}\mu_{(2)}\mu_{(3)} y_3
+q_{23}\lambda_{13}\mu_{(2)} y_3y_4y_3  \\
&  +\nu_3 y_1y_2y_3 +\nu_2 y_3y_4y_3
-q_{12}q_{32}q_{13} \mu_{(3)}\lambda_{24}y_2y_3y_1
+q_{12}^2q_{13}^2 \lambda_{14}\mu_{(3)} y_2y_3y_2 \\
&+q_{12}q_{13}q_{23}^3q_{24}^2 \lambda_{13} y_3y_4y_2y_3y_2
-q_{12}q_{32}q_{14}q_{24}^2 \mu_{(3)} y_4y_2y_1y_2y_3 \\
&+q_{12}q_{13}^2q_{14}q_{23}^3q_{24} \lambda_{13} y_3y_2y_3y_4y_2
-q_{12}q_{13}q_{23}^2q_{24}q_{34} \lambda_{14} y_3y_2y_3y_2y_3 \\
&+q_{12}^2q_{13}q_{24} \lambda_{13} y_2y_3y_4y_2y_1
+q_{12}^2q_{13}^2q_{14} \lambda_{13} y_2y_3y_2y_3y_4 \\
&+q_{12}q_{13}^3q_{14}q_{23}^3q_{24} y_3y_2y_3y_4y_3y_1y_2
-q_{12}q_{13}^2q_{14}q_{23}^2q_{24}^2q_{34} y_3y_4y_2y_3y_1y_2y_3 \\
&+q_{34}\lambda_{24}\mu_{(3)} y_1y_2y_3-q_{12}q_{13}^2q_{23}^2q_{43} y_3y_2y_3y_1y_2y_3y_4.
\end{align*}
Note that $\lambda_{13}\lambda_{24}=0$ since $\chi_{13}\chi_{24}(g_{(1\,4)})=-1$, so either $\chi_{13}\neq\eps$ or else $\chi_{24}\neq\eps$.
\end{enumerate}
The proposition now follows from Proposition \ref{pro:summary}.
\epf

\begin{lem}
For all $j<k$,
$$ \rho( y_{(j\,k)})= y_{(j\,k)}\ot 1+ g_{(j\,k)}\ot x_{(j\,k)}+2\sum\limits_{l=j+1}^{k-1}  y_{(j\,l)}g_{(l+1\,k)}\ot x_{(l+1\,k)}. $$
\end{lem}
\pf
By induction on $k-j$. If $k=j+1$, then
$$ \rho( y_{(j\,j+1)})= y_{(j\,j+1)}\ot 1+ g_{(j\,j+1)}\ot x_{(j\,j+1)}+2  y_{j}g_{j+1}\ot x_{j+1}. $$
by direct computation. If it holds for $k-j$, then
\begin{align*}
\rho & (y_{(j-1\,k)})  = \rho( y_{j-1})\rho( y_{(j\,k)})-\chi_{(j\,k)}(g_{j-1})\rho( y_{(j\,k)})\rho( y_{j-1})\\
& = y_{(j-1\,k)}\ot 1+ \big(1-\chi_{(j\,k)}(g_{j-1})\chi_{j-1}(g_{(j\,k)})\big) \, y_{j-1}g_{(j\,k)}\ot x_{(j\,k)} \\
& +  2\sum_{l=j+1}^{k-1} \big( y_{j-1} y_{(j\,l)} -\chi_{(j\,k)}(g_{j-1})\chi_{j-1}(g_{(l+1\,k)}) y_{(j\,l)} y_{j-1} \big)
g_{(l+1\,k)}\ot x_{(l+1\,k)}  \\
&  + g_{(j\,k)}\ot x_{(j\,k)} + 2\sum\limits_{l=j+1}^{k-1} \chi_{(j\,l)}(g_{j-1}) y_{(j\,l)}g_{j-1}g_{(l+1\,k)}\ot [x_{j-1},x_{(l+1\,k)}]_c.
\end{align*}
Notice that $1-\chi_{(j\,k)}(g_{j-1})\chi_{j-1}(g_{(j\,k)})=2$. For each $j+1\le l\le k-1$,
$$ y_{j-1} y_{(j\,l)} -\chi_{(j\,k)}(g_{j-1})\chi_{j-1}(g_{(l+1\,k)}) y_{(j\,l)} y_{j-1}=y_{(j-1\,l)}, $$
$[x_{j-1},x_{(l+1\,k)}]_c=0$ by Lemma \ref{lem:rels preNichols}, and the inductive step follows.
\epf

\begin{lem}\label{lema:ytilde qconmutan}
For all $j\le k<l$, $ y_{(j\,k)} y_{(j\,l)}=\chi_{(j\,l)}(g_{(j\,k)}) y_{(j\,l)} y_{(j\,k)}$.
\end{lem}
\pf
Set $\mathtt{y}_{j,k,l}=y_{(j\,k)} y_{(j\,l)}-\chi_{(j\,l)}(g_{(j\,k)}) y_{(j\,l)} y_{(j\,k)}$. If $k=j$, $l=j+1$, then
\begin{align*}
\mathtt{y}_{j,j,j+1}&=y_jy_{(j\,j+1)}-q_{jj}q_{(j\,j+1)}y_{(j\,j+1)}y_j \\
&=(\mu_j y_{j+1}-q_{(j\,j+1)}y_jy_{j+1}y_j)+q_{(j\,j+1)}(-q_{(j\,j+1)}\mu_j y_{j+1}+y_jy_{j+1}y_j)\\
&=(1-\mu_jq_{(j\,j+1)}^2)y_{j+1}=(1-\mu_jq_{(j\,j+1)}^2\chi_j^2(g_{j+1}))y_{j+1}=0.
\end{align*}
Assume it holds for all $k',l'$ such that $k'+l'<k+l$. Then
\begin{align*}
\rho (\mathtt{y}_{j,k,l})&= \mathtt{y}_{j,k,l}\ot 1 +\mathrm{(i)} 
+\mathrm{(ii)}+\mathrm{(iii)}+\mathrm{(iv)}+\mathrm{(v)}+\mathrm{(vi)}+\mathrm{(vii)}.
\end{align*}
We compute now the other seven summands. We use repeatedly Lemma \ref{lem:rels preNichols} and inductive hypothesis.
\smallskip

\noindent$\mathrm{(i)} = (1-\chi_{(j\,k)}(g_{(j\,l)})\chi_{(j\,l)}(g_{(j\,k)}))y_{(j\,k)}g_{(j\,l)}\ot x_{(j\,l)}= 2 y_{(j\,k)}g_{(j\,l)}\ot 
x_{(j\,l)}$.

\noindent$\mathrm{(ii)} = 2\sum\limits_{t=1}^{l-1}\big(y_{(j\,k)}y_{(j\,t)}- \chi_{(j\,l)}(g_{(j\,k)})\chi_{(j\,k)}(g_{(t+1\,l)}) 
y_{(j\,t)}y_{(j\,k)}\big)g_{(t+1\,l)}\ot x_{(t+1\,l)}$. If $t>k$, then 
$\chi_{(j\,l)}(g_{(j\,k)})\chi_{(j\,k)}(g_{(t+1\,l)})=\chi_{(j\,t)}(g_{(j\,k)})$ so the summand is zero by Lemma \ref{lem:rels preNichols}. For $t=k$,
$\chi_{(j\,l)}(g_{(j\,k)})\chi_{(j\,k)}(g_{(k+1\,l)})=1$, so the summand is also zero. If $t<k$, then $y_{(j\,t)}y_{(j\,k)}=\chi_{(j\,k)}(g_{(j\,t)}) 
y_{(j\,k)}y_{(j\,t)}$, so we have that

\noindent$\mathrm{(ii)} = 4 \sum\limits_{t=1}^{k-1}y_{(j\,k)}y_{(j\,t)} g_{(t+1\,l)}\ot x_{(t+1\,l)}$. 

\noindent$\mathrm{(iii)} = g_{(j\,k)}g_{(j\,l)}\ot [x_{(j\,k)},x_{(j\,l)}]_c=0$.

\begin{align*}
\mathrm{(iv)} &= 2\sum\limits_{t=1}^{l-1}\chi_{(j\,t)}(g_{(j\,k)}) y_{(j\,t)}g_{(j\,k)}g_{(t+1\,l)}\ot [x_{(j\,k)},x_{(t+1\,l)}]_c \\
&= -4 \sum\limits_{t=1}^{k-1}\chi_{(s+1\,t)}(g_{(k+1\,s)}) y_{(j\,t)}g_{(j\,k)}g_{(t+1\,l)}\ot x_{(j\,l)}x_{(t+1\,k)} \\
&\qquad -2 y_{(j\,k)}g_{(j\,l)}\ot x_{(j\,l)}. \\
\mathrm{(v)} & = 2\sum\limits_{s=1}^{k-1}\chi_{(j\,l)}(g_{(s+1\,k)}) [y_{(j\,s)},y_{(j\,l)}]_c g_{(s+1\,k)}\ot x_{(s+1\,k)}=0.\\
\mathrm{(vi)} &= 2\sum\limits_{s=1}^{k-1} y_{(j\,s)}g_{(s+1\,k)}g_{(j\,l)} \\
& \qquad \qquad \ot \big(x_{(j+1\,s)}x_{(j\,l)} -\chi_{(j\,l)}(g_{(j\,k)})\chi_{(j\,s)}(g_{(j\,l)}) x_{(j\,l)}x_{(s+1\,k)}\big) \\
&= 4 \sum\limits_{t=1}^{k-1}\chi_{(k+1\,l)}(g_{(t+1\,k)}) y_{(j\,t)}g_{(j\,k)}g_{(t+1\,l)}\ot x_{(j\,l)}x_{(t+1\,k)} = - \mathrm{(iv)}-\mathrm{(i)}.\\
\mathrm{(vii)} &= 4 \sum\limits_{t=1}^{l-1} \sum\limits_{s=1}^{k-1}  \chi_{(j\,t)}(g_{(s+1\,k)})y_{(j\,s)}y_{(j\,t)} g_{(s+1\,k)}g_{(t+1\,l)} \ot 
x_{(s+1\,k)}x_{(t+1\,l)}  \\
& \qquad - \chi_{(j\,l)}(g_{(j\,k)})\chi_{(j\,s)}(g_{(t+1\,l)}) y_{(j\,t)} y_{(j\,s)} g_{(t+1\,l)} g_{(s+1\,k)} \ot x_{(t+1\,l)} x_{(s+1\,k)}.
\end{align*}
For $\mathrm{(vii)}$ there are three subcases:
\begin{itemize}[leftmargin=*]
  \item If $t>k$, then $y_{(j\,s)}y_{(j\,t)} = \chi_{(j\,t)}(g_{(j\,s)})y_{(j\,t)}y_{(j\,s)}$ and $x_{(s+1\,k)}x_{(t+1\,l)} = 
\chi_{(t+1\,l)}(g_{(s+1\,k)})x_{(t+1\,l)}x_{(s+1\,k)}$. Hence these summands are 0.
  \item If $t=k$, $x_{(s+1\,k)}x_{(k+1\,l)} = x_{(s+1\,l)}+ \chi_{(k+1\,l)}(g_{(s+1\,k)})x_{(k+1\,l)}x_{(s+1\,k)}$, and $y_{(j\,s)}y_{(j\,k)} = 
\chi_{(j\,k)}(g_{(j\,s)})y_{(j\,k)}y_{(j\,s)}$, so the summand is $-\mathrm{(ii)}$.
  \item For $t<k$, the summands cancel between themselves.
\end{itemize}
Thus $\rho (\mathtt{y}_{j,k,l})= \mathtt{y}_{j,k,l}\ot 1$, so $\mathtt{y}_{j,k,l}\in\k$. Also, 
$\mathtt{y}_{j,k,l}\in\mhA_{\chi_{(j\,k)}\chi_{(j\,l)}}$. As $\chi_{(j\,k)}\chi_{(j\,l)}(g_{(j\,k)}g_{(j\,l)})=-1$, we have that 
$\mhA_{\chi_{(j\,k)}\chi_{(j\,l)}}\cap\k=0$ so $\mathtt{y}_{j,k,l}=0$.
\epf

\begin{lem}\label{lema:ytilde cuadrado}
For all $j<k$,
$$ \rho( y_{(j\,k)}^2)= y_{(j\,k)}^2\ot 1+ g_{(j\,k)}^2\ot x_{(j\,k)}^2+4\sum_{s=j+1}^{k-1}\chi_{(j\,s)}(g_{(s+1\,k)})  y_{(j\,s)}^2g_{(s+1\,k)}^2\ot 
x_{(s+1\,k)}^2.$$
\end{lem}
\pf As $\rho$ is an algebra map,
\begin{align*}
\rho( y_{(j\,k)}^2)&=\Big( y_{(j\,k)}\ot 1+ g_{(j\,k)}\ot x_{(j\,k)}+2\sum_{s=j+1}^{k-1} y_{(j\,s)}g_{(s+1\,k)}\ot x_{(s+1\,k)}\Big)^2.
\end{align*}
By Lemmas \ref{lem:rels preNichols} and \ref{lema:ytilde qconmutan} all the summands $q$-commute.
\epf

\begin{lem}\label{lema:ytilde cuadrado q conmutan}
For all $j<k$ and all $i$, $y_{(j\,k)}^2 y_i = \chi_i(g_{(j\,k)}^2)y_i y_{(j\,k)}^2 $.
\end{lem}
\pf
By induction on $k-j$. If $k=j+1$, then
$$ \rho(y_{(j\,j+1)}^2 y_i - \chi_i(g_{(j\,j+1)}^2)y_i y_{(j\,j+1)}^2)= (y_{(j\,j+1)}^2 y_i - \chi_i(g_{(j\,j+1)}^2) y_i y_{(j\,j+1)}^2) \ot 1 $$
since $x_{(j\,j+1)}^2 x_i = \chi_i(g_{(j\,j+1)}^2) x_i x_{(j\,j+1)}^2$. But $y_{(j\,j+1)}^2 y_i - \chi_i(g_{(j\,j+1)}^2)y_i y_{(j\,j+1)}^2
\in \mhA_{\chi_i\chi_{(j\,j+1)}^2}$ and $\chi_i\chi_{(j\,j+1)}^2(g_ig_{(j\,j+1)}^2)=-1$, so $y_{(j\,j+1)}^2 y_i = \chi_i(g_{(j\,j+1)}^2) 
y_iy_{(j\,j+1)}^2$.

A similar proof follows for the inductive step since for all $j<k$ and all $i$, $x_{(j\,k)}^2 x_i = \chi_i(g_{(j\,k)}^2) x_ix_{(j\,k)}^2$,
see \cite[Proposition 4.1]{Ang}.
\epf

The following theorem shows that \eqref{eqn:condition-recursive} holds for $j=1$, hence \eqref{eqn:condition} holds in general.
\begin{theorem}\label{thm:A-N=2}
Let $\mA=\mA(\boldsymbol\lambda,\boldsymbol\mu,\boldsymbol\nu)$ be the quotient of $\mhA$ by the relations
\begin{align}\label{eq:rels PRV-N2-galois}
y_{(i\,j)}^2&=\mu_{(i\,j)}, \ 1\leq i<j \leq \theta.
\end{align}
Then $\mA\in\Cleft\mH$. As a result,
\begin{align*}
\Cleft'\mH=\{\mA(\boldsymbol\lambda,\boldsymbol\mu,\boldsymbol\nu)|\boldsymbol\lambda,\boldsymbol\mu,\boldsymbol\nu\text{ as in }\eqref{eq:rels-Atheta-N2-galois}\}.
\end{align*}
\end{theorem}
\pf
Indeed these algebras are obtained following \cite[Theorem 4]{G}. As in \emph{loc. cit.} we need to describe the $\mhH$-linear and colinear algebra 
maps $\,^{\co \pi}\mhH\to \mhA$. As $^{\co \pi}\mhH$ is a polynomial ring in the variables $x_{(i\,j)}^2g_{(j\,k)}^{-2}$, it is enough to determine 
the value on $x_{(i\,j)}^2g_{(j\,k)}^{-2}$.
Set $f(x_{(i\,j)}^2g_{(j\,k)}^{-2})=y_{(i\,j)}^2g_{(j\,k)}^{-2}-\mu_{(i\,j)}g_{(j\,k)}^{-2}$. Then $f$ is $\mhH$-colinear by Lemmas \ref{lem:rels 
preNichols} and \ref{lema:ytilde cuadrado}. We claim that $f$ is also $\mhH$-linear. Indeed, for all $g\in H$ and all $1\le k\le \theta$,
\begin{align*}
f(g\cdot x_{(i\,j)}^2g_{(j\,k)}^{-2})&= \chi_{i\,j}^2(g) f(x_{(i\,j)}^2g_{(j\,k)}^{-2}) = g\cdot f(x_{(i\,j)}^2g_{(j\,k)}^{-2}), \\
f(x_k \cdot x_{(i\,j)}^2g_{(j\,k)}^{-2})&= 0 = x_k \cdot f(x_{(i\,j)}^2g_{(j\,k)}^{-2}),
\end{align*}
where the first equality holds by \eqref{eq:condiciones escalares N=2} and the second
by Lemma \ref{lema:ytilde cuadrado q conmutan} and \cite[Proposition 4.1]{Ang}. The claim follows since $\mhH$ is generated by $H$ and the 
$x_k$'s as an algebra. Then $\mhA/\mhA f((^{\co \pi}\mhH)^+)=\mhA/\mhA f((^{\co \pi}\mhH)^+)\mhA=\mA(\boldsymbol\lambda,\bs\mu,\bs\nu)$ is a cleft 
object of $\mH$ by Proposition \ref{pro:summary}.
\epf

\subsection{Liftings}\label{subsec:liftings}

In this subsection we give a presentation for the Hopf algebras
$L(\mhA(\boldsymbol\lambda,\bs\mu,\bs\nu),\mhH)$
and
$L(\mA(\boldsymbol\lambda,\bs\mu,\bs\nu),\mH)$. We apply Proposition 
\ref{pro:summary} together with formula 
\eqref{eqn:u-tilde}.

\begin{pro}\label{pro:lift-pre-N=2}
The Hopf algebra 
$\mhL(\boldsymbol\lambda,\bs\mu,\bs\nu)=L(\mhA(\boldsymbol\lambda,\bs\mu,\bs\nu)
,\mhH)$ is the quotient 
of $\mT(V)$ by relations
relations
\begin{align}
a_{ij} &=\lambda_{ij}(1-g_ig_j); \label{eq:rels-liftings prenichols-N2-1}\\
a_{k}^2&=\mu_{(k)}(1-g_k^2); \label{eq:rels-liftings prenichols-N2-2} \\
[a_{(i-1\,i+1)},a_i]_c &= \nu_{i}(1-g_i^2g_{i-1}g_{i+1}) \label{eq:rels-liftings prenichols-N2-3} \\
& \qquad -4\chi_i(g_{i-1})\mu_{(i)}\lambda_{i-1\,i+1}g_{i-1}g_{i+1}(1-g_i^2).\notag
\end{align}
In particular, it is a cocycle deformation of $\mhH$ with 
$\gr\mhL(\boldsymbol\lambda,\bs\mu,\bs\nu) \simeq \mhH$.
\end{pro}
\pf
We follow Proposition \ref{pro:summary} (c): The $x_{ij}$'s and the $x_{k}^2$'s 
are skew-primitive elements in $T(V)\# H$,
so we quotient $T(V)\#H$ by relations \eqref{eq:rels-liftings prenichols-N2-1} and \eqref{eq:rels-liftings prenichols-N2-2} to obtain the 
corresponding lifting. Again, Proposition \ref{pro:summary} (c), see also 
\cite[Corollary 5.12]{AAnGMV}, applies for the relation $[x_{(i-1i+1)},x_i]_c$ 
since it is primitive,
and $\tilde u=[a_{(i-1\,i+1)},a_i]_c 
+4\chi_i(g_{i-1})\mu_{(i)}\lambda_{i-1\,i+1}g_{i-1}g_{i+1}(1-g_i^2)$ is the 
corresponding skew-primitive element, see \eqref{eqn:u-tilde}.
\epf

Let $i\neq j \in \I$. If $|i-j|\ge2$, then we define recursively scalars $d_{i\,j}(s)$, $b_{i\,j}(s)$, $s\ge 0$, as:
$d_{i\,j}(0)=2\lambda_{i\,j}$, $b_{i\,j}(0)=-2\chi_j(g_{(i\,j)})\lambda_{i\,j}$, and for $s>0$,
\begin{align}\label{eq:corchete yi yjk - formula coef-intro}
d_{i\,j}(s)&= q_{ij} \sum_{l=0}^{s-1} d_{i\,j+1}(l)d_{j\,j+2l+2}(s-l-1),\\
b_{i\,j}(s)&= \sum_{l=0}^{s-1} b_{i+1\,j}(l)d_{i\,i+2l+2}(s-l-1).
\end{align}
If $|i-j|=1$, then we set $d_{i\,j}(s)=b_{i\,j}(s)=0$, for $s\ge0$.
In what follows $y_{(k+1\,k)}:=1$, to simplify the summation formulas.

\begin{rem}\label{rem:aij=0}
Notice that $d_{i\,j}(s)=0$ if $\chi_i\chi_{(j\,j+2s)}\neq\eps$.
\end{rem}

\begin{lem}\label{lem:corchete yi yjk}
Let $j<k$, $i\notin\{j-1,j,\dots,k+1\}$. Then
\begin{align}\label{eq:corchete yi yjk}
[y_i, y_{(j\,k)}]_c &= \sum_{s=0}^{ \frac{k-j}{2}} d_{i\,j}(s) \,  y_{(j+2s+1\,k)}.
\end{align}
\end{lem}

\pf
By induction on $j-k$. If $k=j+1$, then
\begin{align*}
[y_i&,  y_{(j\,j+1)}]_c = \lambda_{i\,j}(1-\chi_{j+1}(g_ig_j))y_{j+1}+\lambda_{i\,j+1}(\chi_{j}(g_i)-\chi_{j+1}(g_j))y_j \\
&= \lambda_{i\,j}(1-\chi_{j+1}(g_ig_j)\chi_i\chi_j(g_{j+1}))y_{j+1}+\lambda_{i\,j+1}(\chi_{j}(g_i)-\chi_i^{-1}(g_j))y_j \\
&= 2\lambda_{i\,j}y_{j+1}.
\end{align*}
The inductive step follows from the following formula:
\begin{align*}
[y_i, & y_{(j-1\,k)}]_c = 2\lambda_{i\,j-1} y_{(j\,k)}+ q_{i\,j-1} \big[ y_{j-1}, [y_i, y_{(j\,k)}]_c \big]_c \\
& = d_{i\,j-1}(0) y_{(j\,k)}+ q_{i\,j-1}\sum_{s=0}^{ \frac{k-j}{2}} d_{i\,j}(s) \,  \big[ y_{j-1},   y_{(j+2s+1\,k)} \big]_c \\
& = d_{i\,j-1}(0) y_{(j\,k)}+ q_{i\,j-1} \hspace{-3pt} \sum_{s=0}^{ \frac{k-j}{2}} d_{i\,j}(s) \hspace{-6pt}  \sum_{t=0}^{ \frac{k-j-1-2s}{2}} 
\hspace{-6pt} d_{j-1\,j+2s+1}(t) \,  y_{(j+2(s+t+1)\,k)}.
\end{align*}
Here we have applied the inductive hypothesis twice.
\epf

\begin{lem}\label{lem:corchetes yk+1, yk contra yjk}
Let $j<k$. Then
\begin{align}\label{eq:corchete yk+1 yjk}
[y_{(j\,k)}, y_{k+1}]_c &= y_{(j\,k+1)}-\sum_{s=1}^{ \frac{k-j}{2}} b_{j\,k+1}(s) \,  y_{(j+2s+1\,k)}, \\
\label{eq:corchete yk yjk}
[y_{(j\,k)}, y_k]_c &= \sum\limits_{s=0}^{ \frac{k-j-1}{2}} b_{j\,k}(s) \,  y_{(j+2s+1\,k)}.
\end{align}
\end{lem}

\pf
First we prove \eqref{eq:corchete yk+1 yjk} by induction on $k-j$. For $k=j+1$ we have
\begin{align*}
[y_{(j\,j+1)}, & y_{j+2} ]_c = [[y_{j},y_{j+1}]_c,y_{j+2}]_c \\
&= y_{(j\,j+2)}-\chi_{j+1}(g_{j}) y_{j+1}[y_{j},y_{j+2}]_c+\chi_{j+2}(g_{j+1}) [y_{j},y_{j+2}]_c y_{j+1} \\
&= y_{(j\,j+2)}+\chi_{j+2}(g_{j+1})(1+\chi_{j+1}(g_{j\,j+2}))\lambda_{jj+2} y_{j+1}  \\
&= y_{(j\,j+2)}-b_{jj+2}(0) y_{j+1}.
\end{align*}
Now assume it holds for $j',k'$ such that $k'-j'<k-j$. Then by inductive hypothesis, Lemma \ref{lem:corchete yi yjk} and using  $\chi_{k+1}^{-1}=\chi_j$ if $\lambda_{jk+1}\neq 0$ we obtain:

\begin{align*}
[y_{(j\,k)},&y_{k+1}]_c = [[y_{j},y_{(j+1\,k)}]_c,y_{k+1}]_c \\
&= [y_{j},[y_{(j+1\,k)},y_{k+1}]_c]_c +\big(\chi_{k+1}(g_{(j+1\,k)})- \chi_{(j+1\,k)}(g_j)\big)\lambda_{jk+1}y_{(j+1\,k)} \\
&= [y_{j}, y_{(j+1\,k+1)}-\sum_{s=1}^{ \frac{k-j}{2}} b_{j+1k+1}(s) \,  y_{(j+2s+1\,k+1)} ]_c \\
& \qquad + \chi_{k+1}(g_{(j+1\,k)})\big(1-\chi_j(g_{(j+1\,k)}) \chi_{(j+1\,k)}(g_j)\big)\lambda_{jk+1}y_{(j+1\,k)}  \\
&= y_{(j\,k+1)} -\sum_{s=1}^{ \frac{k-j}{2}} b_{j+1k+1}(s)  \sum_{t=0}^{ \frac{k-2s-j}{2}}
d_{jj+2s+2t+2}\,  y_{(j+2s+2t+2\,k+1)} \\
& \qquad- b_{jk+1}(0)y_{(j+1\,k)},
\end{align*}

Now we prove \eqref{eq:corchete yk yjk}. For $k=j+1,j+2$ we have
\begin{align*}
[y_{(j\,j+1)},y_{j+1}]_c &= \mu_{(j+1)}(1-q_{(j\,j+1)}^2)y_j=0,\\
[y_{jj+2},y_{j+2}]_c &= [[y_{j},y_{j+1j+2}]_c,y_{k+2}]_c \\
&= \lambda_{jj+2}(\chi_{j+2}(g_{j+1j+2})-\chi_{j+1j+2}(g_{j})) y_{j+1j+2} \\
&=-2\chi_{j+2}(g_{j+1}) \lambda_{jj+2} y_{j+1j+2} = b_{jj+2}(0)y_{j+1j+2}.
\end{align*}
Then we argue by induction in $k-j$ as for \eqref{eq:corchete yk+1 yjk}.
\epf

We define recursively $\tz_{(j\,k)}\in\mhL$ as follows: $\tz_{(j\,j)}=a_j$ and for $j<k$
\begin{multline}\label{eq:def ztilde}
\tz_{(j\,k)} = [a_{j},\tz_{(j+1\,k)}]_c+ d_{jk}(0)\chi_{(j\, k)}(g_j)\tz_{(j+1\,k-1)} g_{jk}\\ + 2 \sum_{t=1}^{ \frac{k-j-1}{2} } d_{j k-2t}(t) 
\chi_{(j+1\, k-2t-1)}(g_{j}) \tz_{(j+1\,k-2t-1)} g_{j}g_{(k-2t\,k)}.
\end{multline}

\begin{rem}\label{rem:a lo sumo un aij}
If $s\neq t$, then $d_{j k-2t}(t)d_{j k-2s}(s)=0$  by Remark \ref{rem:aij=0}.
Thus there is at most one non-trivial summand besides $[a_{j},\tz_{(j+1\,k)}]_c$ in \eqref{eq:def ztilde}.
\end{rem}

\begin{lem}
The $\mhL$-coaction $\up$ of $\mhA$ satisfies:
\begin{align*}
\up(y_{(j\,k)})= \tz_{(j\,k)}\ot 1+ g_{(j\,k)}\ot y_{(j\,k)} + 2 \sum_{s=j}^{k-1} \tz_{(j\,s)}g_{(s+1\,k)}\ot y_{(s+1\,k)}.
\end{align*}
\end{lem}
\pf
Again by induction: the case $k=j+1$ is direct. Now assume it holds for $k-j$. Then we compute
\begin{align*}
\up(&y_{(j-1\,k)}) = \up(y_{j-1})\up(y_{(j\,k)})- \chi_{(j\,k)}(g_{j-1})\up(y_{(j\,k)})\up(y_{j-1})
\\
& = [a_{j-1},\tz_{(j\,k)}]_c \ot 1+ 2 a_{j-1} g_{(j\,k)}\ot y_{(j\,k)} \\
& + 2 \sum_{s=j}^{k-1} [a_{j-1},\tz_{(j\,s)}]_c g_{(s+1\,k)}\ot y_{(s+1\,k)} + g_{(j-1\,k)}\ot y_{(j-1\,k)} \\
& + 2 \sum_{s=j}^{k-1} \chi_{(j\,s)}(g_{j-1}) \tz_{(j\,s)}g_{j-1}g_{(s+1\,k)}\ot [y_{j-1},y_{(s+1\,k)}]_c \\
& = [a_{j-1},\tz_{(j\,k)}]_c \ot 1+ 2 a_{j-1} g_{(j\,k)}\ot y_{(j\,k)} \\
& + 2 \sum_{s=j}^{k-1} [a_{j-1},\tz_{(j\,s)}]_c g_{(s+1\,k)}\ot y_{(s+1\,k)} + g_{(j-1\,k)}\ot y_{(j-1\,k)} \\
& + 2 \sum_{s=j}^{k-2} \chi_{(j\,s)}(g_{j-1}) \tz_{(j\,s)} g_{j-1} g_{(s+1\,k)}\ot \sum_{t=0}^{ \frac{k-s-1}{2}} d_{j-1s+1}(t) \,  y_{(s+2t+2\,k)} \\
& + 2 \chi_{(j\, k-1)}(g_{j-1}) \tz_{(j\,k-1)} g_{j-1k}\ot \lambda_{j-1k},
\end{align*}
by Lemma \ref{lem:corchete yi yjk}. The proof follows by reordering the summands.
\epf

Now, for each $m\geq 1$, consider the $m$-adic approximation $\widehat{\B}_m(V)$. This is the quotient of $T(V)$ by relations
\eqref{eq:rels-Atheta-N2-prenichols} and
\begin{align*}
 & x_{(k\,l)}^2, \qquad 1\leq l-k< m.
\end{align*}
Thus, we obtain a family of cleft objects $\mA_m(\boldsymbol\lambda,\bs\mu,\bs\nu)$ for
$\mH_m=\widehat{\B}_m(V)\#H$ given by the  quotient of $\mT(V)$
by relations \eqref{eq:rels-Atheta-N2-galois} together with
\begin{align}\label{recursive-1}
 & y_{(k\,l)}^2-\mu_{(k\,l)}, \qquad 1\leq l-k< m.
\end{align}
Let $\mL_m(\boldsymbol\lambda,\bs\mu,\bs\nu):=L(\mA_m,\mH_m)$. Notice that $\mL_0=\mhL$.
We keep the name $\up:\mA_m\to \mL_m\ot \mA_m$ for the coaction at each level.

For the next two lemmas we consider a fixed $m$.

\begin{lem}\label{lem:corchete yik contra yjk}
Let $i<j<k$ be such that $k-j<m$. Then there exist $c_{ijk}(s,t)\in\k$ such that
\begin{align}\label{eq:corchete yik contra yjk}
[y_{(i\,k)}, y_{(j\,k)}]_c &= \sum_{i<s<t\le k+1} c_{ijk}(s,t) \,  y_{(t\,k)}y_{(s\,k)}.
\end{align}
\end{lem}
\pf
By induction on $j-i$. If $j=i+1$, then
\begin{align*}
[y_{(j-1\,k)}, y_{(j\,k)}]_c &= \big(y_{j-1} y_{(j\,k)} - \chi_{(j\,k)}(g_{j-1}) y_{(j\,k)}y_{j-1} \big) y_{(j\,k)} \\
& \quad - \chi_{(j\,k)}(g_{(j-1\,k)}) y_{(j\,k)} \big(y_{j-1} y_{(j\,k)} - \chi_{(j\,k)}(g_{j-1}) y_{(j\,k)}y_{j-1} \big) \\
&= \mu_{(j\,k)}(1-\chi_{(j\,k)}^2(g_{j-1})) y_{j-1}=0,
\end{align*}
since $\chi_{(j\,k)}(g_{(j-1\,k)})=\chi_{(j\,k)}(g_{j-1})\chi_{(j\,k)}(g_{(j\,k)})=-\chi_{(j\,k)}(g_{j-1})$.

Now assume it holds for all $i'<j'$ such that $j-i>j'-i'$. Then
\begin{multline*}
[y_{(i-1\,k)}, y_{(j\,k)}]_c = [y_{i-1}, [y_{(i\,k)}, y_{(j\,k)}]_c ]_c- \chi_{(i\,k)}(g_{i-1}) y_{(i\,k)} [y_{i-1},y_{(j\,k)}]_c \\
 + \chi_{(j\,k)}(g_{(i\,k)})  [y_{i-1},y_{(j\,k)}]_c y_{(i\,k)} = \sum_{i<s<t\le k+1} c_{ijk}(s,t) [y_{i-1}, y_{(t\,k)}y_{(s\,k)}]_c \\
+ \sum_{r=0}^{\frac{k-j}{2}} d_{i-1j}(s)  \big( 2 \chi_{(j\,k)}(g_{(i\,k)})y_{(j+2s+1\,k)}y_{(i\,k)}- \chi_{(j\,k)}(g_{i-1}) 
[y_{(i\,k)},y_{(j+2s+1\,k)}]_c\big).
\end{multline*}
We apply the inductive step to express $[y_{(i\,k)},y_{(j+2s+1\,k)}]_c$ as a linear combination of products $y_{(t\,k)}y_{(s\,k)}$. Also,
$$ [y_{i-1}, y_{(t\,k)}y_{(s\,k)}]_c = [y_{i-1}, y_{(t\,k)}]_c y_{(s\,k)} + \chi_{(t\,k)}(g_{i-1})y_{(t\,k)}[y_{i-1}, y_{(s\,k)}]_c . $$
We apply Lemma \ref{lem:corchete yi yjk} and the inductive step to obtain a linear combination of elements $y_{(t\,k)}y_{(s\,k)}$
for $k+1\ge t\ge s>i$. Assume that some $y_{(t\,k)}^2$ appears with non-zero coefficient. Then $\chi_{(i\,k)}\chi_{(j\,k)}=\chi_{(t\,k)}^2$, so 
$\chi_{(i\,t-1)}\chi_{(j\,t-1)}=\eps$, which contradicts
$\chi_{(i\,t+1)}\chi_{(j\,t+1)}(g_{(i\,t+1)}g_{(j\,t+1)})=-1$.
\epf

\begin{lem}
Let $j\le k$. There exist $\mathtt{z}_{(j\,k)}(s,t)\in \mhL$ such that
\begin{align}\notag
\up(y_{(j\,k)}^2)&= \tz_{(j\,k)}^2\ot 1 + 4 \sum_{s=j}^{k-1} \chi_{(j\,s)}(g_{(s+1\,k)})\tz_{(j\,s)}^2g_{(s+1\,k)}^2\ot y_{(s+1\,k)}^2 \\
& \qquad + g_{(j\,k)}^2\ot y_{(j\,k)}^2 + \sum_{i<s<t\le k+1} \mathtt{z}_{(j\,k)}(s,t) \ot  y_{(t\,k)}y_{(s\,k)}. \label{eq:coaccion yjk2}
\end{align}
\end{lem}
\pf
As $\up$ is an algebra map,
\begin{align*}
\up(y_{(j\,k)}^2)= \left(\tz_{(j\,k)}\ot 1+ g_{(j\,k)}\ot y_{(j\,k)} + 2 \sum_{s=j}^{k-1} \tz_{(j\,s)}g_{(s+1\,k)}\ot y_{(s+1\,k)}\right)^2.
\end{align*}
Then we apply Lemma \ref{lem:corchete yik contra yjk} to write the right hand side as a linear combination of elements $y_{(t\,k)}y_{(s\,k)}$
(remember that $y_{(k+1\,k)}=1$).
\epf

Notice that $\mH_{m+1}=\mH_m/I_{m+1}$ is such that $I_{m+1}$ is generated by skew primitive elements \cite[Remark 6.10]{AS2}.
According to Proposition \ref{pro:summary} cf. \eqref{eqn:u-tilde}, to describe $\mL_{m+1}$ as a quotient of $\mL_m$ we need the \emph{deforming elements} $\uvi_{(j\,k)}$ defined by the equation: 
\begin{align}\label{eqn:dl-2}
\tz_{(j\,k)}^2\ot 1-\up(y_{(j\,k)})^2=\uvi_{(j\,k)}\ot 1.
\end{align}

Recall the definition of $\tz_{(j\,k)}$ in \eqref{eq:def ztilde}.

\begin{rem}\label{rem:zeta}
As in the case of $x_{(jk)}$, $y_{(jk)}$, we define recursively $a_{(jj)}=a_j$, $a_{(jk)}=[a_j,a_{(j+1\,k)}]_c$. By induction we see that
\begin{align*}
\tz_{(jk)}&=a_{(jk)}+ \mbox{ other terms with factors }a_{(st)}, \, t-s<k-j.
\end{align*}
For example,
\begin{align*}
\tz_{(12)}&=a_{(12)}, & \tz_{(13)}=a_{(13)}+2\lambda_{13}q_{12}a_2g_1g_3.
\end{align*}
\end{rem}

\begin{theorem}\label{thm:L(A,H)-N=2}
The algebra 
$\mL(\boldsymbol\lambda,\bs\mu,\bs\nu):=L(\mA(\boldsymbol\lambda,\bs\mu,\bs\nu)
,\mH)$ is the quotient of 
$\mhL(\boldsymbol\lambda,\bs\mu,\bs\nu)$ by
\begin{align}\label{eq:rels-liftings-N2}
\tz_{(j\,k)}^2&=\mu_{(j\,k)}(1-g_{(j\,k)}^2)+ \uvi_{(j\,k)},
\end{align}
where $\uvi_{(j\,k)}$ is defined recursively as: $\uvi_{(kk)}=0$, $k\in \I$, and for $k>j$
\begin{align*}
\uvi_{(j\,k)}&=-4\sum_{j\leq p<k} \chi_{p+1,k}(g_{j,p})  \mu_{(p+1\,l)}
\Big(\uvi_{(j\,p)}+\mu_{(j\,p)}(1-g_{(j\,p)}^2)\Big)g_{(p+1\,k)}^2.
\end{align*}
\end{theorem}
\pf We prove the statement by induction in $m=k-j$. We work over $\mH_m$, $\mA_m$, $\mL_m$. Then $x_{(j\,k)}^2$ is primitive and 
$\gamma(x_{(j\,k)}^2)=y_{(j\,k)}^2$.
Set $\pi_m:\mhL\twoheadrightarrow\mL_m$ the canonical projection. By \eqref{eq:coaccion yjk2},
\begin{align*}
\up(y_{(j\,k)}^2)&= \big(\tz_{(j\,k)}^2- \uvi_{(j\,k)}\big)\ot 1+ g_{(j\,k)}^2\ot y_{(j\,k)}^2 \\
& \qquad + \sum_{i<s<t\le k+1} \pi_m\big(\mathtt{z}_{(j\,k)}(s,t)\big) \ot  y_{(t\,k)}y_{(s\,k)}.
\end{align*}
By Proposition \ref{pro:summary}, $\pi_m\big(\mathtt{z}_{(j\,k)}(s,t)\big)=0$ and the theorem follows.
\epf

\begin{exa}\label{ej:N=2,theta=3}
For $\theta=2,3$ we have the following relations
\begin{align*}
\tz_{(12)}^2& -\mu_{(12)}(1-g_{12}^2)- \uvi_{(12)} \\
&= a_{(12)}^2-\mu_{(12)}(1-g_{12}^2)-4q_{21}\mu_{(1)}\mu_{(2)}(1-g_1^2)g_2^2, \\
\tz_{(13)}^2& -\mu_{(13)}(1-g_{(1\,3)}^2)- \uvi_{(13)} \\
&= (a_{(13)}+2\lambda_{13}q_{12}a_2g_{13})^2 -\mu_{(13)}(1-g_{(1\,3)}^2)- \uvi_{(13)} \\
&=a_{(13)}^2 + 2q_{12}\lambda_{13}\nu_2(1-g_{(13)}g_2)g_{13} +4q_{12}^2\lambda_{13}^2\mu_{(2)}(1-g_2^2)g_{13}^2 \\
& \qquad  -\mu_{(13)}(1-g_{(1\,3)}^2)  - 4q_{31}q_{21} \mu_{(23)}\mu_{(1)}(1-g_1^2) g_{23}^2  \\
& \qquad -4q_{31}q_{32}\mu_{(3)} \Big(\mu_{(12)}(1-g_{12}^2)+4q_{21}\mu_{(1)}\mu_{(2)}(1-g_1^2)g_2^2 \Big)  g_3^2 .
\end{align*}
\end{exa}

\begin{rem}\label{rem:no todos 0}
We set $G=(\Z/2n\Z)^3$ for some $n\ge 2$, $g_i$, $i=1,2,3$, the generators of each cyclic factor. Set $H=\k G$. Given the matrix
$\mathbf{q}=\left(\begin{smallmatrix} -1 & 1 & -1 \\ -1 & -1 & -1 \\ -1 & 1& -1 \end{smallmatrix}\right)$,
there exist $\chi_i\in\widehat{\Gamma}$, $i=1,2,3$, such that $\chi_j(g_i)=q_{ij}$, so $V$ is realized in $\ydh$. Notice that $\chi_i^2=\eps$, $i=1,2,3$,
$\chi_1\chi_3=\eps$, so the scalars $\mu_{(i)}$, $\lambda_{13}$ can be simultaneously non-zero by \eqref{eq:condiciones escalares N=2}
for the Yetter-Drinfeld module $V$ with basis $v_i\in V_{g_i}^{\chi_i}$.

But for $\theta\ge4$ we have that $\lambda_{13}\lambda_{24}=0$. Indeed $\chi_{(14)}(g_{(14)})=-1$, so
either $\chi_1\chi_3\neq\eps$ or else $\chi_2\chi_4\neq\eps$.

Also $\nu_2\nu_3=0$. Otherwise $\chi_1\chi_2^3\chi_3=\chi_2\chi_3^2\chi_4$, which implies that $\chi_{12}=\chi_{34}$, so 
$\chi_{(14)}(g_{(14)})=\chi_{12}^2(g_{12})\chi_{34}^2(g_{34})=1$, a contradiction.
\end{rem}

\begin{rem}
The computation of $\tz_{(jk)}^2$, $j<k$, involves the computation of $[a_{(rs)},a_{(r's')}]_c$ for $r\le s$, $r'\le s'$.
A general abstract formula involves all the scalars $\mu_{rs}$, $\nu_s$, $\lambda_{st}$. However not all of them can be non-zero simultaneously,
see Remark \ref{rem:no todos 0}.
For example,
\begin{align*}
[a_1,a_{(35)}]_c &= 2\lambda_{13}a_{45}-2\chi_{34}(g_1)\lambda_{15}a_{34}g_{15}+4\chi_3(g_1)\lambda_{14}\lambda_{35}(1-g_{35}) \\
[a_{(14)},a_3]_c &= 2\lambda_{13}\chi_3(g_{(24)})a_{(24)} - 2 \nu_3 a_1 g_{2343} + 8 \lambda_{24}\mu_{(3)} \chi_3(g_2) a_1g_{24}(g_3^2-1), \\
[a_{(14)},a_2]_c &= 2\lambda_{24}a_{(13)}g_{24} + 2\chi_2(g_4) \nu_2 a_4 + 4 \lambda_{24} \chi_3(g_{12}) a_3a_{12}g_{24} \\
&  + 4 \lambda_{24} \chi_{23}(g_1)a_{23}a_1g_{24} -8 \lambda_{24} \chi_{23}(g_{12})a_3a_2a_1 g_{24}.
\end{align*}
In the first identity we necessarily have $\lambda_{13}\lambda_{14}=\lambda_{15}\lambda_{14}=0$. Similar conditions follow for the other two 
identities.
\end{rem}

\begin{exa}\label{ej:not-in-H-N2}
Set $\theta=5$ and consider the braiding matrix
$$\left(\begin{smallmatrix}
\text{--}1 & \text{--}1 & 1  & \text{--}1 & 1 \\
1  & \text{--}1 & 1  & 1  & 1 \\
1  & \text{--}1 & \text{--}1 & 1  & \text{--}1\\
\text{--}1 & 1  & \text{--}1 & \text{--}1 & \text{--}1\\
1  & 1 & \text{--}1 & 1  & \text{--}1
\end{smallmatrix}\right).$$
Thus $\chi_{13}$, $\chi_{15}$, $\chi_{24}$, $\chi_{(13)}\chi_2$, $\chi_{(24)}\chi_3\neq \eps$. We may assume there are $g_i$, $\chi_j$ are such that
\begin{align*}
\chi_{14}&=\chi_{35}=\chi_{(35)}\chi_4=\chi_i^2= \eps, & 1\le i&\le5.
\end{align*}
Notice that $\tz_{(13)}=a_{(13)}$, $\tz_{(24)}=a_{(24)}$, but
\begin{align*}
\tz_{(14)} &= a_{(14)}-4\lambda_{14}a_{23}g_{14}, &
\tz_{(15)} &= a_{(15)}+4\lambda_{14}\lambda_{35}a_2g_{35}+4\lambda_{35}a_4a_{12}g_{35}, \\
\tz_{(25)} &= a_{(25)}+4\lambda_{35}a_4a_2g_{35}.
\end{align*}

We compute the relations involving $\tz_{(i,i+1)}^2$, $\tz_{(i,i+2)}^2$ as in Example \ref{ej:N=2,theta=3}.
For the other cases we need some auxiliary computations as
\begin{align*}
[a_{(14)},a_{23}]_c&= [a_2,a_{(25)}]_c = [a_{(25)},a_4]_c = 0, \\
[a_{(15)},a_2]_c&= 8 \lambda_{14}\lambda_{35}\mu_{(2)} (1-g_2^2)(1+g_{14})g_{35}, \\
[a_{12},a_{(15)}]_c &= [a_1, [a_2,a_{(15)}]_c]_c=0, \\
[a_{(15)},a_4]_c&= 8\lambda_{35}\mu_{(4)}a_{12}g_{35}(g_4^2-1)-2\lambda_{14}a_{(25)}g_{14}-2\nu_4a_{12}g_{(35)}g_4.
\end{align*}
We have that
\begin{align*}
\tz_{(14)}^2 &= a_{(14)}^2+16\lambda_{14}^2 \big(\mu_{(23)}(1-g_{23}^2)-2\mu_{(2)}\mu_{(3)}(1-g_2^2)g_3^2\big)g_{14}^2, \\
\tz_{(25)}^2 &= a_{(25)}^2+16\lambda_{35}^2 \mu_{(4)}\mu_{(2)}(1-g_4^2) (1-g_2^2) g_{35}^2, \\
\tz_{(15)}^2 &= a_{(15)}^2-16\lambda_{14}^2\lambda_{35}^2\mu_{(2)}(1-g_2^2)g_{35}^2
\\ & \quad -16\lambda_{35}^2\mu_{(4)}(1-g_4^2)
\big(\mu_{(12)}(1-g_{12}^2)+2\mu_{(1)}\mu_{(2)}(1-g_1^2)g_2^2\big)g_{35}^2  \\
& \quad +8 \lambda_{35}a_4a_{(15)}a_{12}g_{35} + 32\lambda_{14}^2\lambda_{35}^2 \mu_{(2)} (1-g_2^2)(1+g_{14})g_{35}^2  \\
& \quad -8\nu_4\lambda_{35}\big(\mu_{(12)}(1-g_{12}^2)+2\mu_{(1)}\mu_{(2)}(1-g_1^2)g_2^2\big)g_{(35)}^2  \\
& \quad - 8 \lambda_{14}\lambda_{35}a_{(25)}a_{12}g_{14}g_{35}-32 \lambda_{14}\lambda_{35}^2 a_4a_2a_{12} g_{35}^2,
\end{align*}
so last relation can be deformed in higher strata of the coradical filtration.
\end{exa}

\section{The case $N=3$}\label{sec:N=3}

Let $H$, $V$ as in \S \ref{subsec:gral-context}, $\Gamma$ as in \eqref{eqn:Gamma}. Assume that $V$ is of type $A_\theta$, $\theta\in\N$, associated to (primitive) cubic root of unity $\xi$. Let $\B(V)$ be the corresponding 
Nichols algebra.
In this section we compute the liftings of $V$. We show that all of them arise as cocycle deformations of $\B(V)\# H$.

\subsection{Cleft objects}\label{subsec:galois-3}

Let us set, following Proposition \ref{pro:summary}, $\mtA=\mtA(\boldsymbol\lambda)$
the quotient of $T(V)\# H$ by the
relations
\begin{align}\label{eq:rels-Atheta-N3-galois}
y_{ij} &= 0, \ i < j - 1; & y_{iij} &= \lambda_{iij}, \ \vert j - i\vert = 1
\end{align}
for some family of scalars $\boldsymbol\lambda=(\lambda_{iij})$ satisfying
\begin{align}\label{eqn:restriction-1}
 \lambda_{iij}=0 \quad \text{ if } \chi_{iij}\neq\eps.
\end{align}
Here, we have renamed the basis $\{x_1,\dots,x_\theta\}$ of $V$ by $\{y_1,\dots,y_\theta\}$.

If $\mtA\neq 0$, then these algebras are cleft objects for $\mtH$, by Proposition \ref{pro:summary}.
The coaction
$\rho:\mtA\to \mtA\ot\mtH$ given by
$$
\rho(y_i)=y_i\ot 1+ g_i\ot y_i, \qquad i \in \I.
$$
We will show that $\mtA(\boldsymbol\lambda)\neq 0$ for every
$\boldsymbol\lambda$ satisfying \eqref{eqn:restriction-1}. In particular, this will show that \eqref{eqn:condition-recursive} holds for $j=0$.
First, we develop some
technicalities about these scalars.

\begin{lem}\label{lem:l112}
\begin{enumerate}
\item Let $1\leq i\neq j< \theta$, $\mid i-j\mid=1$. If $\chi_{iij}=\eps$, then
$q_{ij}=q_{ji}=\xi$.
 \item\label{eq:l1l2=0} Let $1\leq i< \theta-2$, $j=i+1$, $k=i+2$.  If
$\chi_{iij}=\eps$ or
$\chi_{ijj}=\eps$, then $\chi_{jjk}\neq \eps$ and $\chi_{jkk}\neq \eps$.
\item\label{eq:l1l3=0} Let $1\leq i< \theta-3$, $j=i+1$, $k=i+2$, $l=i+3$.  If
$\chi_{iij}=\eps$ or
$\chi_{ijj}=\eps$, then $\chi_{kkl}\neq \eps$ and $\chi_{kll}\neq \eps$.
\end{enumerate}
\end{lem}
\pf
We set $i=1$ to simplify the notation. For (1), observe that $\chi_{112}=\eps$ gives
$1=\chi_{112}(g_1)=\xi^2q_{12}$ and $1=\chi_{112}(g_2)=\xi^2q_{21}$. Hence
$q_{21}=q_{12}=\xi$. Idem for $\chi_{122}=\eps$.
For (2), we have
$$
\chi_{112}(g_{223})\chi_{223}(g_{112})=(q_{31}q_{13})^2(q_{21}q_{12})^4(q_{23}q_{32})q_{22}^4=\xi^2.
$$
Hence, if $\chi_{112}=\eps$, then $\chi_{223}\neq\eps$. The other combinations follow analogously.
For (3), it follows that $\chi_{112}(g_{334})\chi_{334}(g_{112})=\xi$. Thus if $\chi_{112}=\eps$, then $\chi_{334}\neq\eps$. A similar computation
yields the
other combinations.
\epf

\begin{pro}\label{pro:galois}
Let $\boldsymbol\lambda=(\lambda_{iij})$ satisfy
\eqref{eqn:restriction-1}. Then  $\mtA(\boldsymbol\lambda)\neq 0$. Hence
Hence $\mtA(\boldsymbol\lambda)\in\Cleft\mtH$ and, in particular,
\begin{align*}
\Cleft'\mtH=\{\mtA(\boldsymbol\lambda)|\boldsymbol\lambda\text{ as in }\eqref{eqn:restriction-1}\}.
\end{align*}
\end{pro}
\pf
Set $\mtA=\mtA(\boldsymbol\lambda)$. Observe that
$\mtA\simeq \mtE_\theta\# H$, for
$$
\mtE_\theta=\k\langle y_1, \dots, y_{\theta}\mid y_{ij}, \, i < j - 1; \ y_{iij} -
\lambda_{iij}, \, \vert j - i\vert = 1\rangle
$$
as the ideal of $T(V)$ generated by \eqref{eq:rels-Atheta-N3-galois} is an object in $\ydh$.

Hence we need to show that $\mtE_\theta\neq0$. We see this by induction on $\theta\geq
2$. For $\theta=2$, we distinguish three cases, namely
\begin{align*}
 &\mathrm{(i)}\, \lambda_{112}=\lambda_{122}=0; &&\mathrm{(ii)}\, \lambda_{112}\neq 0,
\lambda_{122}=0; &&\mathrm{(iii)}\, \lambda_{112}\lambda_{122}\neq0.
\end{align*}
Case (i) is clear, as $\mtE_2$ is the distinguished pre-Nichols algebra of type
$A_2$, cf. p. \pageref{eq:rels-Atheta-c}. For both cases
(ii) and (iii) notice that, as $q_{12}=q_{21}=\xi$ by Lemma \ref{lem:l112}, the defining
relations become:
\begin{align*}
&\lambda_{112}=y_1^2y_2+y_1y_2y_1+y_2y_1^2, && \lambda_{122}=y_2^2y_1+y_2y_1y_2+y_1y_2^2.
\end{align*}
Now case (iii) follows by observing that $\alpha:\mtE_2\to \k$ given by
\begin{align}\label{eqn:repr1}
&\alpha(y_1)=\left(\frac{\lambda_{112}^2}{3\lambda_{122}}\right)^{\frac{1}{3}},
&&\alpha(y_2)=\left(\frac{\lambda_{122}^2}{3\lambda_{112}}\right)^{\frac{1}{3}}
\end{align}
is a well defined one-dimensional representation.

For case (ii), we have the representation $\alpha:\mtE_2\to \k^{3\times 3}$:
\begin{align}\label{eqn:repr2}
&\alpha(y_1)=\left(\begin{smallmatrix}
                 0&1&0\\ 0&0&1\\0&0&0
                 \end{smallmatrix}
\right),
&&\alpha(y_2)=\left(\begin{smallmatrix}
                  -\lambda_{112}&0&-\lambda_{112}\\ 0&0&0\\ \lambda_{112}&0&\lambda_{112}
                 \end{smallmatrix}
\right).
\end{align}

We turn now to the general case $\theta\geq 3$. If $\lambda_{112}=\lambda_{122}=0$, then
we have an algebra isomorphism
$$
\mtE_\theta/\langle y_1 \rangle \simeq \mtE_{\theta-1},\qquad  \bar{y_i}\mapsto a_{i-1},
\quad 2\leq i\leq \theta
$$
where $a_i$, $1\leq i\leq\theta-1$ stand for the generators of $\mtE_{\theta-1}$. Hence
we may assume $\lambda_{112}\neq 0$. Thus it follows that
$\lambda_{223}=\lambda_{233}=0$ and (if $\theta\geq 4$) $\lambda_{334}=\lambda_{344}=0$,
by Lemma \ref{lem:l112}. In particular, if $\theta=3$, there is an isomorphism
$$
\mtE_3/\langle y_3 \rangle \simeq \mtE_{2},\qquad  \overline{y_i}\mapsto a_{i},
\quad 1\leq i\leq 2
$$
which shows that $\mtE_3\neq 0$. Let then $\theta\geq 4$ and assume $\mtE_\vartheta\neq0$
for every $\vartheta<\theta$. For $i\in\N$, we set $i^*=i+3$.
Let $\mQ=(q_{ij})_{1\leq i,j\leq\theta}$ be the braiding matrix and consider the
submatrix $\mQ'=(q'_{ij})_{1\leq i,j\leq\theta-3}$, $q'_{ij}=q_{i^*j^*}$ and the
subfamily $\boldsymbol\lambda'=(\lambda_{iij}')$, $\lambda'_{iij}=\lambda_{i^*i^*j^*}$.
Consider the corresponding algebra
$\mtE_{\theta-3}(\boldsymbol\lambda')$, with generators $a_1,\dots, a_{\theta-3}$. We
denote this algebra by $B$ and rename the generators $w_{i+3}:=a_i$, $1\leq i\leq
\theta-3$. That is, $B$ is the algebra generated by $w_4,\dots, w_\theta$ with relations
defined by a subset of the relations in \eqref{eq:rels-Atheta-N3-galois}. Let $A$ be the
algebra $\mtE_2(\lambda_{112},\lambda_{122})$, with generators $s_1,s_2$. Let us also set
$E=\mtE_\theta/\langle y_3\rangle$. We will show that $E\neq 0$, hence $\mtE_\theta\neq
0$. Let us denote by $y_1,y_2,y_4,\dots, y_\theta\in E$ the images of the corresponding
generators of $\mtE_\theta$.

Let $\alpha:A\to k^{m\times m}$ be the representation defined in \eqref{eqn:repr1} if
$\lambda_{122}\neq0$ or the representation in \eqref{eqn:repr2} if
$\lambda_{122}=0$, with $m=1$ or $m=3$, respectively. Set, accordingly, $M_1=\alpha(s_1),
M_2=\alpha(s_2)$.
Let us consider the vector space $W=B\ot V$. We define
$\varrho:E \to \End(W)$ given by
\begin{align*}
 \varrho(y_j)(b\ot v)=
\begin{cases}
g_j\cdot b\ot M_j\,v, & j<3\\
w_jb\ot v, & j>3.
\end{cases}
\end{align*}
This is a well-defined representation. For instance:
\begin{align*}
&\varrho(y_{112}-\lambda_{112})(w_k\ot
v)=(\chi_{k}(g_{112})\lambda_{112}-\lambda_{112})w_k\ot v=0,\\
&\varrho(y_1y_{j}-q_{jk}y_jy_1)(w_k\ot
v)=(q_{1j}q_{jk}-q_{jk}q_{1j})w_jw_k\ot
M_1\,v=0, \quad j>3,\\
&\varrho(y_{jjl}-\lambda_{jjl})(w_k\ot
v)=w_{jjl}w_k\ot v-\lambda_{jjl}w_k\ot v=0, \quad 3<j=l-1.
\end{align*}
Hence the proposition follows.
\epf


We need to compute $\rho(y_{(kl)}^3)=\rho(y_{(kl)})^3$, $1\leq k<l\leq
\theta$. We have:
\begin{align}\label{eqn:l=1}
 \rho(y_i)^3=y_i^3\ot 1 + g_i^3\ot x_i^3, \quad i \in \I
\end{align}
as these elements are skew-primitives for the coaction. Now, for every
$k,l$,
\begin{align*}
\rho(y_{(kl)})&=y_{(kl)}\ot 1+g_{(kl)}\ot
x_{(kl)}+(1-\xi^2)\sum_{k\leq p<l}y_{(kp)}g_{(p+1l)}\ot x_{(p+1l)},
\end{align*}
again as relation $y_{ij} = 0$, $i < j - 1$ holds in $\mtA$.

Consider a family of indeterminate variables ${\bf t}=(\tt_{iij})$ and let $R=\k[{\bf t}]$ be the polynomial
ring on those variables. Let us denote by $\mtH_R$ and $\mtA_R$ the $R$-algebras defined by the same relations as $\mtH$ and $\mtA$.
Given a family $\boldsymbol\lambda$ as above, we consider the evaluation $\tt_{iij}\mapsto \lambda_{iij}$. Thus, we have $\mtA_R\ot_R\,\k=\mtA$,
$\mtH_R\ot_R\,\k=\mtH$.
Let $R^+$ be
the ideal generated by ${\bf t}$. By \eqref{eq:rels-Atheta-N3-galois} we have:
\begin{align*}
y_{(kr)}y_{(ks)}=\chi_{(ks)}(g_{(kr)})y_{(ks)}y_{(kr)}+R^+\mtA_R.
\end{align*}
If $\lambda_{iij}=0$  for
every $1\leq i,j\leq
\theta$, then $\mtA\simeq \mtH$ and thus
\begin{align*}
\rho(y_{(k\,l)}^3)=y_{(k\,l)}^3\ot 1+ g_{(k\,l)}^3\ot x_{(k\,l)}^3
+\sum_{k\leq p<l}C{p}y_{(k\,p)}^3g_{(p+1\,l)}^3\ot x_{(p+1\,l)}^3,
\end{align*}
for $B_p$ as in \eqref{eqn:Cp}.
Hence, 
\begin{multline*}
\rho(y_{(k\,l)}^3)=y_{(k\,l)}^3\ot 1+ g_{(k\,l)}^3\ot x_{(k\,l)}^3
\\+\sum_{k\leq p<l}C_{p}y_{(k\,p)}^3g_{(p+1\,l)}^3\ot x_{(p+1\,l)}^3+R^+\mtA_R\ot \mtH_R.
\end{multline*}
Hence, in the computation of $\rho(y_{(k\,l)}^3)$ in the general case, we need to focus
on the terms in which a scalar $\lambda_{***}$ may appear. See Example \ref{exa:ABC} for
$\theta=2$. This example shows the
philosophy behind our calculations. Also, it introduces the notation $\rightsquigarrow $ in \eqref{eqn:approx} to keep only the terms with a
factor $\lambda_{***}$.

\begin{exa}\label{exa:ABC}
We will show that
\begin{align*}
 \rho(y_{12}^3)=y_{12}^3\ot 1+g_{12}^3\ot x_{12}^3+(1-\xi^2)^3\chi_1(g_2)^3y_1^3g_2^3\ot x_2^3.
\end{align*}
Hence, we can take a section $\gamma:\mtH\to\mtA$ such that $\gamma(x_{12}^3)=y_{12}^3$.

Let us compute $\rho(y_{12})^3$. Set
\begin{align*}
&A=y_{12}\ot 1, && B=g_{12}\ot x_{12}, && C=y_1g_2\ot x_2.
\end{align*}
Hence $\rho(y_{12})=A+B+(1-\xi^2)C$. As said, we need to focus on the terms in which a
factor $\lambda_{***}$ may appear. These are related with the fact that $y_1$ appears to
the left of $y_{12}$ and are precisely:
\begin{align*}
 &CAB, && CBA, && BCA,  &CAA, && ACA,  &&CAC, && CCA.
\end{align*}
Now, for instance
\begin{align*}
CAB&=y_1g_2y_{12}g_{12}\ot x_2x_{12}=\chi_{12}(g_2)y_1y_{12}g_{122}\ot x_2x_{12}\\
&=\lambda_{112}\chi_{12}(g_2)g_{122}\ot
x_2x_{12}+\chi_{12}(g_2)\chi_{12}(g_1)y_{12}y_1g_{122}\ot x_2x_{12}.
\end{align*}
We only need to keep the term involving $\lambda_{112}$. Hence, we write
\begin{align}\label{eqn:approx}
CAB\rightsquigarrow  \lambda_{112}\chi_{12}(g_2)g_{122}\ot
x_2x_{12}=\lambda_{112}\xi^2g_{122}\ot x_2x_{12}.
\end{align}
as $\lambda_{112}\chi_{12}(g_2)=\lambda_{112}\xi^2$.  We will do this for every term. We need the following equalities:
\begin{align*}
y_1y_{12}&=\lambda_{112}+\chi_{12}(g_1)y_{12}y_1\rightsquigarrow  \lambda_{112};\\
y_1y_{12}^2&=\lambda_{112}(1+\xi^2)y_{12}+\chi_{12}(g_1)^2y_{12}^2y_1\rightsquigarrow
\lambda_{112}(1+\xi^2)y_{12};\\
y_{12}y_1y_{12}&=\lambda_{112}y_{12}+\chi_{12}(g_1)y_{12}^2y_1\rightsquigarrow
\lambda_{112}y_{12}\\
y_1y_{12}y_1&=\lambda_{112}y_1+\chi_{12}(g_1)y_{12}y_1^2\rightsquigarrow  \lambda_{112}y_1\\
y_1^2y_{12}&=\lambda_{112}(1+\xi^2)y_1+\chi_{12}(g_1)^2y_{12}y_1^2\rightsquigarrow
\lambda_{112}(1+\xi^2)y_1.
\end{align*}
We have:
\begin{align*}
&CAB\rightsquigarrow  \lambda_{112}\xi^2g_{122}\ot x_2x_{12}; &&
CBA\rightsquigarrow  \lambda_{112}g_{122}\ot x_2x_{12};\\
&BCA\rightsquigarrow  \lambda_{112}\xi g_{122}\ot x_2x_{12}.
\end{align*}
Thus $CAB+CBA+BCA\rightsquigarrow  0$.
\begin{align*}
&CAA\rightsquigarrow  \lambda_{112}\xi(1+\xi^2)y_{12}g_2\ot x_2; &&
ACA\rightsquigarrow  \lambda_{112}\xi^2y_{12}g_2\ot x_2.
\end{align*}
Thus $CAA+ACA\rightsquigarrow  0$.
\begin{align*}
&CAC\rightsquigarrow  \lambda_{112}y_1g_2^2\ot x_2^2; &&
CCA\rightsquigarrow  \lambda_{112}\xi^2(1+\xi^2)y_1g_2^2\ot x_2^2.
\end{align*}
Thus $CAC+CCA\rightsquigarrow  0$.
Therefore,
$$
\rho(y_{12})^3=y_{12}^3\ot 1+g_{12}^3\ot
x_{12}^3+(1-\xi^2)^3\chi_1(g_2)^3y_{1}^3g_{2}^3\ot
x_{2}^3.
$$
\qed
\end{exa}

From now on we consider the case $\theta\geq 3$.
We will collect some technical
identities needed to compute $\rho(y_{(kl)})^3$ in a series of general
lemmas.

\begin{lem}\label{lem:tech1}
The following identities hold in $\mtA$:
$$
[y_{(1\,l)},y_{2}]_c=\lambda_{122}(1-\xi^2)\chi_2(g_{(3\,l)})y_{(3\,l)},
\qquad l\geq 3.$$
\end{lem}
\pf
Assume  first $l=3$. We have two cases, namely $\chi_{122}=\eps$ or not, (in
which
case it is possible to have $\chi_{223}=\eps$). We proceed as in  \cite[Lemma 1.11]{AS2}.
In the first case, we have the lemma. In the second, we get $
[y_{(13)},y_{2}]_c=\lambda_{223}(1-\chi_{223}(g_1))x_1=0$, hence the lemma also
holds
(as
$\lambda_{122}=0$). For the general case we get:
\begin{align*}
[y_{(1\,l)},y_{2}]_c&=[[y_{(1\,3)},y_{(4\,l)}]_c,y_2]_c\\
&=\chi_2(g_{(4\,l)})[y_{(1\,3)},y_2]_cy_{
(4\,l)}-\chi_{(4\,l)}(g_{(1\,3)})y_{(4\,l)}[y_{(1\,3)},y_2]_c\\
&=\lambda_{122}q_{32}(1-\xi^2)\chi_2(g_{(4\,l)})
\Big(y_3y_{(4\,l)}-\chi_{(4\,l)}(g_3)\chi_{(4\,l)}(g_{122})y_{(4\,l)}y_3\Big)\\
&=\lambda_{122}(1-\xi^2)\chi_2(g_{(3\,l)})y_{(3\,l)},
\end{align*}
as
$\lambda_{122}\chi_{(5\,l)}(g_{122})=\lambda_{122}$, $1=\chi_2(g_{(4\,l)})\chi_{(4\,l)}(g_2)$
plus q-Jacobi \eqref{eq:qjacob}.
\epf

\begin{lem}\label{lem:tech2}
The following identities hold in $\mtA$, for $1\leq p\leq l$:
\item\label{it4}
\begin{align*}
 [y_{p},y_{(1\,l)}]_c=\begin{cases}
                       \lambda_{112}(1-\xi^2)y_{(3\,l)}, & p=1;\\
\lambda_{122}(1-\xi)y_{(3\,l)}, & p=2;\\
0, & 2<p<l;\\
(1-\xi)y_ly_{(1\,l)}, & p=l.
                      \end{cases},\qquad  l\geq 3.
\end{align*}
\end{lem}
\pf Using that $\lambda_{112}\chi_{112}=\lambda_{112}\eps$,
\begin{align*}
[y_{1},y_{(1\,l)}]_c&=\lambda_{112}(1-\chi_{(3\,l)}(g_{112}))y_{(3\,l)}\\
&=\lambda_{112}(1-\chi_{(3\,l)}(g_{112})\chi_{112}(g_{(3\,l)}))y_{(3\,l)} =\lambda_{112}
(1-\xi^2)y_{(3\,l)}.
\end{align*}
If $p=2$, it follows by Lemma \ref{lem:tech1} that
$[y_2,y_{(1\,l)}]_c=\lambda_{122}(1-\xi)y_{(3\,l)}$.

Assume $2<p<l$. Then
\begin{align*}
[y_{p},y_{(1\,l)}]_c&=[y_p,[y_{(1\,p-2)},y_{(p-1\,l)}]_c]_c=-\chi_{(p-1\,l)}(g_{(1\,p-2)})
[ y_{p},y_{(p-1\,l)}]_cy_{(1\,p-2)}\\
&\qquad\qquad\qquad\qquad  +
\chi_{(1\,p-2)}(g_{p})y_{(1\,p-2)}[y_{p},y_{(p-1\,l)}]_c\\
&\overset{\text{case }p=2}{=}
-\chi_{(p-1\,l)}(g_{(1\,p-2)})\lambda_{p-1pp}(1-\xi)y_{(p+1\,l)}y_{(1\,p-2)}\\
&\qquad\qquad\qquad\qquad +
\chi_{(1\,p-2)}(g_{p})\lambda_{p-1pp}(1-\xi)
y_{(1\,p-2)}y_{(p+1\,l)}\\
&=\lambda_{p-1pp}(1-\xi)\chi_{(1\,p-2)}(g_{p})\Big(y_{(1\,p-2)}y_{(p+1\,l)}\\
&\qquad\qquad\qquad -
\chi_{(p+1\,l)}(g_{(1\,p-2)})
\chi_{p-1pp}(g_{(1\,p-2)})y_{(p+1\,l)}y_{(1\,p-2)}\Big)\\
&=\lambda_{p-1pp}(1-\xi)\chi_{(1\,p-2)}(g_{p})[y_{(1\,p-2)},y_{(p+1\,l)}]_c=0.
\end{align*}
Finally, if $p=l$, using that $
 [y_l,y_{l-1\,l}]_c=(1-\xi)y_ly_{l-1\,l}-\lambda_{l-1ll}\xi^2$ and q-Jacobi
\eqref{eq:qjacob} we arrive to
\begin{align*}
 [y_{l},y_{(1\,l)}]_c&=-\chi_{l-1l}(g_{(1\,l-2)})\Big((1-\xi)y_ly_{l-1\,l}-\lambda_{l-1ll}
\xi^2\Big)y_{(1\,l-2)}\\
&\qquad +\chi_{(1\,l-2)}(g_{l})y_{(1\,l-2)}\Big((1-\xi)y_ly_{l-1\,l}-\lambda_{l-1ll}
\xi^2\Big)\\
&=(1-\xi)y_ly_{(1\,l)}\\
&\qquad +\lambda_{l-1ll}\xi^2\chi_{l-1l}(g_{(1\,l-2)})(1-\chi_{(1\,l-2)}
(g_l)\chi_l(g_{(1\,l-2)}))y_{(1\,l-2)},
\end{align*}
and the lemma follows using that
$\lambda_{l-1ll}\chi_{l-1l}(g_{(1\,l-2)})^{-1}=\lambda_{l-1ll}\chi_{l}(g_{(1\,l-2)})$ and
$\chi_{l}(g_{(1\,l-2)})\chi_{(1\,l-2)}(g_l)=1$.
\epf

\begin{rem}\label{rem:1}
If $l=2$, then $[y_{1},y_{(1\,l)}]_c=\lambda_{112}$.
\end{rem}

\begin{rem}\label{rem:2}
If $2< p\leq l$, then
$$[y_{(1\,l)},y_p]_c=0.$$
Indeed, if $p=l$,
\begin{align*}
[y_{(1\,l)},y_l]_c&=(1-\chi_l(g_{(1\,l)})\chi_{(1\,l)}(g_l))y_{(1\,l)}y_l-\chi_l(g_{(1\,l)
})[y_l,y_{(1\,l)}]_c\\
&=(1-\xi)y_{(1\,l)}y_l-\chi_l(g_{(1\,l)})(1-\xi)y_ly_{(1\,l)}\\
&=(1-\xi)[y_{(1\,l)},y_l]_c.
\end{align*}
If $1<p<l$ this follows from
$$
[y_{(1\,l)},y_p]_c=(1-\chi_p(g_{(1\,l)}
)\chi_{(1\,l)}(g_p))y_{(1\,l)}y_p-\chi_p(g_{(1\,l)})[y_p,y_{(1\,l)}]_c=0.
$$
\end{rem}

\begin{rem}\label{rem:tech1}
 Let $1\leq p\leq l-2$. Then
\begin{align}\label{eqn:ind}
[y_{(1\,p+1)},y_{(1\,l)}]_c=\chi_{(1\,l)}(g_{p+1})[[y_{(1\,p)
} ,y_{(1\,l)}]_c, y_{p+1}]_c.
\end{align}
Hence,
\begin{align*}
 [y_{(1\,p)},y_{(1\,l)}]_c=\chi_{(1\,l)}(g_{(2\,p)})[\dots
[y_{1},y_{(1\,l)}]_c,y_2]_c,\dots , y_{p-1}]_c, y_{p}]_c.
\end{align*}
Indeed, using q-Jacobi \eqref{eq:qjacob} we see that \eqref{eqn:ind} holds:
\begin{align*}
[y_{(1\,p+1)}&,y_{(1\,l)}]_c=[y_{(1\,p)},[y_{p+1},y_{(1\,l)}]_c]_c+\chi_{(1\,l)}(g_{p+1})[
y_{(1\,p)},y_{(1\,l)}]_cy_{p+1}\\
&-\chi_{p+1}(g_{(1\,p)})y_{p+1}[y_{(1\,p)},y_{(1\,l)}]_c=\chi_{(1\,l)}(g_{p+1})[[y_{(1\,p)
} ,y_{(1\,l)}]_c, y_{p+1}]_c,
\end{align*}
as $[y_{(1\,p)},[y_{p+1},y_{(1\,l)}]_c]_c=0$. This last equality is clear if
$p\geq 2$ by Lemma \ref{lem:tech2} that also yields
$[y_{1},[y_{2},y_{(1\,l)}]_c]_c=\lambda_{122}(1-\xi)(y_1y_{(3\,l)}-\chi_{122}(g_1)
\chi_{(3\,l)} (g_1))=\lambda_{122}(1-\xi)[y_1,y_{(3\,l)}]_c=0.$
\end{rem}

\begin{lem}\label{lem:tech2b}
The following identities hold in $\mtA$.
\begin{enumerate}
\item\label{it8} $[y_{(1\,l)},y_{(k\,p)}]_c=0$, for $3\leq k\leq p\leq l$.
\item\label{it9}
$[y_{(1\,p)},y_{(3\,l)}]_c=\chi_{(3\,p)}(g_{(1\,p)})(1-\xi^2)y_{(3\,p)}y_{(1\,l)}$, for
$3\leq p< l$.
\end{enumerate}
\end{lem}

\pf
\eqref{it8} Fix $k$. Recall that $[y_{(1\,l)},y_j]_c=0$, for
$3\leq j\leq l$, by Remark \ref{rem:2}. In particular,
$[y_{(1\,l)},y_k]_c=0$. Now, using induction on $p$ and q-Jacobi
\eqref{eq:qjacob},
\begin{multline*}
[y_{(1\,l)},y_{(k\,p)}]_c=-\chi_{p}(g_{(k\,p-1)})[y_{(1\,l)},y_{p}]_cy_{(k\,p-1)}\\+
 \chi_{(k\,p-1)}(g_{(1\,l)})y_{(k\,p-1)}[y_{(1\,l)},y_{p}]_c=0.
\end{multline*}
\eqref{it9} We have, using q-Jacobi \eqref{eq:qjacob} and Item \eqref{it8} for $k=3$:
\begin{align*}
 [y_{(1\,p)},y_{(3\,l)}]_c&=[y_{(1\,p)},[y_{(3\,p)},y_{(p+1\,l)}]_c]_c\\
&=-\chi_{(p+1\,l)}(g_{(3\,p)})y_{(1\,l)}y_{(3\,p)}+\chi_{(3\,p)}(g_{(1\,p)})
y_{(3\,p)}y_{(1\,l)}\\
&=\chi_{(3\,p)}(g_{(1\,p)})(1-\chi_{(p+1\,l)}(g_{(3\,p)})\chi_{(3\,p)
} (g_{(p+1\,l) }))y_{(3\,p)}y_ {(1\,l)}
\end{align*}
and \eqref{it9} follows as $\chi_{(p+1\,l)}(g_{(3\,p)})\chi_{(3\,p)
} (g_{(p+1\,l) })=\xi^2$.
\epf

\begin{lem}\label{lem:tech1-gral}
The following identities hold in $\mtA$:
$$[y_{(1\,l)},y_{(2\,l)}]_c=-3\lambda_{122}\chi_2(g_{(1\,l)})y_{(3\,l)}^2,
\qquad l\geq 3.$$
 \end{lem}
\pf
Follows using Lemma \ref{lem:tech1} combined with Lemma \ref{lem:tech2b} \eqref{it8}:
\begin{align*}
[y_{(1\,l)},y_{(2\,l)}]_c&=
[y_{(1\,l)},[y_2,y_{(3\,l)}]_c]_c=[[y_{(1\,l)},y_2]_c,y_{(3\,l)}]_c\\
&=\lambda_{122}(1-\xi^2)\chi_2(g_{(3\,l)})(1-\chi_{(3\,l)}(g_{122}g_{(3\,l)}))y_{(3\,l)}
^2\\
&=\lambda_{122}(1-\xi^2)\chi_2(g_{(3\,l)})(1-\xi^2)y_{(3\,l)}^2=-3\lambda_{122}
\xi^2\chi_2(g_{(3\,l)})y_{(3\,l)}^2,
\end{align*}
and thus the lemma follows as
$\lambda_{122}\chi_2(g_{(1\,l)})=\lambda_{122}\xi^2\chi_2(g_{(3\,l)})$.
\epf

\begin{lem}\label{lem:tech3}
The following identities hold in $\mtA$:
\begin{enumerate}
\item\label{it5} $[y_{1},y_{(1\,l)}]_c=
\lambda_{112}(1-\xi^2)y_{(3\,l)}$, $l\geq 3$.
\item\label{it6} $[y_{(12)},y_{(1\,l)}]_c=-3\xi^2
\lambda_{112}\chi_{(1\,l)}(g_2)y_{(3\,l)}y_2+\lambda_{112}(1-\xi)y_{(2\,l)}$,
$l\geq 3$.
\item\label{it7}  For $3\leq p <l$:
\begin{align*}
[y_{(1\,p)},y_{(1\,l)}]_c=-3\xi^2&\lambda_{112}\chi_{(1\,l)}(g_{(2\,p)})
y_{(3\,l)}y_{(2\,p)} +3\lambda_{112}\chi_{1}(g_{(3\,p)})y_{(3\,p)}y_{(2\,l)}.
\end{align*}
\end{enumerate}
\end{lem}
\pf \eqref{it5} is Lemma \ref{lem:tech2} for $p=1$. For \eqref{it6} we use
\eqref{eqn:ind}
and \eqref{it5} to get:
\begin{align*}
 [y_{12},y_{(1\,l)}]_c=\chi_{(1\,l)}(g_2)[[y_1,y_{
(1\,l)}]_c,y_2]_c=\lambda_{112}(1-\xi^2)\chi_{(1\,l)}(g_2)[y_{(3\,l)},y_2]_c.
\end{align*}
Now
$[y_{(3\,l)},y_2]_c=(1-\xi^2)y_{(3\,l)}y_2-\chi_2(g_{(3\,
l) })y_{(2\,l)}$ and \eqref{it6} follows using the equality
$\lambda_{112}\chi_2(g_{(3\,l)})\chi_{(1\,l) }(g_2)=\lambda_{112}\xi$.

For \eqref{it7}, we use q-Jacobi \eqref{eq:qjacob} and Lemma \ref{lem:tech2b} to get
\begin{align*}
[y_{(1\,p)}&,y_{(1\,l)}]_c=\chi_{(1\,l)}(g_{(3\,p)})[y_{12},y_{(1\,l)}]_cy_{(3\,p)}
-\chi_{(3\,p)}(g_{12})y_{(3\,p)}[y_{12},y_{(1\,l)}]_c\\
&=-3\xi^2\lambda_{112}\chi_{(1\,l)}(g_2)\Big(\chi_{(1\,l)}(g_{(3\,p)})
y_{(3\,l)}y_2y_{(3\,p)}-\chi_{(3\,p)}(g_{12})y_{(3\,p)} y_{(3\,l)}y_2\Big)\\
&\qquad +\lambda_{112}(1-\xi)\Big( \chi_{(1\,l)}(g_{(3\,p)}) y_{(2\,l)}y_{(3\,p)}
-\chi_{(3\,p)}(g_{12})y_{(3\,p)}y_{(2\,l)}\Big)\\
&=-3\xi^2\lambda_{112}\chi_{(1\,l)}(g_2)\chi_{(1\,l)}(g_{(3\,p)})
y_{(3\,l)}y_{(2\,p)}\\
&\qquad +3\xi^2\lambda_{112}\chi_{(1\,l)}(g_2)\chi_{(3\,p)}(g_{12})[y_{(3\,p)},
y_{(3\,l)}]_cy_2\\
&\qquad +\lambda_{112}(1-\xi)\chi_{(1\,l)}(g_{(3\,p)})[y_{(2\,l)},y_{(3\,p)}]_c
\\
&\qquad +\lambda_{112}(1-\xi)\chi_{(1\,l)}(g_{(3\,p)})\chi_{(3\,p)}(g_{(2\,l)})(1
-\xi^2)y_{(3\,p)}y_{(2\,l)}.
\end{align*}
We use this equality and Lemma \ref{lem:l112} to deduce $\lambda_{112}[y_{(3\,p)},
y_{(3\,l)}]_c=0$. We use this fact together with Lemma \ref{lem:tech2b} \eqref{it8} and
Lemma
\ref{lem:tech1} to get
\begin{align}
\notag  [y_{(1\,l)},&y_{(2\,p)}]_c=[[y_{(1\,l)},y_2]_c,y_{(3\,p)}]_c\\
\notag  &=\lambda_{122}(1-\xi^2)
\chi_2(g_{(3\,l)})(1-\chi_{(3\,p)}(g_{122}g_{(3\,l)})
\chi_{(3\,l)}(g_{(3\,p)}))y_{(3\,l)}y_{(3\,p)}\\
\label{eqn:1-2} &\qquad -\lambda_{122}(1-\xi^2)
\chi_2(g_{(3\,l)})\chi_{(3\,p)}(g_{(3\,l)})\xi^2[y_{(3\,p)},y_{(3\,l)}]_c\\
\notag &=-3\xi^2\lambda_{122}\chi_2(g_{(3\,l)})y_{(3\,l)}y_{(3\,p)}.
\end{align}
In particular, $\lambda_{112}[y_{(2\,l)},y_{(3\,p)}]_c=0$ by Lemma \ref{lem:l112}. Hence
\eqref{it7}
follows.
\epf

\begin{rem}\label{rem:1l-2p}
We have, for $2\leq p\leq l$:
\begin{align*}
 [y_{(1\,l)},y_{(2\,p)}]_c=\begin{cases}
\lambda_{122}(1-\xi^2)\chi_2(g_{(3\,l)})y_{(3\,l)}, &p=2\\
\lambda_{122}(1-\xi^2)^2\chi_2(g_{(3\,l)})y_{(3\,l)}y_{(3\,p)}, &p>2.
                          \end{cases}
\end{align*}
For $p=2$ this is Lemma \ref{lem:tech1}. For $p>2$ this is \eqref{eqn:1-2}.
\end{rem}

\begin{pro}\label{pro:right-coaction}
For $C_p$ as in \eqref{eqn:Cp} we have
 \begin{align*}
 \rho(y_{(k\,l)})^3=y_{(k\,l)}^3\ot 1+g_{(k\,l)}^3\ot x_{(k\,l)}^3
+\sum_{k\leq p<l}C_py_{(k\,p)}^3g_{(p+1\,l)}^3\ot x_{(p+1\,l)}^3.
 \end{align*}
\end{pro}
\pf
Let us set $k=1<l$ to simplify the notation. We may assume $l\geq 3$ as case $l=1$ is
\eqref{eqn:l=1} and case $l=2$ is Example \ref{exa:ABC}.
Set
\begin{align*}
A=y_{(1\,l)}\ot 1, && B=g_{(1\,l)}\ot x_{(1\,l)}, && X_p=y_{(1\,p)}g_{(p+1\,l)}\ot
x_{(p+1\,l)},
1\leq p<l.
\end{align*}
In what respects to commutation rules, we may set,
without lack of rigour, $X_l:=A$, as with the convention $g_{(l+1\,l)}=x_{(l+1\,l)}=1$ it
becomes
$$X_l=y_{(1\,l)}g_{(l+1\,l)}\ot x_{(l+1\,l)}=y_{(1\,l)} \ot1.$$
As in Example \ref{exa:ABC}, we need
to focus on the terms of $(A+B+(1-\xi^2)\sum X_p)^3$
involving a factor $\lambda_{***}$. These are divided into three big groups,
namely:
\begin{enumerate}
\item[(G1)] For every pair $p<q$, terms $XYZ$ involving $X,Y,Z\in \{B,X_p,X_q\}$, all
different, $X_p$ to the left of $X_q$.
\item[(G2)] For every pair $p<q$, terms $XYZ$ involving $X,Y,Z\in \{X_p,X_q\}$, not all
equal and with a factor $X_p$ to the left of $X_q$.
\item[(G3)] For every triple $p<q<r$, terms $XYZ$  from \emph{distinct} $X,Y,Z\in
\{X_p,X_q,X_r\}$ and with $X_p$ to the left of $X_{q}$ or $X_r$ or with $X_q$ to the left
of $X_r$.
\end{enumerate}
Since our aim is to show that there is no term involving a factor $\lambda_{***}$, we may
further restrict these groups, as the other resulting combinations
provide equivalent terms. For instance, we have
\begin{align*}
X_pX_l&=y_{(1\,p)}g_{(p+1\,l)}y_{(1\,l)}g_{(l+1\,l)}\ot
x_{(p+1\,l)}x_{(l+1\,l)}\\
&=y_{(1\,p)}g_{(p+1\,l)}y_{(1\,l)}\ot
x_{(p+1\,l)}\\
&=\chi_{(1\,l)}(g_{(p+1\,l)})y_{(1\,p)}y_{(1\,l)}g_{(p+1\,l)}\ot
x_{(p+1\,l)}.
\end{align*}
\begin{align*}
X_p&X_q=y_{(1\,p)}g_{(p+1\,l)}y_{(1\,q)}g_{(q+1\,l)}\ot
x_{(p+1\,l)}x_{(q+1\,l)}\\
&=\chi_{(1\,q)}(g_{(p+1\,l)})\chi_{(q+1\,l)}(g_{(p+1\,l)})y_{(1\,
p) } y_ { (1\, l) } g_ { (p+1\, l) }g_{(q+1\,l)} \ot x_{(q+1\,l)}x_{(p+1\,l)}\\
&=\chi_{(1\,l)}(g_{(p+1\,l)})y_{(1\,
p) } y_ { (1\, l) } g_ { (p+1\, l) }g_{(q+1\,l)} \ot x_{(q+1\,l)}x_{(p+1\,l)}.
\end{align*}
Hence we restrict to the following subgroups:
\begin{enumerate}
\item[(G1')] For every $p<l$, terms $XYZ$ involving $X,Y,Z\in \{B,X_p,A\}$, all
different, $X_p$ to the left of $A$.
\item[(G2')] For every $p<l$, terms $XYZ$ involving $X,Y,Z\in \{X_p,A\}$, not all
equal and with a factor $X_p$ to the left of $A$.
\item[(G3')] For every pair $p<q<l$, terms $XYZ$ arising from 
\emph{distinct} $X,Y,Z\in
\{X_p,X_q,A\}$ and with $X_p$ to the left of $X_{q}$ or $A$ or with $X_q$ to the left
of $A$.
\end{enumerate}
We start with group (G1'): notice that, for any $p$:
\begin{multline*}
X_pAB+X_pBA+BX_pA\\
 =(1+\xi+\xi^2)\chi_{(1\,l)}(g_{p+1\,l})y_{(1\,p)}y_{(1\,l)}g_{(1\,l)}g_{(p+1\,l)}\ot
x_{(p+1\,l)}x_{(1\,l)}=0.
\end{multline*}
We now proceed to group (G2'), {\it i.e.} . terms of the form $X_pAX_p,
X_pX_pA$ and $X_pAA, AX_pA$. We further divide:
this group into
\begin{enumerate}
\item[(G2'.1)] Factors arising from $\{A,X_1\}$.
\item[(G2'.2)] Factors arising from $\{A,X_2\}$.
\item[(G2'.3)] Factors arising from $\{A,X_p\}$, $p\geq 3$.
\end{enumerate}
The computations for item (G2'.1) are analogous to the ones in Example \ref{exa:ABC}, and
we get that the factor involving $\lambda_{***}$ is zero. For (G2'.2),
we need the following computations:
\begin{align*}
&y_{12}y_{(1\,l)}^2\rightsquigarrow
-3\lambda_{112}(1+\xi)\chi_{(3\,l)}(g_{2})
y_{(3\,l)}y_2y_{(1\,l)}+\lambda_{112}(\xi^2-\xi)y_{(2\,l)}y_{(1\,l)};\\
&y_{(1\,l)}y_{12}y_{(1\,l)}\rightsquigarrow
-3\lambda_{112}\chi_{(3\,l)}(g_{12})y_{(3\,l)}y_2y_{(1\,l) }+\lambda_{112}
(\xi-\xi^2)\chi_{(2\,l)}(g_1)y_{ (2\,l)}y_{(1\,l)};\\
&
y_{12}y_{(1\,l)}y_{12}\rightsquigarrow  -3\xi^2
\lambda_{112}\chi_{(1\,l)}(g_2)y_{(3\,l)}y_2y_{12}+\lambda_{112}(1-\xi)y_{(2\,l)
}y_{12};\\
&y_{12}^2y_{(1\,l)}\rightsquigarrow  3
\lambda_{112}\xi\chi_{(2\,l)}(g_{122})y_{(3\,l)}y_2y_{12}
+\lambda_{112}
\chi_ {(2\,l)}(g_{12})(\xi^2-\xi)y_{(2\,l)}y_{12}.
\end{align*}
We have $X_2AA+AX_2A\rightsquigarrow  0$:
\begin{align*}
&X_2AA=\chi_{(1\,l)}(g_{(3\,l)})^2y_{12}y_{(1\,l)}^2g_{(3\,l)}\ot x_{(3\,l)}\\
& \rightsquigarrow  \lambda_{112}\Big(3
y_{(3\,l)}y_2y_{(1\,l)}g_{(3\,l)}
-(1-\xi)\chi_{2}(g_{(3\,l)})y_{(2\,l)}y_{
(1\,l)}g_{(3\,l)}\Big)\ot
x_{(3\,l)}\\
&AX_2A=\chi_{(1\,l)}(g_{(3\,l)})y_{(1\,l)}y_{12}y_{(1\,l)}g_{(3\,l)}\ot
x_{(3\,l)}\\
& \rightsquigarrow  \lambda_{112}\Big(-3y_{(3\,l)}y_2y_{(1\,l) }g_{(3\,l)}+
(1-\xi)\chi_{2}(g_{(3\,l)})y_{ (2\,l)}y_{(1\,l)}g_{(3\,l)}\Big)\ot x_{(3\,l)}
\end{align*}
Also, $X_2AX_2+X_2X_2A\rightsquigarrow  0$:
\begin{align*}
& X_2AX_2=
\chi_{12}(g_{(3\,l)})\chi_{(1\,l)}(g_{(3\,l)})y_{12}y_{(1\,l)}y_{12}g_{(3\,l)}^2\ot
x_{(3\,l)}^2\\
& \qquad \rightsquigarrow  \lambda_{112}\Big(-3\xi
y_{(3\,l)}y_2y_{12}g_{(3\,l)}^2+(\xi-\xi^2)\chi_2(g_{(3\,l)})y_{(2\,l)}y_{12}g_{
(3\,l) }^2\Big)\ot
x_{(3\,l)}^2 \\
 &X_2X_2A=
\chi_{12}(g_{(3\,l)})\chi_{(1\,l)}(g_{(3\,l)})^2y_{12}^2y_{(1\,l)}
g_{(3\,l)}^2\ot
 x_{(3\,l)}^2\\
 &\qquad \rightsquigarrow  \lambda_{112}\Big( 3
 \xi y_{(3\,l)}y_2y_{12}g_{(3\,l)}^2+
 \chi_2(g_{(3\,l)})(\xi^2-\xi)y_{(2\,l)}y_{12}g_{(3\,l)}^2\Big)\ot
 x_{(3\,l)}^2.
\end{align*}
We move onto (G2'.3). We need:
\begin{align*}
&y_{(1\,p)}y_{(1\,l)}^2\rightsquigarrow  
\chi_{(1\,l)}(g_{(1\,p)})y_{(1\,l)}y_{(1\,p)}y_{(1\,l)} -3\xi^2\lambda_{112}
\chi_{(1\,l)}(g_{(2\,p)})y_{(3\,l)}y_{(2\,p)}y_{(1\,l)}\\
&\qquad
+3\lambda_{112}(1-\xi^2)\chi_{(1\,l)}(g_{(2\,p)})\chi_{(3\,p)}(g_2)y_{(3\,l)}y_{(3\,p)}
y_2y_{(1\,l)}.\\
&y_{(1\,l)}y_{(1\,p)}y_{(1\,l)}\rightsquigarrow  -3\xi^2\lambda_{112}
\chi_{(3\,l)}(g_{(1\,l)})y_{(3\,l)}
y_ {(2\, p) }y_{(1\, l) }  \\
&\qquad
+3\lambda_{112}(1-\xi^2)\chi_{(3\,p)}(g_2)\chi_{(3\,l)}(g_{(1\,l)
} )y_{(3\,l)}y_{(3\,p)}
y_2y_{(1\,l)}. \\
&y_{(1\,p)}y_{(1\,l)}y_{(1\,p)}\rightsquigarrow   -3\xi^2\lambda_{112}
\chi_{(1\,l)}(g_{(2\,p)})y_{(3\,l)}y_{(2\,p)}y_{(1\,p)}\\
&\qquad
+3\lambda_{112}(1-\xi^2)\chi_{(1\,l)}(g_{(2\,p)})\chi_{(3\,p)}(g_2)y_{(3\,l)}y_{(3\,p)}
y_2y_{(1\,p)}.\\
&y_{(1\,p)}^2y_{(1\,l)}\rightsquigarrow \chi_{(1\,l)}(g_{(1\,p)})y_{(1\,p)}y_{(1\,l)}y_{(1\,p)}\\
&\qquad -3\xi^2\lambda_{112}
\chi_{(1\,l)}(g_{(2\,p)})\chi_{(3\,l)}(g_{(1\,p)})\chi_{(2\,p)}(g_{(1\,p)})y_{(3\,l)}
 y_{(2\,p)}y_{(1\,p)}\\
&\qquad
-3\lambda_{112}(1-\xi)\chi_{(1\,l)}(g_{(2\,p)})\chi_{(3\,l)}(g_{(1\,p)})\chi_{(3\,p)}
(g_{(1\,p)})y_{(3\,l)}y_{ (3\,p)}y_2y_{(1\,p)}.
\end{align*}
Now, $X_pAA+AX_pA$ equals
\begin{align*}
\chi_{(1\,l)}(g_{(p+1\,l)})\Big(\chi_{(1\,l)}(g_{(p+1\,l)})y_{(1\,p)}y_{(1\,l)}
^2+y_{(1\,l)}y_{(1\,p)}y_{(1\,l)}\Big)g_{(p+1\,l)}\ot
x_{(p+1\,l)}
\end{align*}
and the factor between brackets, according to the equalities above is
\begin{align*}
&\rightsquigarrow (1+\xi)y_{(1\,l)}y_{(1\,p)}y_{(1\,l)} -3\xi^2\lambda_{112}
\chi_{(1\,l)}(g_{(2\,l)})y_{(3\,l)}y_{(2\,p)}y_{(1\,l)}\\
&\qquad +3\lambda_{112}(1-\xi^2)\chi_{(1\,l)}(g_{(2\,l)})\chi_{(3\,p)}(g_2)y_{(3\,l)}
y_{(3\, p) }y_2y_{(1\,l)}\\
&\rightsquigarrow  -3(1+\xi^2)\lambda_{112}
\chi_{(3\,l)}(g_{(1\,l)})y_{(3\,l)}
y_ {(2\, p) }y_{(1\, l) }  \\
&\qquad +3\lambda_{112}(\xi-\xi^2)\chi_{(3\,p)}(g_2)\chi_{(3\,l)}(g_{(1\,l)
} )y_{(3\,l)}y_{(3\,p)}
y_2y_{(1\,l)} \\
&\qquad -3\xi^2\lambda_{112}
\chi_{(1\,l)}(g_{(2\,l)})y_{(3\,l)}y_{(2\,p)}y_{(1\,l)}\\
&\qquad +3\lambda_{112}(1-\xi^2)\chi_{(1\,l)}(g_{(2\,l)})\chi_{(3\,p)}(g_2)y_{(3\,l)}
y_{(3\, p) }y_2y_{(1\,l)}\\
&=-3\lambda_{112}\chi_{(1\,l)}(g_{(2\,l)})(\xi^2+(1+\xi^2)\chi_{(3\,l)}(g_1)\chi_1(g_{(2\,
l) } ))y_{ (3\, l) } y_{ (2\, p) } y_{ (1\, l) } \\
&\qquad
+3\lambda_{112}\chi_{(1\,l)}(g_{(2\,l)})\chi_{(3\,p)}(g_2)(1-\xi^2+(\xi-\xi^2)\chi_{(3\,
l)}(g_1)\chi_1(g_{(2\,
l) } ))\\
&\qquad \hspace*{8cm} \times y_{(3\,l)}y_{(3\, p) }y_2y_{(1\,l)}\\
&=-3\lambda_{112}\chi_{(1\,l)}(g_{(2\,l)})(1+\xi+\xi^2)y_{(3\,l)}y_{(2\,p)} y_{(1\,l)}\\
&\qquad
+3\lambda_{112}\chi_{(1\,l)}(g_{(2\,l)})\chi_{(3\,p)}(g_2)(1-\xi^2+(\xi-\xi^2)\xi)y_{(3\,
l)}
y_{(3\, p) }y_2y_{(1\,l)}=0.
\end{align*}
Here we use that $\lambda_{112}\chi_{12}^{-1}=\lambda_{112}\chi_1$ and
$\lambda_{112}\chi_1(g_2)=\xi$.
Analogously, if $\alpha=\chi_{(1\,p)}(g_{(p+1\,l)})\chi_{(1\,l)}(g_{(p+1\,l)})$, then
$X_pAX_p+X_pX_pA$ is
\begin{align*}
\alpha\Big(\chi_{(1\,l)}(g_{(p+1\,l)})
y_{(1\,p)}^2y_{(1\,l)}
+y_{(1\,p)}y_{(1\,l)}y_{(1\,p)}\Big)g_{(p+1\,l)}^2\ot
x_{(p+1\,l)}^2
\end{align*}
and the factor between brackets is now
\begin{align*}
&\rightsquigarrow -3(1+\xi^2)\lambda_{112}
\chi_{(1\,l)}(g_{(2\,p)})y_{(3\,l)}y_{(2\,p)}y_{(1\,p)}\\
&\qquad +3\lambda_{112}(\xi-\xi^2)\chi_{(1\,l)}(g_{(2\,p)})\chi_{(3\,p)}(g_2)y_{(3\,l)}y_{
(3\,  p) }
y_2y_{(1\,p)}\\
&\qquad -3\xi^2\lambda_{112}
\chi_{(1\,l)}(g_{(2\,l)})\chi_{(3\,l)}(g_{(1\,p)})\chi_{(2\,p)}(g_{(1\,p)})y_{(3\,l)}
 y_{(2\,p)}y_{(1\,p)}\\
&\qquad -3\lambda_{112}(1-\xi)\chi_{(1\,l)}(g_{(2\,l)})\chi_{(3\,l)}(g_{(1\,p)
})\chi_{(3\,p)}(g_{(1\,p)})y_{(3\,l)}y_{ (3\,p)}y_2y_{(1\,p)}\\
&=-3\lambda_{112}\chi_{(1\,l)}(g_{(2\,p)})(1+\xi+\xi^2)y_{(3\,l)}y_{(2\,p)}y_{(1\,p)}\\
&\qquad
+3\lambda_{112}\chi_{(1\,l)}(g_{(2\,p)})\chi_{(3\,p)}(g_2)(1-\xi)(\xi-\xi^4)y_{(3\,l)}
y_{(3\,p)}y_2y_{(1\,p)}=0.
\end{align*}
We are left with group (G3'). Again, we subdivide it:
\begin{enumerate}
 \item[(G3'.1)] Case $p=1$, $q=2$.
 \item[(G3'.2)] Case $p=1$, $q\geq 3$.
 \item[(G3'.3)] Case $p=2$, $q\geq 3$.
 \item[(G3'.4)] Case $p\geq 3$.
\end{enumerate}
For (G3'.1), we have
$$X_1X_2X_l+X_1X_lX_2+X_lX_1X_2+X_2X_lX_1+X_2X_1X_l={\bf
Y}_1g_{(2\,l)}g_{(3\,l)}\ot x_{(3\,l)}x_{(2\,l)},$$
\begin{align*}
{\bf Y}_1&=
\chi_{(1\,l)}(g_{(2\,l)}g_{(3\,l)})\chi_{12}(g_{(2\,l)})\chi_{(3\,l)}(g_{(2\,l)})y_{1}y_{
12 } y_ { (1\, l) }\\
&\qquad
+\chi_{(1\,l)}(g_{(2\,l)})\chi_{12}(g_{(2\,l)})\chi_{(3\,l)}(g_{(2\,l)})y_{1}y_{(1\,l) }
y_ {
12}
\\
&\qquad +\chi_{12}(g_{(2\,l)})\chi_{(3\,l)}(g_{(2\,l)})y_{(1\,l)}y_{1}y_{12}
+\chi_{(1\,l)}(g_{(3\,l)})\chi_{1}(g_{(3\,l)})y_{12}y_{(1\,l)}y_{1}\\
&\qquad
+\chi_{(1\,l)}(g_{(2\,l)}g_{(3\,l)})\chi_{1}(g_{(3\,l)})y_{12}y_{1}y_{(1\,l)}\\
&=\xi^2\chi_{1}(g_{(3\,l)})\chi_{112}(g_{(2\,l)})y_{1}
y_{12} y_ { (1\, l) }+\xi^2\chi_{112}(g_{(2\,l)})\chi_{(3\,l)}(g_2)y_{1}y_{(1\,l)}y_{
12}\\
&\qquad
+\xi\chi_{1}(g_{(2\,l)})y_{(1\,l)}y_{1}y_{12}+\xi\chi_{112}(g_{(3\,l)})y_{12}y_{(1\,l)}y_{
1}\\
&\qquad + \xi^2\chi_{1}(g_{(2\,l)})\chi_{112}(g_{(3\,l)})y_{12}y_{1}y_{
(1\,l)}.
\end{align*}

Hence we need:
\begin{align*}
 &y_{1}y_{12}y_{(1\,l)}\rightsquigarrow
\chi_{12}(g_1)y_{12}y_1y_{(1\,l)}+\lambda_{112}y_{(1\,l)};\\
 &y_{1}y_{(1\,l)}y_{12}\rightsquigarrow
\lambda_{112}\chi_{(1\,l)}(g_1)y_{(1\,l)}+\lambda_{112}(1-\xi^2)y_{(3\,l)}y_{12}
;\\
 &y_{(1\,l)}y_{1}y_{12}\rightsquigarrow  \lambda_{112}y_{(1\,l)};\\
 &y_{12}y_{(1\,l)}y_{1}\rightsquigarrow
-3\xi^2
\lambda_{112}\chi_{(1\,l)}(g_2)y_{(3\,l)}y_2y_{1}+\lambda_{112}(1-\xi)y_{(2\,l)}
y_{1};\\
 &y_{12}y_{1}y_{(1\,l)}\rightsquigarrow
-3\xi^2
\lambda_{112}\chi_{(1\,l)}(g_{12})y_{(3\,l)}y_2y_{1}+\lambda_{112}
(1-\xi^2)y_
{(1\,l)}\\
&\qquad
+\lambda_{112}(1-\xi)\chi_{(1\,l)}(g_1)y_{(2\,l)}y_{1}+\lambda_{112}(1-\xi^2)
\chi_{(3\,l)}(g_{12})y_ {(3\,l)}y_{12}.
\end{align*}
Using the above identities,
\begin{align*}
{\bf Y}_1&\rightsquigarrow  \Big(\xi^2\chi_{1}(g_{(3\,l)})\chi_{112}(g_{(2\,l)})
\chi_{12}(g_1)+\xi^2\chi_{1}(g_{(2\,l)})\chi_{112}(g_{(3\,l)})\Big)y_{12}y_1y_{(1\,l)}\\
&\qquad + \lambda_{112}(\xi^2-\xi)\chi_{1}(g_{(3\,l)})y_{(1\,l)}
+\lambda_{112}(\xi^2-\xi)\chi_{(3\,l)}(g_2)y_{(3\,l)}y_{12}\\
&\qquad
-3\xi^2\lambda_{112}\chi_{(3\,l)}(g_2)y_{(3\,l)}y_2y_{1}+\lambda_{112}(\xi-\xi^2)y_{(2\,l)
}y_{1}\\
&\rightsquigarrow  -3(1+\xi)
\lambda_{112}\chi_{(3\,l)}(g_{2})y_{(3\,l)}y_2y_{1}+\lambda_{112}
(\xi-\xi^2)\chi_{1}(g_{(3\,l)})y_
{(1\,l)}\\
&\qquad
+\lambda_{112}(\xi^2-\xi)y_{(2\,l)}y_{1}
+\lambda_ { 112}(\xi-\xi^2)
\chi_{(3\,l)}(g_{2})y_{(3\,l)}y_{12}\\
&\qquad + \lambda_{112}(\xi^2-\xi)\chi_{1}(g_{(3\,l)})y_{(1\,l)}
+\lambda_{112}(\xi^2-\xi)\chi_{(3\,l)}(g_2)y_{(3\,l)}y_{12}\\
&\qquad
-3\xi^2\lambda_{112}\chi_{(3\,l)}(g_2)y_{(3\,l)}y_2y_{1}+\lambda_{112}(\xi-\xi^2)y_{(2\,l)
}y_{1}\\
&=-3(1+\xi+\xi^2)\lambda_{112}\chi_{(3\,l)}(g_{2})y_{(3\,l)}y_2y_{1}\\
&\qquad +\lambda_{112}((\xi^2-\xi)+(\xi-\xi^2))\chi_{(3\,l)}(g_2)y_{(3\,l)}y_{12}\\
&\qquad +\lambda_{112}((\xi^2-\xi)+(\xi-\xi^2))y_{(2\,l)}y_{1}\\
&\qquad +\lambda_{112}\chi_{1}(g_{(3\,l)})((\xi-\xi^2)+(\xi^2-\xi))y_{(1\,l)}=0.
\end{align*}
Now we turn to (G3'.2): we have, for $q\geq 3$,
\begin{multline*}
X_1X_qX_l+X_1X_lX_q+X_lX_1X_q+X_qX_lX_1+X_qX_1X_l\\={\bf
Y}_2g_{(2\,l)}g_{(q+1\,l)}\ot x_{(q+1\,l)}x_{(2\,l)},
\end{multline*}
for
\begin{align*}
{\bf Y}_2&=\xi^2\chi_{1}(g_{(2\,l)})^2\chi_{(1\,l)}(g_{(q+1\,l)})y_1y_{(1\,q)}y_{(1\,l)}
+\xi^2\chi_{1}(g_{(2\,l)})^2y_1y_{(1\,l)}y_{(1\,q)}\\
&\qquad +\xi\chi_{1}(g_{(2\,l)})y_{(1\,l)}y_1y_{(1\,q)}
+\chi_{112}(g_{(q+1\,l)})\chi_{(3\,l)}(g_{(q+1\,l)})
y_{(1\,q)}y_{(1\,l)}y_1\\
&\qquad
+\xi\chi_{1}(g_{(2\,l)})\chi_{112}(g_{(q+1\,l)})\chi_{(3\,l)}(g_{(q+1\,l)})
y_{(1\,q)}y_1y_{(1\,l)}
\end{align*}
and thus we need
\begin{align*}
 &y_1y_{(1\,q)}y_{(1\,l)}\rightsquigarrow
\chi_{(1\,q)}(g_1)y_{(1\,q)}y_1y_{(1\,l)}+\lambda_{112}(1-\xi^2)y_{(3\,q)}y_{(1\,l)};\\
 &y_1y_{(1\,l)}y_{(1\,q)}\rightsquigarrow
\lambda_{112}(1-\xi^2)\chi_{(1\,l)}(g_1)\chi_{(3\,q)}(g_{(1\,l)})y_{(3\,q)}y_{(1\,l)}\\
&\qquad \qquad \qquad
+\lambda_{112}(1-\xi^2)y_{(3\,l)}y_{(1\,q)};\\
 &y_{(1\,l)}y_1y_{(1\,q)}\rightsquigarrow
\lambda_{112}(1-\xi^2)\chi_{(3\,q)}(g_{(1\,l)})y_{(3\,q)}y_{(1\,l)};\\
 &y_{(1\,q)}y_{(1\,l)}y_1\rightsquigarrow
-3\xi^2\lambda_{112}\chi_{(1\,l)}(g_{(2\,q)})
y_{(3\,l)}y_{(2\,q)}y_1
+3\lambda_{112}\chi_{1}(g_{(3\,q)})y_{(3\,q)}y_{(2\,l)}y_1;\\
 &y_{(1\,q)}y_1y_{(1\,l)}\rightsquigarrow  -3\xi^2\lambda_{112}\chi_{(1\,l)}(g_{(1\,q)})
y_{(3\,l)}y_{(2\,q)}y_1  -3\lambda_{112}\chi_{(3\,q)}(g_{12})y_{(3\,q)}y_{(1\,l)}\\
&\,
+3\lambda_{112}\chi_{(1\,l)}(g_{1})\chi_{1}(g_{(3\,q)})y_{(3\,q)}y_{(2\,l)}y_1
+\lambda_{112}(1-\xi^2)\chi_{(3\,l)}(g_{(1\,q)})y_{(3\,l)}y_{(1\,q)}.
\end{align*}
That is
\begin{align*}
{\bf Y}_2&\rightsquigarrow
-3\lambda_{112}(1+\xi+\xi^2)\chi_{(3\,l)}(g_{2})
\chi_{12}(g_{(3\,q)})
y_{(3\,l)}y_{(2\,q)}y_1\\
&\qquad +3\lambda_{112}(1+\xi+\xi^2)\xi\chi_{(3\,l)}(g_{(q+1\,l)})
\chi_{1}(g_{(3\,q)})y_{(3\,q)}y_{(2\,l)}y_1\\
&\qquad
+\lambda_{112}(\xi-\xi^2)\chi_{(3\,l)}(g_{2})(1-\chi_{112}(g_{(3\,l)}))
y_{(3\,l)}y_{(1\,q)}\\
&\qquad-3\lambda_{112}\chi_{(3\,q)}(g_{(1\,l)})\chi_{1}(g_{(2\,l)})(1+\xi+\xi^2)y_{(3\,q)}y_{(1\,l)}\rightsquigarrow  0.
\end{align*}
For (G3'.3), we have, for $q\geq 3$,
\begin{multline*}
X_2X_qX_l+X_2X_lX_q+X_lX_2X_q+X_qX_lX_2+X_qX_2X_l\\={\bf
Y}_3g_{(3\,l)}g_{(q+1\,l)}\ot x_{(q+1\,l)}x_{(3\,l)},
\end{multline*}
for
\begin{align*}
{\bf Y}_3&=\chi_{12}(g_{(3\,l)})^2\chi_{(1\,q)}(g_{q+1\,l})y_{12}y_{(1\,q)}y_{(1\,l)}
+\xi^2\chi_{12}(g_{(3\,l)})^2y_{12}y_{(1\,l)}y_{(1\,q)}\\
&\qquad
+\xi\chi_{12}(g_{(3\,l)})y_{(1\,l)}y_{12}y_{(1\,q)}+\xi\chi_{112}(g_{(q+1\,l)})\chi_{(2\,
q) } (g_{q+1\,l})y_{(1\,q)}y_{(1\,l)
}y_{12}\\
&\qquad
+\xi^2\chi_{12}(g_{(3\,l)})\chi_{(2\,q)}(g_{q+1\,l})\chi_{112}(g_{(q+1\,l)})y_{(1\,q)}
y_{12}y_{(1\,l)}.
\end{align*}
Hence we need:
\begin{align*}
 &y_{12}y_{(1\,q)}y_{(1\,l)}\rightsquigarrow  \chi_{(1\,q)}(g_{12})y_{(1\,q)}y_{12}y_{(1\,l)}
-3\lambda_{112}\xi\chi_{(3\,q)}(g_2)y_{(3\,q)}y_2y_{(1\,l)}\\
&\qquad +\lambda_{112}(1-\xi)y_{(2\,q) }y_{(1\,l)};\\
 &y_{12}y_{(1\,l)}y_{(1\,q)}\rightsquigarrow
\chi_{(1\,l)}(g_{12})y_{(1\,l)}y_{12}y_{(1\,q)}\\
&\qquad - 3\lambda_{112}\xi\chi_{(3\,l)}(g_2)
y_{(3\,l)}y_2y_{(1\,q)}\\
&\qquad +\lambda_{112}(1-\xi)
y_{(2\,l) }y_{(1\,q)};\\
&y_{(1\,l)}y_{12}y_{(1\,q)}\rightsquigarrow  -3\lambda_{112}\xi
\chi_{(3\,q)}(g_{122})\chi_{(3\,q)}(g_{(q+1\,l)})\chi_{2}(g_{(3\,l)})
y_{(3\,q)}y_2y_{(1\,l)}\\
&\qquad+\lambda_{112} (1-\xi)\chi_{ (2\,q)}(g_{(1\,l)})y_{ (2\,q)}y_{(1\,l)}; \\
 &y_{(1\,q)}y_{(1\,l)}y_{12}\rightsquigarrow  -3\lambda_{112}\xi^2\chi_{(1\,l)}(g_{(2\,q)})
y_{(3\,l)}y_{(2\,q)}y_{12} \\
&\qquad +3\lambda_{112}\chi_{1}(g_{(3\,q)})y_{(3\,q)}y_{(2\,l)}y_{12};
\end{align*}
\begin{align*}
 &y_{(1\,q)}y_{12}y_{(1\,l)}\rightsquigarrow
-3\lambda_{112}\xi^2\chi_{(1\,l)}(g_{(3\,q)})\chi_{(1\,l)}(g_{122})
y_{(3\,l)}y_{(2\,q)}y_{12} \\
&\qquad
+3\lambda_{112}\xi\chi_{1}(g_{(3\,q)})\chi_{(3\,l)}(g_{12})y_{(3\,q)}y_{(2\,l)}y_{12}\\
&\qquad -3\xi^2
\lambda_{112}\chi_{(1\,l)}(g_2)\chi_{(3\,l)}(g_{(1\,q)})\chi_2(g_{(1\,q)})y_{(3\,l)}y_2y_{
(1\,q)}\\
&\qquad+\lambda_{112}(1-\xi)\chi_{(2\,l)}(g_{(1\,q)})y_{(2\,l)}y_{(1\,q)}\\
&\qquad -3\xi^2
\lambda_{112}\chi_{(1\,l)}(g_2)\chi_{(3\,q)}(g_{(1\,q)})\chi_2(g_{(1\,l)})
(1-\xi^2)y_{(3\,q)} y_2y_{(1\,l)}\\
&\qquad+3\lambda_{112}\xi^2\chi_{(3\,q)}(g_1)y_{(2\,q)}y_{(1\,l)}.
\end{align*}
We have used that, by q-Jacobi \eqref{eq:qjacob} and Remark
\ref{rem:1l-2p} we have
\begin{align*}
[y_{(1\,q)},y_{(2\,l)}]_c&=[y_1,[y_{(1\,q)},y_{(2\,l)}]_c]_c
+\chi_{(2\,l)}(g_{(2\,q)})[y_{(1\,l)},y_{(2\,q)}]_c\\
&\qquad
-\chi_{(2\,q)}(g_1)(1-\xi)y_{(2\,q)}y_{(1\,l)}\\
&=[y_1,[y_{(1\,q)},y_{(2\,l)}]_c]_c
-3\lambda_{122}\xi^2\chi_{(2\,l)}(g_{(3\,q)})y_{(3\,l)}y_{(3\,q)}\\
&\qquad -\chi_{(2\,q)}
(g_1)(1-\xi)y_{(2\,q)}y_{(1\, l)}
\end{align*}
and thus, combining Lemma \ref{lem:l112} and
Lemma \ref{lem:tech3} \eqref{it7}:
\begin{multline}\label{eqn:1q-2l}
\lambda_{112}[y_{(1\,q)},y_{(2\,l)}]_c=
-3\lambda_{112}\lambda_{122}\xi^2\chi_{(2\,l)}(g_{(3\,q)})y_{(3\,l)}y_{(3\,q)}\\
 -\lambda_{112}\chi_{(2\,q)}
(g_1)(1-\xi)y_{(2\,q)}y_{(1\, l)}.
\end{multline}
Hence,
\begin{align*}
{\bf Y}_3&\rightsquigarrow
-3\lambda_{112}\chi_{122}(g_{(q+1\,l)})\chi_{(3\,q)}(g_{q+1\,l} )
(1+\xi+\xi^2)y_{(3\,q)}y_2y_{(1\,l)}\\
&\qquad+3\lambda_{112}(1+\xi+\xi^2)\chi_{2}(g_{(3\,l)})\chi_{(1\,q)}(g_{(q+1\,l)})
y_{(2\,q)}y_{(1\,l)} \\
&\qquad -3\lambda_{112}(1+\xi+\xi^2)y_{(3\,l)}y_2y_{(1\,q)}\\
&\qquad+\lambda_{112}((1-\xi^2)-(1-\xi^2))\chi_{2}(g_{(3\,l)})
y_{(2\,l)}y_{(1\,q)}\\
 &\qquad
-3\lambda_{112}(1+\xi+\xi^2)\chi_{1}(g_{(3\,q)})
y_{(3\,l)}y_{(2\,q)}y_{12} \\
&\qquad
+3\lambda_{112}(1+\xi+\xi^2)\chi_{1}(g_{(3\,q)})\chi_
{(2\,q)}(g_{q+1\,l})
y_{(3\,q)}y_{(2\,l)}y_{12}\rightsquigarrow 0.
\end{align*}
Finally, for (G3'.4), we have, for $3\leq p<q$,
 \begin{multline*}
 X_pX_qX_l+X_pX_lX_q+X_lX_pX_q+X_qX_lX_p+X_qX_pX_l\\={\bf
Y}_4g_{(p+1\,l)}g_{(q+1\,l)}\ot x_{(q+1\,l)}x_{(p+1\,l)},
                    \end{multline*}
\begin{align*}
{\bf Y}_4&=
\chi_{(1\,l)}(g_{(p+1\,l)})^2\chi_{(1\,l)}(g_{(q+1\,l)})y_{(1\,p)}y_{(1\,q)}y_{(1\,l)}\\
&\qquad
+\chi_{(1\,l)}(g_{(p+1\,l)})^2y_{(1\,p)}y_{(1\,l)}y_{(1\,q)}+\chi_{(1\,l)}(g_{(p+1\,l)})y_
{(1\,l)}y_{(1\,p)}y_{(1\,q)}\\
&\qquad
+\chi_{(1\,l)}(g_{(q+1\,l)})\chi_{(1\,p)}(g_{(q+1\,l)})y_{(1\,q)}y_{(1\,l)}y_{(1\,p)}\\
&\qquad +\chi_{(1\,l)}(g_{(p+1\,l)})\chi_{(1\,l)}(g_{(q+1\,l)})\chi_{(1\,p)}(g_{(q+1\,l)})
y_{(1\,q)}y_{(1\,p)}y_{(1\,l)}.
\end{align*}
Hence we need:
\begin{align*}
 &y_{(1\,p)}y_{(1\,q)}y_{(1\,l)}\rightsquigarrow
\chi_{(1\,q)}(g_{(1\,p)})y_{(1\,q)}y_{(1\,p)}y_{(1\,l)}\\
&\qquad -3\xi^2\lambda_{112}\chi_{(1\,q)}
(g_{(2\,p)})
y_{(3\,q)}y_{(2\,p)}y_{(1\,l)}
+3\lambda_{112}\chi_{1}(g_{(3\,p)})y_{(3\,p)}y_{(2\,q)}y_{(1\,l)};\\
 &y_{(1\,p)}y_{(1\,l)}y_{(1\,q)}\rightsquigarrow
\chi_{(1\,l)}(g_{(1\,p)})y_{(1\,l)}y_{(1\,p)}y_{(1\,q)} \\
&\qquad-3\xi^2\lambda_{112}\chi_{(1\,l)}
(g_{(2\,p)})
y_{(3\,l)}y_{(2\,p)}y_{(1\,q)}
+3\lambda_{112}\chi_{1}(g_{(3\,p)})
y_{(3\,p)}y_{(2\,l)}y_{(1\,q)};\\
 &y_{(1\,l)}y_{(1\,p)}y_{(1\,q)}\rightsquigarrow  -3\xi^2\lambda_{112}\chi_{(1\,q)}(g_{(2\,p)})
\chi_{(3\,q)}(g_{(1\,l)})\chi_{(2\,p)}(g_{(1\,l)})
y_{(3\,q)}y_{(2\,p)}y_{(1\,l)} \\
&\qquad
+3\lambda_{112}\chi_{1}(g_{(3\,p)})\chi_{(3\,p)}(g_{(1\,l)})\chi_{(2\,q)}(g_{(1\,l)})y_{
(3\,p)}y_{(2\,q)}y_{(1\,l)};
\\
 &y_{(1\,q)}y_{(1\,l)}y_{(1\,p)}\rightsquigarrow
-3\xi^2\lambda_{112}\chi_{(1\,l)}(g_{(2\,q)})
y_{(3\,l)}y_{(2\,q)}y_{(1\,p)} \\
&\qquad +3\lambda_{112}\chi_{1}(g_{(3\,q)})y_{(3\,q)}y_{(2\,l)}y_{(1\,p)};\\
 &y_{(1\,q)}y_{(1\,p)}y_{(1\,l)}\rightsquigarrow
-3\xi^2\lambda_{112}\chi_{(1\,l)}(g_{(2\,q)})
\chi_{(1\,l)}(g_{(1\,p)})y_{(3\,l)}y_{(2\,q)}y_{(1\,p)} \\
&\qquad
+3\lambda_{112}\chi_{1}(g_{(3\,q)})\chi_{(1\,l)}(g_{(1\,p)})y_{(3\,q)}y_{(2\,l)}y_{(1\,p)}
\\
&\qquad -3\xi^2\lambda_{112}\chi_{(1\,l)}
(g_{(2\,p)})
\chi_{(3\,l)}(g_{(1\,q)})\chi_{(2\,p)}(g_{(1\,q)})y_{(3\,l)}y_{(2\,p)}y_{(1\,q)} \\
&\qquad
-3\xi^2\lambda_{112}\chi_{(1\,l)}(g_{(2\,p)})
\chi_{(3\,q)}(g_{(1\,q)})(1-\xi^2)\chi_{(2\,p)}(g_{(1\,l)})y_{(3\,q)}y_{(2\,p)}y_{(1\,l)}
\\
&\qquad
+3\lambda_{112}\chi_{1}(g_{(3\,p)})\chi_{(3\,p)}(g_{(1\,q)})\chi_{(2\,l)}(g_{(1\,q)})
y_{(3\,p)}y_{(2\,l)}y_{(1\,q)}\\
&\qquad
-3\lambda_{112}\chi_{1}(g_{(3\,p)})\chi_{(3\,p)}(g_{(1\,q)})
\chi_{(2\,q)}
(g_1)(1-\xi)y_{(3\,p)}y_{(2\,q)}y_{(1\,l)}.
\end{align*}
Thus, we get to
\begin{align*}
{\bf
Y}_4&\rightsquigarrow -3\lambda_{112}\xi^2\chi_{(1\,l)}(g_{(p+1\,l)}^2g_{(q+1\,l)}
)\chi_{(1\,q)}
(g_{(2\,p)})(1+\xi+\xi^2)y_{(3\,q)}y_{(2\,p)}y_{(1\,l)} \\
&\quad
+3\lambda_{112}\chi_1(g_{(2\,p)})\chi_{(1\,l)}(g_{(q+1\,l)}g_{(p+1\,l)}^2)(1+\xi+\xi^2)y_{
(3\,p)}y_{(2\,q)}y_{(1\,l)} \\
&\quad -3\lambda_{112}\chi_{1}(g_{(2\,l)})
\chi_{(1\,l)}(g_{(p+1\,l)})(1+\xi+\xi^2)y_{(3\,l)}y_{(2\,p)}y_{(1\,q)} \\
&\quad +3\lambda_{112}\chi_{1}(g_{(3\,p)})
\chi_{(1\,l)}(g_{(p+1\,l)})^2(1+\xi+\xi^2)y_{(3\,p)}y_{(2\,l)}y_{(1\,q)} \\
&\quad
-3\lambda_{112}\chi_{1}(g_{(2\,l)})\chi_{(1\,p)}(g_{(q+1\,l)})(1+\xi+\xi^2)
y_{(3\,l)}y_{(2\,q)}y_{(1\,p)}\\
&\quad
+3\lambda_{112}\chi_{(1\,l)}(g_{(q+1\,l)})\chi_{(1\,p)}(g_{(q+1\,l)}
)\chi_1(g_{(3\,q)})(1+\xi+\xi^2)
y_{(3\,q)}y_{(2\,l)}y_{(1\,p)}\\
&\quad \rightsquigarrow 0,
\end{align*}
which establishes the lemma.
\epf

\begin{lem}\label{lem:i-kl3}
We have $[y_i,y_{(k\,l)}^3]_c=0$ for every $1\leq i\leq\theta$, $1\leq k\leq
l\leq\theta$.
\end{lem}
\pf
We show this by induction on $l-k$. If $l-k=0$, it is straightforward that
$[y_i,y_{k}^3]_c=0$ for $\mid k-i\mid>1$ by \eqref{eq:rels-Atheta-N3-galois}. This is also
clear if $k=i$. If $k=i+1$, say $i=1, k=2$, we get:
\begin{align*}
[y_1,y_2^3]_c&=y_{12}y_2^2+\chi_2(g_1)y_2y_{12}y_2+\chi_{2}(g_1)^2y_2^2y_{12}\\
&=\chi_2(g_{12})y_2y_{12}y_2+\lambda_{122}(1+\xi)y_2
+\chi_{2}(g_1)^2(1+\xi)y_2^2y_{12}\\
&=\lambda_{122}
(1+\xi+\xi^2)y_2
+\chi_{2}(g_1)^2(1+\xi+\xi^2)y_2^2y_{12}=0.
\end{align*}
Case $i=k+1$ is analogous. Fix $k$ and assume $[y_i,y_{(k\,p)}^3]_c=0$ for every $k\leq
p\leq l$, every $1\leq k\leq l\leq\theta$. It follows from Proposition
\ref{pro:right-coaction} that
\begin{multline*}
 \rho([y_i,y_{(k\,l)}^3]_c)=[y_i,y_{(k\,l)}^3]_c \ot
1+g_ig_{(k\,l)}^3\ot [x_i,x_{(k\,l)}^3]_c\\
+\sum_{k\leq p<l}C_p [y_i,y_{(k\,p)}^3]_cg_{(p+1\,l)}^3\ot
x_{(p+1\,l)}^3\\
+\sum_{k\leq p<l}C_p\gamma_p^3 y_{(k\,p)}^3g_ig_{(p+1\,l)}^3\ot
[x_i,x_{(p+1\,l)}^3]_c,
 \end{multline*}
for $C_p$ as in \eqref{eqn:Cp}  and
$\gamma_p=\chi_{(k\,p)}(g_i)$.
By induction, we have $[y_i,y_{(k\,p)}^3]_c=0$ for every $k\leq p<l$ while
$[x_i,x_{(k\,l)}^3]_c=[x_i,x_{(p+1\,l)}^3]_c=0$ for every $k\leq p<l$ by
\cite[Proposition 4.1]{Ang}. That is,
\begin{equation*}
 \rho([y_i,y_{(k\,l)}^3]_c)=[y_i,y_{(k\,l)}^3]_c \ot 1,
\end{equation*}
{\it i.e.}  $[y_i,y_{(k\,l)}^3]_c\in \mtA^{\co\mtH}=\k$. Set
$\k\ni s:=[y_i,y_{(k\,l)}^3]_c$. Now,
\begin{equation*}
s=g_i[y_i,y_{(k\,l)}^3]_cg_i^{-1}=\xi\chi_{(k\,l)}(g_i)^3[y_i,y_{(k\,l)}^3]_c.
\end{equation*}
Hence $s=0$ if $\chi_{(k\,l)}(g_i)^3\neq\xi^2$. On the other hand,
\begin{equation*}
s=g_{(k\,l)}^3[y_i,y_{(k\,l)}^3]_cg_{(k\,l)}^{-3}=\chi_i(g_{(k\,l)})^3
[y_i,y_{(k\,l)}^3] _c
\end{equation*}
and thus $s=0$ if $\chi_i(g_{(k\,l)})^3\neq 1$. But we cannot have both
$\chi_i(g_{(k\,l)})^3= 1$ and $\chi_{(k\,l)}(g_i)^3=\xi^2$ as it
contradicts $1=\chi_{(k\,l)}(g_i)^3\chi_i(g_{(k\,l)})^3$. Therefore $s=0$ and the
lemma follows.
\epf

The following shows \eqref{eqn:condition-recursive} for $j=1$ and thus \eqref{eqn:condition} in general.
\begin{theorem}\label{thm:A-N=3}
Let $\mA=\mA(\boldsymbol\lambda,\boldsymbol\mu)$ be
the algebra quotient of $\mT(V)$ by
\begin{align}\label{eq:rels-Atheta-N3-galois-b}
y_{ij} &= 0, \ i < j - 1\in\I; & y_{iij} &= \lambda_{iij}, \ i,j\in\I, \vert j - i\vert = 1;\\
y_{(k\,l)}^3 &= \mu_{(k\,l)}, \  k\leq l\in\I.
\end{align}
for families of scalars $\boldsymbol\lambda=(\lambda_{iij})_{i,j}$ and
$\boldsymbol\mu=(\mu_{(k\,l)})_{k,l}$ satisfying \eqref{eqn:restriction-1} and
\begin{align}\label{eqn:restriction-2}
 \mu_{(k\,l)}=0 \quad \text{ if } \chi_{(k\,l)}^3\neq\eps.
\end{align}
Then $\mA\in\Cleft\mH$. In particular,
\begin{align*}
\Cleft'\mH=\{\mA(\boldsymbol\lambda,\boldsymbol\mu)|\boldsymbol\lambda\text{ as in }\eqref{eqn:restriction-1},\boldsymbol\mu \text{ as in }\eqref{eqn:restriction-2}\}.
\end{align*}
\end{theorem}
\pf
By \cite{Ang}, $X=\,^{\co \mH}\mtH$ is the polynomial algebra in the variables
$$\text{x}_{kl}:=g_{(k\,l)}^{-3}x_{(k\,l)}^3, \qquad 1\leq k\leq l\leq\theta.$$ We will
show that the $\mtH$-colinear algebra maps $f:X\to \mtA$ generated by
$$\text{x}_{kl}\mapsto \text{y}_{kl}-\text{g}_{kl},$$ for
$\text{y}_{kl}:=g_{(k\,l)}^{-3}y_{(k\,l)}^3$ and
$\text{g}_{kl}:=\mu_{(k\,l)}g_{(k\,l)}^{-3}$  are also
$\mtH$-linear, when we consider the right adjoint action $\cdot:X\ot \mtH\to X$ and the
Miyashita-Ulbrich action $\leftharpoonup:\mtA\ot \mtH\to \mtA$. We have, for $h\in H$:
\begin{align*}
 f(\text{x}_{kl}\cdot
h)&=\chi_{(k\,l)}(h)^3(\text{y}_{kl}-\text{g}_{kl})=f(\text{x}_{kl})\leftharpoonup h,
\end{align*}
as $\chi_{(k\,l)}(g_i)^3\mu_{(k\,l)}=\mu_{(k\,l)}$ by \eqref{eqn:restriction-2}. Also, by
\cite[Proposition 4.1]{Ang}:
\begin{align*}
 \text{x}_{kl}\cdot
x_i&=(1-\chi_{(k\,l)}(g_i)^3\chi_i(g_{(k\,l)})^3)\text{x}_{kl}x_i=0.
\end{align*}
On the other hand, $(\text{y}_{kl}-\text{g}_{kl})\leftharpoonup
x_i=0$, by Lemma \ref{lem:i-kl3}. Then, the theorem follows by \cite[Theorem 8]{G}, see also \cite[Theorem 3.3]{AAnGMV}.
\epf

\subsection{Liftings}\label{subsec:liftings-3}

Let $\mtL=\mtL(\boldsymbol\lambda)$ be the quotient of $T(V)\# H$ by the relations
\begin{align}\label{eq:rels-Ltheta-N3-lifting}
a_{ij} &= 0, \ i < j - 1; & a_{iij} &= \lambda_{iij}(1-g_{iij}), \ \vert j - i\vert = 1
\end{align}
for some family of scalars $\boldsymbol\lambda=(\lambda_{iij})$ satisfying
\eqref{eqn:restriction-1} and normalized by
\begin{align}\label{eqn:normalization-1}
 \lambda_{iij}=0 \quad \text{ if } g_{iij}=1.
\end{align}
Here we rename the basis $\{x_1,\dots,x_\theta\}$ of $V$ by $\{a_1,\dots,a_\theta\}$.
\begin{rem}
Observe that normalization \eqref{eqn:normalization-1} is not necessary when $\theta\geq 3$. Take, for simplicity, $i=1, j=2$. Then we have
$ \xi^2=\chi_{112}(g_3)\chi_3(g_{112})$.
\end{rem}

Recall the definition of the distinguished pre-Nichols algebra $\tobat(V)$
from p. \pageref{eq:rels-Atheta-c} and the
bosonization $\mtH=\tobat(V)\# H$. 
\begin{pro}\label{pro:lift-pre-N=3}
Let $\boldsymbol\lambda=(\lambda_{iij})$ satisfy
\eqref{eqn:restriction-1} and \eqref{eqn:normalization-1}. Then
\begin{enumerate}
\item 
$\mtL(\boldsymbol\lambda)=L(\mtA(\boldsymbol\lambda),
\mtA(\boldsymbol\lambda))$; hence
$\mtL(\boldsymbol\lambda)$ is a cocycle deformation of $\mtH$.
 \item $\mtL(\boldsymbol\lambda)$ is a pointed Hopf algebra with
$\gr \mtL=\mtH$. 
\end{enumerate}
\end{pro}
\pf
Follows directly from Proposition \ref{pro:summary} (c); see also the case $N=2$ in Proposition 
\ref{pro:lift-pre-N=2}.
\epf

Let $\up:\mtA(\boldsymbol\lambda)\to
\mtL(\boldsymbol\lambda)\ot \mtA(\boldsymbol\lambda)$ be the coaction. We have
\begin{align*}
\up (y_{(kl)})&=a_{(kl)}\ot 1+g_{(kl)}\ot
y_{(kl)}+(1-\xi^2)\sum_{k\leq p<l}a_{(kp)}g_{(p+1l)}\ot y_{(p+1l)}.
\end{align*}

We proceed to describe the algebra $L(\mA(\boldsymbol\lambda,\boldsymbol\mu),\mH)$. As $\mH$ is obtained from $\mtH$ as a quotient of an ideal not
generated exclusively by (skew-)primitive elements, we thus need to prepare the 
setting
accordingly, cf. (the proof of) Proposition \ref{pro:summary} (a).

For each $m\geq 1$, consider the $m$-adic approximation $\widehat{\B}_m(V)$ to $\B(V)$. This is the quotient of $T(V)$ by relations
\eqref{eq:rels-Atheta-a} and \eqref{eq:rels-Atheta-b3} together with
\begin{align}\label{recursive-2}
 & x_{(k\,l)}^3, \qquad 1\leq l-k< m.
\end{align}
Thus, we obtain a family of cleft objects $\mA_m(\boldsymbol\lambda,\boldsymbol\mu_m)$ for
$\mH_m=\widehat{\B}_m(V)\# H$ given by the  quotient of $\mT(V)$
by relations \eqref{eq:rels-Atheta-N3-galois} for each together with
\begin{align*}
 & y_{(k\,l)}^3-\mu_{(k\,l)}, \qquad 1\leq l-k< m.
\end{align*}
Here $\boldsymbol\lambda=(\lambda_{iij})_{i,j}$ satisfies \eqref{eqn:restriction-1}
and $\boldsymbol\mu_m=(\mu_{(k\,l)})_{k\leq l}$ satisfies \eqref{eqn:restriction-2}.

Now, fix $\boldsymbol\lambda,\boldsymbol\mu_m$ and set
$\mA_m=\mA_m(\boldsymbol\lambda,\boldsymbol\mu_m)$. Let $\mL_m(\boldsymbol\lambda,\boldsymbol\mu_m):=L(\mA_m,\mH_m)$. Notice that $\mL_0=\mtL$.
We keep the name $\up :\mA_m\to \mL_m\ot \mA_m$ for the coaction at each level.
Thus $\mH_{m+1}=\mH_m/I_{m+1}$ is such that $I_{m+1}$ is generated by skew primitive elements \cite[Remark 6.10]{AS2}. Hence, by Proposition
\ref{pro:summary}, $\mL_{m+1}$ is the quotient of $\mL_m$ by the ideal generated by
\begin{align}\label{eqn:z=u}
 a_{(k\,l)}^3-\sigma_{(k\,l)}-\mu_{(k\,l)}(1-g_{(k,l)}^3),
\end{align}
where, according to \eqref{eqn:u-tilde}, the \emph{deforming elements} $\sigma_{(k\,l)}$ are defined by:
\begin{align}\label{eqn:dl}
\begin{split}
a_{(k\,l)}^3\ot 1-\up (y_{(k\,l)})^3=\sigma_{(k\,l)}\ot 1.
\end{split}
\end{align}
We give a description of these elements in Proposition \ref{pro:uil}.

In this way we obtain a description of the full algebra $\mL=L(\mA(\boldsymbol\lambda,\boldsymbol\mu),\mH)$ in the final step of this procedure.
We further normalize $\boldsymbol\mu$ by
\begin{align}\label{eqn:normalization-2}
 \mu_{(k\,l)}=0 \quad \text{ if } g_{(k\,l)}^3=1.
\end{align}

We illustrate this situation in the following two examples.
\begin{exa}\label{exa:ABC-lift-1}
$\mL_1$ is the quotient of $\mT(V)$ by relations \eqref{eq:rels-Ltheta-N3-lifting} and
$$
a_k^3=\mu_{(k)}(1-g_k^3).
$$
In particular, $\sigma_{(k\,k)}=0$, $k\in\I$.
\pf
Let $m=0$. The elements $y_k^3$, generating $I_1$, satisfy:
$$
 \up (y_k)^3=a_k^3\ot 1 + g_k^3\ot y_k^3.
$$
Hence $\uvi_{(k)}=0$ and the statement follows.
\epf
\end{exa}
The following example contains the spirit of our computations ahead.
\begin{exa}\label{exa:ABC-lift-2}
$\mL_2$ is the quotient of $\mT(V)$ by relations \eqref{eq:rels-Ltheta-N3-lifting} and
\begin{align*}
a_k^3&=\mu_{(k)}(1-g_k^3), \quad 1\leq k\leq \theta;\\
a_{kk+1}^3&=\mu_{(kk+1)}(1-g_k^3g_{k+1}^3)-\mu_{(k+1)}\mu_{(k)}(1-\xi)^3\chi_k(g_{k+1})^3 (1-g_k^3)g_{k+1}^3\\
&\qquad -\lambda_{kk+1k+1}\lambda_{kkk+1}\xi^2
(1-g_k^2g_{k+1})g_kg_{k+1}^2 \quad 1\leq k< \theta.
\end{align*}
\pf
Set $m=1$. We have already described $\mL_1$ in Example \ref{exa:ABC-lift-1}. We need to compute the elements $\uvi_{(k,k+1)}$, $k<\theta$.
It will be enough to understand $\up (y_{12})^3$. Set
\begin{align*}
&A=a_{12}\ot 1, && B=g_{12}\ot y_{12}, && C=a_1g_2\ot y_2,
\end{align*}
so that $\up (y_{12})=A+B+(1-\xi^2)C$. As before, we focus on the terms in which
a factor $\lambda_{***}$ may appear. These are related with two possible facts:
\begin{enumerate}
 \item Fact {\bf A}: $a_1$
appears to the left of $a_{12}$, that is:
\begin{align*}
 &CAB, && CBA, && BCA,  &CAA, && ACA,  &&CAC, && CCA.
\end{align*}
 \item Fact {\bf B}: $y_{12}$ appears to the left of $y_{2}$, that is:
\begin{align*}
 &ABC, && BAC, && BCA,  &BCB, && BBC,  &&BCC, && CBC.
\end{align*}
\end{enumerate}
We have
\begin{align*}
CAB&= a_1g_2a_{12}g_{12}\ot y_2y_{12}=\chi_{12}(g_2)a_1a_{12}g_{122}\ot
y_2y_{12}\\
&\rightsquigarrow  \lambda_{112}\xi^2(1-g_{122})g_{122}\ot y_2y_{12};\\
CBA & \rightsquigarrow  \lambda_{112}(1-g_{122})g_{122}\ot y_2y_{12};\\
BCA & =\chi_1(g_2)
a_1a_{12}g_{122}\ot y_2y_{12}+\lambda_{122}\xi^2 a_1a_{12}g_{122}\ot 1\\
&\rightsquigarrow  \lambda_{112}\xi
(1-g_{112})g_{122}\ot y_2y_{12}+\lambda_{122}\xi a_{12}a_1g_{122}\ot
1\\
&\quad  +\lambda_{122}\lambda_{112}\xi^2 (1-g_{112})g_{122}\ot 1;\\
CAA&=\chi_{12}(g_1)^2 a_1a_{12}^2g_2\ot y_2\rightsquigarrow
\lambda_{112}(1+\xi)a_{12}(1-g_{112})\ot y_2;\\
ACA&=\chi_{12}(g_2)a_{12}a_1a_{12}g_2\ot y_2\rightsquigarrow
\lambda_{112}\xi^2a_{12}(1-g_{112})\ot y_2;\\
CAC&=\chi_{112}(g_2)a_1a_{12}a_1g_2^2\ot y_2^2\rightsquigarrow  \lambda_{112}a_1(1-g_{112})g_2^2\ot
y_2^2;\\
CCA&=\chi_{112}(g_2)\chi_{12}(g_2)a_1^2a_{12}g_2^2\ot y_2^2\rightsquigarrow
\lambda_{112}(\xi+\xi^2)a_1(1-g_{112})g_2^2\ot y_2^2.
\end{align*}
On the other hand, we get
\begin{align*}
 ABC&=\chi_1(g_{12})a_{12}a_1g_{122}\ot y_{12}y_2\rightsquigarrow
\lambda_{122}\xi^2a_{12}a_1g_{122}\ot 1;\\
BAC&\rightsquigarrow  \lambda_{122} a_{12}a_1g_{122}\ot 1;\\
BCB&= \chi_1(g_{12})a_1g_{122}g_{12}\ot y_{12}y_2y_{12}\rightsquigarrow
\lambda_{122}\xi^2 a_1g_{122}g_{12}\ot y_2;\\
BBC&=\chi_1(g_{12})^2 a_1g_{122}g_{12}\ot y_{12}^2y_2\rightsquigarrow
\lambda_{122}(1+\xi) a_1g_{122}g_{12}\ot y_{12};\\
BCC&=\chi_1(g_{122})\chi_1(g_{12})a_1^2g_{122}g_2\ot y_{12}y_2^2
\rightsquigarrow  \lambda_{122}(\xi+\xi^2)a_1^2g_{122}g_2\ot y_2;\\
CBC&=\lambda_{122}a_1^2g_{122}g_2 \ot y_2.
\end{align*}
Hence,
\begin{align*}
 CAB+CBA+BCA\rightsquigarrow \lambda_{122}\xi a_{12}a_1g_{122}\ot
1  +\lambda_{122}\lambda_{112}\xi^2 (1-g_{112})g_{122}\ot 1
\end{align*}
and $CAA+ACA\rightsquigarrow 0$, $CAC+CCA\rightsquigarrow  0$. Also, we have
\begin{align*}
 ABC+BAC&\rightsquigarrow
(1+\xi^2)\lambda_{122}a_{12}a_1g_{122}\ot 1;\\
BCB+BBC&\rightsquigarrow  0, \qquad BCC+CBC\rightsquigarrow  0.
\end{align*}
Therefore,
\begin{align*}
 \up (y_{12})^3&=a_{12}^3\ot 1+g_{12}^3\ot y_{12}^3+(1-\xi)^3\chi_1(g_2)^3 a_1^3g_2^3\ot
y_2^3\\
& \hspace*{4cm }+\lambda_{122}\lambda_{112}\xi^2 (1-g_{112})g_{122}\ot 1\\
&=a_{12}^3\ot 1+g_{12}^3\ot y_{12}^3+\mu_{(2)}\mu_{(1)}(1-\xi)^3\chi_1(g_2)^3 (1-g_1^3)g_2^3\ot 1\\
& \hspace*{4cm }+\lambda_{122}\lambda_{112}\xi^2 (1-g_{112})g_{122}\ot 1.
\end{align*}
In particular, as $\mu_1\chi_1(g_2)^3=1$,
\begin{align*}
 \uvi_{(1,2)}=-\mu_{(2)}\mu_{(1)}(1-\xi)^3 (1-g_1^3)g_2^3-\xi^2\lambda_{122}\lambda_{112} (1-g_1^2g_2)g_1g_2^2.
\end{align*}
The statement follows.
\epf
\end{exa}

\begin{rem}\label{rem:BDR}
When $\theta=2$, then $\mL_2$ as in Example \ref{exa:ABC-lift-2} is a lifting of type $A_2$. It coincides with the liftings found in \cite{BDR} for this type.
\end{rem}

\subsection{The deforming elements}\label{sec:uil}

The expressions for both $\sigma_{(i)}:=\sigma_{(i\,i)}$ and 
$\sigma_{(i\,i+1)}$ follow from Examples \ref{exa:ABC-lift-1} and 
\ref{exa:ABC-lift-2}. Namely,
\begin{align*}
\sigma_{(i)}&=0,\quad i\in\I,\\
\sigma_{(i\,i+1)}&=-\mu_{(i+1)}\mu_{(i)}(1-\xi)^3\chi_i(g_{i+1})^3 (1-g_i^3)g_{i+1}^3\\
&\qquad\qquad -\lambda_{ii+1i+1}\lambda_{iii+1}\xi^2
(1-g_i^2g_{i+1})g_ig_{i+1}^2, \quad i< \theta\in\I.
\end{align*}
For the general case of $\sigma_{(i\,l)}$, $i,l \in\I$, we proceed in a similar fashion. 

We first define $\uvi_{(i\,l)}(\bs\mu)$ and $h_{il}(\bs\lambda)$ in $\k \Gamma$. 
We set $\uvi_{(ii)}=0$, and, recursively,
\begin{align}\label{eqn:def-gral-1}
\uvi_{(i\,l)}(\bs\mu)&=-\sum_{i\leq p<l} C_p\mu_{(p+1\,l)}
\Big(\uvi_{(i\,p)}+\mu_{(i\,p)}(1-g_{(i\,p)}^N)\Big)g_{(p+1\,l)}^N.
\end{align}
Now, we set $h_{ii}(\bs\lambda)=h_{ii+1}(\bs\lambda)=0$  and, for $l\geq i+2$,
\begin{align}\label{eqn:ul-2} h_{il}(\bs\lambda)&=-9\mu_{(i+2\,l)}\lambda_{ii+1i+1}\lambda_{iii+1}(1-g_{iii+1})g_{ii+1i+1}g_{(i+2\,l)}^3.
\end{align}

Next, for $i\leq p<l$, we set $q=p+1$, $r=p+2$ and consider the following elements in $T(V)\# H$:
\begin{multline}\label{eqn:monomials-1}
\varsigma^{p}_i(\bs\lambda,\bs\mu)=\lambda_{qrr}\Big(
\xi^2a_{(i\,p)}a_{(i\,q)}a_{(i\,r)}
+\chi_{r}(g_{(1\,p)})
a_{(i\,p)}a_{(i\,r)}a_{(i\,q)}\\
+a_{(i\,r)}a_{(i\,p)}a_{(i\,q)}\Big).
\end{multline}
Let us fix $s_p=-3(1-\xi^2)$, $p<l-2$, $s_{l-2}=1$, and set
$$
d_{i\,l}(p)=\chi_{(i\,q)}(g_{(q\,l)}g_{(r+1\,l)})
\chi_{(i\,p)}(g_{(r+1\,l)})s_p.
$$
Finally, we consider:
\begin{align}\label{eqn:monomials-intro-1}
\varsigma_{il}(\bs\lambda,\bs\mu)&=-3\xi^2\sum_{i\leq p<l}\mu_{(p+3\,l)}\chi_{r}(g_{(p+3\,l)})
d_{il}(p)\varsigma^p_i(\bs\lambda,\bs\mu)g_{qrr}g_{(p+3\,l)}^3.
\end{align}
Recall that $g_{(l+1\,l)}=1$; also we set $\mu_{(l+1\,l)}:=1$. 
\begin{rem}
Observe that nor $\varsigma^{p}(\bs\lambda,\bs\mu)$ neither $\varsigma_{il}(\bs\lambda,\bs\mu)$ are expressed in the PBW basis. This is an arduous computation that we perform in full generality in \S \ref{sec:uil}, see Corollary \ref{cor:suma-a} for a complete answer.
\end{rem}

\begin{pro}\label{pro:uil}
Let $i,l$ be as above. Then 
\begin{align}\label{eqn:uil}
\sigma_{(i\,l)}(\bs\lambda,\bs\mu)=\uvi_{(i\,l)}(\bs\mu)+h_{il}(\bs\lambda)+\varsigma_{il}(\bs\lambda,\bs\mu).
\end{align}
\end{pro}

See below for a proof. As a result, we have the following.

\begin{theorem}\label{thm:L(A,H)-N=3}
The Hopf algebra
$L(\mA(\boldsymbol\lambda,\boldsymbol\mu),\mH):=\mL(\boldsymbol\lambda,\boldsymbol\mu)$ is
the quotient of $\mT(V)$ by relations \eqref{eq:rels-Ltheta-N3-lifting} and
\begin{equation*}
 a_{(i\,l)}^3=\mu_{(i\,l)}(1-g_{(i\,l)}^3)+\sigma_{(i\,l)}(\boldsymbol\lambda,\boldsymbol\mu),
\end{equation*}
for $\sigma_{(i\,l)}(\boldsymbol\lambda,\boldsymbol\mu)$ as in \eqref{eqn:uil}.
\end{theorem}
\pf
Follows by Proposition \ref{pro:summary} (c), see \eqref{eqn:z=u}.
\epf

See Example \ref{ej:not-in-H-N3} below for a lifting of a concrete $V$. 
Next, we prove Proposition \ref{pro:uil}:
\pf
We take $i=1$ to ease up the notation, so $l\geq 3$. Set
\begin{align*}
A=a_{(1\,l)}\ot 1, && B=g_{(1\,l)}\ot y_{(1\,l)}, && X_p=a_{(1\,p)}g_{(p+1\,l)}\ot
y_{(p+1\,l)},
1\leq p<l
\end{align*}
so $\up (y_{(1\,l)})=A+B+(1-\xi^2)\sum_{1\leq p<l} X_p$.
We will also denote $X_l:=A$, by identifying as usual $
 g_{(l+1\,l)}:=1$, $y_{(l+1\,l)}:=1$. Finally, set
\begin{align}\label{eqn:dl-1}
 &\sigma_{p}:= \sigma_{(1\,p)}, \qquad 1\leq p\leq l.
\end{align}
As in Example \ref{exa:ABC-lift-2}, we need
to focus on the terms of $(A+B+(1-\xi^2)\sum X_p)^3$
involving a factor $\lambda_{***}$, as by \cite[Remark
6.10]{AS2} we have, for $C_p$ as in \eqref{eqn:Cp}:
\begin{multline}\label{eqn:plus-terms-with-lambda}
\up (y_{(1\,l)}^3)=a_{(1\,l)}^3\ot 1+ g_{(1\,l)}^3\ot y_{(1\,l)}^3 \\
+\sum_{1\leq p<l}C_{p}a_{(1\,p)}^3g_{(p+1\,l)}^3\ot y_{(p+1\,l)}^3\\
+\text{terms involving a factor }\lambda_{***}.
\end{multline}
Combining this 
with the recursive deformation procedure following
\cite[Corollary 5.12]{AAnGMV}, {\it i.e.} we assume $y_{(p+1\,l)}^3=\mu_{(p+1\,l)}$, we obtain
\begin{align}\label{eqn:dl-paso1}
 \sigma_l&=\eqref{eqn:def-gral-1}- \text{terms involving a factor }\lambda_{***}.
\end{align}

As in Proposition
\ref{pro:right-coaction},
we consider the cases (here we need to distinguish a factor $A$ from a factor $X_p$, identified previously):
\begin{enumerate}
\item[(L1)] For every $p<q$, terms $XYZ$ involving $X,Y,Z\in\{B,X_p,X_q\}$, all
different, $X_p$ to the left of $X_q$.
\item[(L2)] For every pair $p<q$, terms $XYZ$ involving $X,Y,Z\in \{X_p,X_q\}$, not all
equal and with a factor $X_p$ to the left of $X_q$.
\item[(L3)] For every triple $p<q<r$, terms $XYZ$ involving 
\emph{distinct} $X,Y,Z\in
\{X_p,X_q,X_r\}$ and with $X_p$ to the left of $X_{q}$ or $X_r$ or with $X_q$ to the left
of $X_r$.
\end{enumerate}
However, as Example \ref{exa:ABC-lift-2} illustrates, we also need to consider:
\begin{enumerate}
\item[(L4)] Terms $ABX_p$ and $BAX_p$, $1\leq p<l$.
\item[(L5)] Terms $BX_pB$ and $BBX_p$, $1\leq p<l$.
\item[(L6)] Terms $BX_qX_p$ and $X_qBX_p$, $1\leq p\leq q<l$.
\end{enumerate}
\begin{rem}\label{rem:p+3}
In cases (L1) and (L2) it is enough to consider $q<l$, as a factor $X_l=A$ will not contribute to $\sigma_l$. Case (L3) is different, and we will
take this difference into account: the main difference lays in the fact that the factors $X_l$ -unlike $X_p$,
$p<l$- are not multiplied by $(1-\xi^2)$. Hence commutativity computations 
follow rather smoothly, 
and we only have to recall this in the final expression.
\end{rem}

\begin{claim}\label{claim:1}
Cases (L4), (L5) and  (L6) do not contribute to  $\sigma_l$.
\end{claim}
These cases easily follow from Lemmas
\ref{lem:tech1} and \ref{lem:tech2b}. In (L4) we have $BAX_p=\xi ABX_p$ and:
\begin{align*}
BAX_p+ABX_p\rightsquigarrow
\begin{cases}
 0, &p>1;\\
3\xi^2\lambda_{122}\chi_{12}(g_{(1\,l)})a_{(1\,l)}a_{1}g_{(1\,l)}g_{(2\,l)
} \ot
y_{(3\,l)}^2, &p=1.
\end{cases}
\end{align*}
In (L5) we get:
\begin{align*}
 BX_pB\rightsquigarrow
\begin{cases}
 0, &p>1;\\
-3\lambda_{122}\chi_{12}( g_{(1\,l)})a_{1}g_{(2\,l)}g_{(1\,l)}^2\ot
y_{(3\,l)}^2y_{(1\,l)}, &p=1.
\end{cases}
\end{align*}
\begin{align*}
 BBX_p\rightsquigarrow
\begin{cases}
 0, &p>1;\\
3\lambda_{122}\chi_{12}( g_{(1\,l)})a_{1}g_{(2\,l)}g_{(1\,l)}^2\ot
y_{(3\,l)}^2y_{(1\,l)}, &p=1.
\end{cases}
\end{align*}
In particular, $BX_pB+BBX_p\rightsquigarrow 0$. For (L6), as $q+1\geq 3$,
\begin{align*}
BX_qX_p\rightsquigarrow  \begin{cases}
                0, &p>1;\\
-3\xi\lambda_{122}\chi_{12}(g_{(1\,l)})\chi_{1}(g_{(q+1\,l)})
a_{(1\,q)}a_{1}
g_{(1\,l)}g_{(q+1\,l)}g_{ (2\,l)}\\
\hspace*{6.6cm}\ot y_{(q+1\,l)}y_{(3\,l)}^2&p=1.
\end{cases}
\end{align*}
\begin{align*}
X_qBX_p\rightsquigarrow  \begin{cases}
                0, &p>1;\\
-3\lambda_{122}\chi_{12}(g_{(1\,l)})\chi_1(g_{(q+1\,l)})a_{(1\,q)}
a_{1}g_{(q+1\,l)}g_{(1\,l)}g_{(2\,l)} \\
\hspace*{6.6cm} \ot y_{(q+1\,l)}y_{(3\,l)}^2,
&p=1.
\end{cases}
\end{align*}
Hence $X_qBX_p=\xi BX_qX_p$. In particular, they do not contribute to $\sigma_l$ and the claim follows.

We deal with cases (L1), (L2), (L3) using the identities developed in \S \ref{sec:technical}.
We need to take into account
Equation
\eqref{eq:conm-lift-1}.

\begin{claim}\label{claim:2}
Case (L1) contributes to  $\sigma_l$ with \eqref{eqn:ul-2}.
\end{claim}
We have to analyze terms
$BX_pX_q$, $X_pBX_q$, $X_pX_qB$, $p<q$. Now, if $p>1$, as
$[y_{(1\,l)},y_{(p+1\,l)}]_c=[y_{(1\,l)},y_{(q+1\,l)}]_c=0$ it follows
\begin{multline*}
B^{p,q}:=BX_pX_q + X_pBX_q +
X_pX_qB\\=(1+\xi+\xi^2)\chi_{(1\,q)}(g_{p+1\,l})a_{(1\,p)}a_{(1\,q)}
g_{(p+1\,l)}g_{(q+1\,l)}g_{(1\,l)}\\ \ot y_{(p+1\,l)}y_{(q+1\,l)}y_{(1\,l)}=0.
\end{multline*}
If $p=1$, then still $[y_{(1\,l)},y_{(q+1\,l)}]_c=0$ as $q+1\geq 3$ and using Lemma
\ref{lem:tech1} to compute $[y_{(1\,l)},y_{(2\,l)}]_c$ we get:
\begin{align*}
B^{1,q}&:= BX_pX_q + X_pBX_q +
X_pX_qB \\ &\rightsquigarrow
-3\lambda_{122}\chi_{12}(g_{(1\,l)})\chi_{(1\,q)}(g_{(1\,l)})
\chi_{(1\,q)}(g_{(2\,l)})
a_{1}a_{(1\,q)}g_{(2\,l)}g_{(q+1\,l)}g_{(1\,l)}\\ &\hspace*{9cm}\ot
y_{(3\,l)}^2y_{(q+1\,l)}.
\end{align*}
Hence, as $\lambda_{122}[y_{(3\,l)},y_{(q+1\,l)}]_c=0$, $q>2$:
\begin{align*}
 B^{1,q}\rightsquigarrow
\begin{cases}
 -3\xi\lambda_{122}\lambda_{112}(1-g_{112})g_{(1\,l)}g_{(2\,l)}g_{(3\,l)}\ot
y_{(3\,l)}^3, & q=2;\\
3(1-\xi)\lambda_{122}\lambda_{112}\chi_{12}(g_{(3\,l)})
a_{(3\,q)}g_{(1\,l)}g_{(2\,l)}g_{(q+1\,l)}\\
\hspace*{7cm} \ot
y_{(q+1\,l)}y_{(3\,l)}^2, & q\geq 3.
\end{cases}
\end{align*}
Notice that in this way $B^{1,2}$ will contribute to $\sigma_l$, as by the induction process we have $y_{(3\,l)}^3=\mu_{(3\,l)}$, that is
we get a term \eqref{eqn:ul-2}.

\begin{claim}\label{claim:3}
Case (L2) does not contribute to  $\sigma_l$.
\end{claim}
We have to deal with terms
$X_pX_q^2$, $X_qX_pX_q$, $X_pX_qX_p$
and $X_p^2X_q$. According to \eqref{eq:conm-lift-1}, we have to distinguish cases
\begin{enumerate}
\item[(L2i)]$p+1<q< l$,
\item[(L2ii)] $p+1=q<l-1$,
\item[(L2iii)]  $p+1=q=l-1$.
\end{enumerate}
In case (L2i) we have
\begin{align*}
X_pX_q^2&=\chi_{(1\,q)}(g_{(q+1\,l)}g_{(p+1\,l)}^2)a_{(1\,p)}a_{(1\,q)}^2
g_{(q+1\,l)}^2g_{(p+1\,l)}\ot y_{(p+1\,l)}y_{(q+1\,l)}^2\\
&=\chi_{(1\,q)}(g_{(q+1\,l)})\chi_{(1\,l)}(g_{(p+1\,l)})^2a_{(1\,p)}a_{(1\,q)}^2
 g_{(q+1\,l)}^2g_{(p+1\,l)}\ot y_{(q+1\,l)}^2y_{(p+1\,l)}
\end{align*}
which does not contribute to $\uvi_l$. For (L2ii):
\begin{align*}
X_p&X_q^2\rightsquigarrow   \chi_{(1\,q)}(g_{(q+1\,l)}g_{(p+1\,l)})\chi_{(1\,l)}(g_{(p+1\,l)})a_{(1\,p)}a_{(1\,q)}^2
g_{(q+1\,l)}^2g_{(p+1\,l)}\\
&\hspace*{6cm}\ot y_{(q+1\,l)}[y_{(p+1\,l)},y_{(q+1\,l)}]_c\\
&\quad +\chi_{(1\,q)}(g_{(q+1\,l)}g_{(p+1\,l)}^2)a_{(1\,p)}a_{(1\,q)}^2
g_{(q+1\,l)}^2g_{(p+1\,l)}\\
&\hspace*{6cm}\ot [y_{(p+1\,l)},y_{(q+1\,l)}]_cy_{(q+1\,l)}\\
&=-3\xi^2\lambda_{p+1p+2p+2}\chi_{(p+2\,l)}(g_{(p+3\,l)})\chi_{(1\,p+1)}(g_{(p+2\,l)}g_{(p+1\,l)})\chi_{(1\,l)}(g_{(p+1\,l)})
\\
&\hspace*{4cm}a_{(1\,p)}a_{(1\,p+1)}^2g_{(p+2\,l)}^2g_{(p+1\,l)}\ot y_{(p+2\,l)}y_{(p+3\,l)}^2\\
&\quad -3\xi^2\lambda_{p+1p+2p+2}\chi_{(p+2\,l)}(g_{(p+3\,l)})\chi_{(1\,p+1)}(g_{(p+2\,l)}g_{(p+1\,l)}^2)\\
&\hspace*{4cm}a_{(1\,p)}a_{(1\,p+1)}^2
g_{(p+2\,l)}^2g_{(p+1\,l)}\ot y_{(p+3\,l)}^2y_{(p+2\,l)}
\end{align*}
and we see that this does not contribute to $\uvi_l$, using Lemma \ref{lem:l112} to deduce $\lambda_{p+1p+2p+2}[y_{(p+2\,l)},y_{(p+3\,l)}]_c=0$.
The same holds for (L2iii), as in this case:
\begin{align*}
X_pX_q^2\rightsquigarrow  \lambda_{l-1ll}(1+\xi^2)\chi_{(1\,l-1)}(g_{(l-1\,l)}^2g_{l})a_{(1\,l-2)}a_{(1\,l-1)}^2
g_{l}^2g_{(l-1\,l)}\ot y_{l}.
\end{align*}

The same holds for the combinations $X_qX_pX_q$, $X_pX_qX_p$ and $X_p^2X_q$.

\begin{claim}\label{claim:4}
If $p<3$, then case (L3) does not contribute to  $\sigma_l$. 

If $p\geq 3$, then case (L3) contributes to  $\sigma_l$ with \eqref{eqn:monomials-1}.
\end{claim}
Here we deal with terms $X_pX_qX_r$, $X_pX_rX_q$, $X_qX_rX_p$, $X_rX_pX_q$, $X_qX_pX_r$, $p<q<r\leq l$, which we denote
by $C^{x,y,z}$, $x,y,z\in\{p,q,r\}$. By the computations above and the commutation rule \eqref{eq:conm-lift-1} we see that we will get a factor
contributing to $\sigma_l$  if and only if
$$
q=p+1,\qquad  r=p+2,
$$
and $p$ is on the left of $q$.

Thus we are left with cases $C^{p,q,r}$, $C^{p,r,q}$, $C^{r,p,q}$, $q=p+1$, $r=p+2$. 
Set:
\begin{align}\label{eqn:alpha_p}
  c_p=\begin{cases}
 -3\lambda_{p+1p+2p+2}\chi_{p+2}(g_{(p+3\,l)})\mu_{(p+3\,l)}, & p\leq l-3\\
\lambda_{p+1p+2p+2}, & p=l-2.
          \end{cases}
\end{align}
Then, for each of these terms, the corresponding factor in $\k$ that arises in the second tensorand ($\chi_{(l+1\,l)}:=\eps$) is:
\begin{align*}
&C^{p,q,r}\dashrightarrow   c_p,
&&C^{p,r,q}\dashrightarrow   c_p\chi_{(r+1\,l)}(g_{(p+1\,l)}),
&&C^{r,p,q}\dashrightarrow   c_p.
\end{align*}
Set, with the convention, for the case $r=l$, $g_{(l+1\,l)}=1$, $\chi_{(l+1\,l)}=\eps$,
\begin{align}\label{eqn:omegas}
\omega_{x,y,z}&=\chi_{(1\,z)}(g_{(y+1\,l)}g_{(x+1\,l)})\chi_{(1\,y)}(g_{(x+1\,l)}), \qquad x,y,z\in\{p,q,r\}.
\end{align}
Set also, $g_{p,q,r}:=g_{(p+1\,l)}g_{(p+2\,l)}g_{(p+3\,l)}$ and let us set:
\begin{align}\label{eqn:apqr}
\begin{split}
a_{p,q,r}&:=a_{(1\,p)}a_{(1\,q)}a_{(1\,r)}, \qquad 
a_{p,r,q}:=a_{(1\,p)}a_{(1\,r)}a_{(1\,q)}, \\
a_{r,p,q}&:=a_{(1\,r)}a_{(1\,p)}a_{(1\,q)}.
\end{split}
\end{align}
Set, cf. Remark \ref{rem:p+3}:
\begin{align}\label{eqn:xip}
\Xi_p&=\begin{cases}
             (1-\xi^2)^3=3(\xi-\xi^2), & p<l-2;\\
(1-\xi^2)^2=-3\xi^2, & p=l-2,
            \end{cases}
\end{align}
Hence the contribution of these terms to $\sigma_l$  is
\begin{align}\label{eqn:monomials}
 c_p\Xi_p\Big(\omega_{p,q,r}a_{p,q,r}
+\omega_{p,r,q}\chi_{(p+3\,l)}(g_{(q\,l)})
a_{p,r,q}
+\omega_{r,p,q}a_{r,p,q}\Big)
g_{p,q,r}.
\end{align}
Notice that
\begin{align*}
 c_p\omega_{p,q,r}&= c_p\xi^2\chi_{(1\,q)}(g_{(p+1\,l)})
\chi_{(1\,q)}(g_{(r+1\,l)})\chi_{(1\,p)}(g_{(r+1\,l)})\\
 c_p\omega_{p,r,q}\chi_{(p+3\,l)}(g_{(q\,l)})
&= c_p\chi_{(1\,q)}(g_{ (p+1\,l)})
\chi_{(1\,q)}(g_{(r+1\,l)})\chi_{(1\,p)}(g_{(r+1\,l)})
\chi_r(g_{(1\,p)})\\
 c_p\omega_{r,p,q}&=
 c_p\chi_{(1\,q)}(g_{ (p+1\,l)})
\chi_{(1\,q)}(g_{(r+1\,l)})\chi_{(1\,p)}(g_{ (r+1\,l)}).
\end{align*}
Set 
\begin{align}\label{eqn:widehatalpha_p}
 d'_p
&= \begin{cases}-3\mu_{(p+3\,l)}
\chi_{(1\,q)}(g_{qr})
\chi_{qqr}(g_{(r+1\,l)}),&p<l-2;\\
\chi_{(1\,q)}(g_{ (q\,l)})
,&p=l-2,
\end{cases}
\end{align}
and $d_p=\lambda_{qrr}d'_p$. Observe that
\begin{align*}
 c_p\chi_{(1\,q)}(g_{ (p+1\,l)})
\chi_{(1\,q)}(g_{(r+1\,l)})\chi_{(1\,p)}(g_{ (r+1\,l)})= d_p.
\end{align*} 
To see this, we use the identity
\begin{multline*}
\mu_{(p+3\,l)}\lambda_{qrr}\chi_{(1\,q)}(g_{ (q\,l)})\chi_{(1\,r)}(g_{(r+1\,l)})\chi_{(1\,p)}(g_{ (r+1\,l)})\\
=\mu_{(p+3\,l)}\lambda_{qrr}\chi_{(1\,p)}(g_{ (r+1\,l)})^3
\chi_{(1\,q)}(g_{qr})
\chi_{qqr}(g_{(r+1\,l)})
\end{multline*}
and $\mu_{(p+3\,l)}\lambda_{qrr}\chi_{(1\,p)}(g_{ (r+1\,l)})^3=
\mu_{(p+3\,l)}\lambda_{qrr}\chi_{ (r+1\,l)}(g_{(1\,p)})^{-3}=\mu_{(p+3\,l)}$.

Hence \eqref{eqn:monomials} becomes
\begin{align*}
 d_p\Xi_p\Big(\xi^2a_{p,q,r}
+\chi_{r}(g_{(1\,p)})
a_{p,r,q}
+a_{r,p,q}\Big)
g_{p,q,r}
= d'_p\Xi_p\varsigma^p
g_{p,q,r}.
\end{align*}
Adding all of these terms and reordering the scalars, 
we get \eqref{eqn:monomials-intro-1}.

\smallbreak

Finally, adding all of these contributions, we obtain $\sigma_{(i\,l)}(\boldsymbol\lambda,\boldsymbol\mu)$ as in \eqref{eqn:uil} and the proposition follows.
\epf

\begin{exa}\label{ej:not-in-H-N3}
Set $\theta=5$, so $\I=\I_5$ and consider the braiding matrix $$\qb=\left(\begin{smallmatrix}
\xi&\xi&1&1&1\\
\xi&\xi&\xi^2&1&1\\
1&1&\xi&1&1\\
1&1&\xi^2&\xi&\xi\\
1&1&1&\xi&\xi
\end{smallmatrix}\right).$$ 
Set $G=(\Z/3n\Z)^5$, $n\geq 2$, so $G$ is an abelian group such that $V\in\ydh$, $H=\k G$. 
Indeed, let $g_i$, $i\in\I$, the generators of each cyclic factor. Let $q\in\G'_{3n}$ with $q^n=\xi$. Observe that $\widehat{G}$ is generated by $\varphi_i$, $i\in\I$, with $\varphi_i(g_i)=q$ and $\varphi_i(g_j)=1$, $i\neq j\in\I$.
A principal realization is given by $((g_i,\chi_i)\big)_{i\in\I}$, for $\chi_1=\chi_2=\varphi_1^n\varphi_2^n$, 
$\chi_4=\chi_5=\varphi_4^n\varphi_5^n$, and $\chi_3=\varphi_2^{2n}\varphi_3^n\varphi_4^{2n}$.   
In particular $\chi_{112}=\chi_{455}=\eps$. 

Let us choose $\boldsymbol\lambda$ such that
all $\lambda_{iij}=0$ except $\lambda_{112}$, $\lambda_{455}$. Choose 
$\boldsymbol\mu$ with $\mu_{(k\,l)}=0$ for every $1\leq k\leq l\leq 5$.
Then $\mL(\boldsymbol\lambda,\boldsymbol\mu)$ is the algebra generated by $\Gamma$ and $a_1,\dots,a_5$ satisfying:
\begin{align*}
a_{ik}&=0, \quad  \mid i-k\mid>1, \qquad  a_{(k\,l)}^3=0, \quad  \mid 
k-l\mid<4,\\
a_{iiij}&=\begin{cases}
          \lambda_{112}(1-g_1^2g_2), & i=1, j=2,\\
\lambda_{455}(1-g_4g_5^2), & i=4, j=5,\\
0, & \text{else},
         \end{cases} \qquad \mid i-j\mid=1,\\
a_{(1\,5)}^3&=9\lambda_{112}\lambda_{455}\sum_{\sigma\in\Sym_3}(-1)^{|\sigma|}h_{\sigma,1}a_{(3\,\sigma(5))}a_{(2\,\sigma(4))}a_{(1\,\sigma(3))}g_4^2g_{5}.
\end{align*}
The scalars $h_{\sigma,1}\in\k^\times$, $\sigma\in\Sym_3$, are as in Corollary \ref{cor:suma-a} and can be explicitly computed from the matrix $\qb$. 
\end{exa}

\begin{rem}\label{rem:not-in-H-N3}
If $\theta\geq 4$, then the relations of the Nichols algebra $\B(V)$ become deformed in the lifting $\mL$ by elements in the group algebra
$\k\Gamma\leq H$, as in the case $\ord(\xi)>3$ of \cite{AS2}, as Lemma \ref{lem:l112}
only allows a single pair $(\lambda_{kkk+1},\lambda_{kk+1k+1})$ to
have a nonzero entry. If $\theta\geq 5$ the relations may be deformed in higher strata of the coradical filtration, as in Example \ref{ej:not-in-H-N3}.
\end{rem}
Now, we  give a full description of $\varsigma^{p}_{i}(\bs\lambda,\bs\mu)$, cf. \eqref{eqn:monomials-1}, as a linear combination in the PBW basis.
We set $q=p+1$, $r=p+2$ and $j=i+1$, $k=i+2$. We consider the action of $\s_3$ on $\{r,q,p\}$ by  
\begin{align*}
(12)(r)=q, \qquad (23)(q)=p. 
\end{align*}
\begin{cor}\label{cor:suma-a}
If $p=i,j$, then $\varsigma^{p}_i(\bs\lambda,\bs\mu)=0$. When $p>i+2$, 
\begin{align}\label{eqn:ul-3}
\begin{split}
\varsigma^{p}_{i}(\bs\lambda,\bs\mu)=-3\lambda_{qrr}&\lambda_{qqr}\chi_{(i\,p)}(g_{q})a_{(i\,p)}^3g_{qqr}\\
&-3\lambda_{qrr}\lambda_{iij}\sum_{\sigma\in\Sym_3}(-1)^{|\sigma|}h_{\sigma,i}a_{(k\,\sigma(p))}a_{(j\,\sigma(q))}a_{(i\,\sigma(r))}.
\end{split}
\end{align}
for $h_{\sigma,i}\in\k$, $\sigma\in\Sym_3$, given by:
\begin{align*}
h_{\id,i}&=\xi \chi_{qqr}(g_{(i\,p)})\chi_{(i\,r)}(g_{(j\,q)}), &h_{(12),i}&=(\xi^2-1)\chi_{qqr}(g_{(i\,p)})\chi_{i}(g_{(k\,q)}), \\
h_{(23),i}&=\xi \chi_{r}(g_{i})\chi_{i}(g_{(j\,p)}), & h_{(13),i}&=\xi (\xi-2)\chi_{(k\,p)}(g_{ij}), \\
h_{(123),i}&=2\chi_{r}(g_{(i\,p)})\chi_{i}(g_{(k\,p)}), &h_{(132),i}&=\xi^2\chi_{(k\,q)}(g_{(i\,r)})\chi_{(j\,p)}(g_{r}). 
\end{align*}
\end{cor}

\pf
First, we show that $\varsigma^{p}(\bs\lambda,\bs\mu)$ equals
\begin{align}\label{eqn:explicit-1}
\lambda_{qrr}\xi^2a_{(i\,p)}[a_{(i\,q)},a_{(i\,r)}]_c-\lambda_{qrr}\xi^2\chi_{r}(g_{(i\,p)})[a_{(i\,p)},a_{(i\,r)}]_ca_{(i\,q)}.
\end{align}
In particular, by Lemma \ref{lem:l112} and Corollary \ref{cor:p-p+1}, we have
$\lambda_{qrr}\varsigma^{p}(\bs\lambda,\bs\mu)=0$ if $p-i<3$.
Hence $\varsigma_{(i\,i)}=\varsigma_{(i\,i+1)}=0$.
Indeed, it follows that:
\begin{align*}
\varsigma^{p}&(\bs\lambda,\bs\mu)
=
\lambda_{qrr}\alpha_p\chi_{(i\,q)}(g_{(i\,p)})\,a_{(i\,r)}a_{(i\,q)}a_{(i\,p)}\\
&\,
+\lambda_{qrr}\Big(\xi^2\chi_{(i\,r)}(g_{(i\,q)}g_{(i\,p)})
+1+\chi_{r}(g_{(i\,p)})\chi_{(i\,r)}(g_{(i\,p)})\Big)
a_{(i\,r)}[a_{(i\,p)},a_{(i\,q)}]_c\\
&\,
+\xi^2\lambda_{qrr}a_{(i\,p)}[a_{(i\,q)},a_{(i\,r)}]_c
\\
&\,
+\lambda_{qrr}\Big(\chi_{p+2}(g_{(i\,p)})+\xi^2\chi_{(i\,r)}(g_{(i\,q)})\Big)[a_{(i\,p)},a_{(i\,r)}]_ca_{(i\,q)},
\end{align*}
for $\alpha_p=1+\xi^2\chi_{(i\,r)}(g_{(i\,q)}g_{(i\,p)})
+\chi_{p+2}(g_{(i\,p)})\chi_{(i\,r)}(g_{(i\,p)})$. Notice that $\alpha_p=1+\xi\chi_{r}(g_{q})\chi_{qrr}(g_{(i\,p)})
+\xi\chi_{qrr}(g_{(i\,p)})$ and thus it follows that $\lambda_{qrr}\alpha_p=0$ as $\lambda_{qrr}\chi_{qrr}=\lambda_{qrr}\eps$,  $\lambda_{qrr}\chi_{r}(g_{q})=\lambda_{qrr}\xi$ and $1+\xi+\xi^2=0$.
On the other hand, we use $\lambda_{qrr}\chi_{(i\,r)}(g_{(i\,q)})=\lambda_{qrr}\xi^2\chi_{r}(g_{(i\,p)})$ to simplify the coefficients of third and fourth summands. As for the second,
we have 
\begin{align*}
\lambda_{qrr}\chi_{(i\,r)}(g_{(i\,q)}g_{(i\,p)})&=\lambda_{qrr}\xi^2\chi_{r}(g_{(i\,p)})\chi_{(i\,r)}(g_{(i\,p)})\\
&=\lambda_{qrr}\chi_{r}(g_{(i\,p)})\chi_{qr}(g_{(i\,p)})=\lambda_{qrr}.
\end{align*}
Also, $\lambda_{qrr}\chi_{r}(g_{(i\,p)})\chi_{(i\,r)}(g_{(i\,p)})=\lambda_{qrr}\xi\chi_{qrr}(g_{(i\,p)})=\lambda_{qrr}$ and thus the coefficient is $\lambda_{qrr}(1+\xi+\xi^2)=0$. Hence, we have \eqref{eqn:explicit-1}.

Next we show \eqref{eqn:ul-3}, using Proposition \ref{pro:zpqr}.  We have:
\begin{align*}
\varsigma^{p}&(\bs\lambda,\bs\mu)=
\lambda_{qrr}\chi_{q}(g_{(i\,p)})(1+\xi^2\chi_{qrr}(g_{(i\,p)})+\xi) a_{(i\,r)}a_{(i\,q)}a_{(i\,p)}\\
&-3\xi \lambda_{qrr}\lambda_{qqr}\chi_{(i\,r)}(g_{q})a_{(i\,p)}^3g_{qqr}\\
&-3\xi \lambda_{qrr}\lambda_{iij}\chi_{qqr}(g_{(i\,p)})\chi_{(i\,r)}(g_{(j\,q)})a_{(k\,r)}a_{(j\,q)} a_{(i\,p)}\\
&-3\xi \lambda_{qrr}\lambda_{iij}\Big(\chi_{qqr}(g_{(i\,p)})\chi_{(k\,r)}(g_{(j\,q)})\chi_{(j\,p)} (g_{(i\,q)})\\
&\qquad +\xi\chi_{(i\,r)}(g_{(j\,p)})\chi_{r}(g_{(i\,p)})\Big)a_{(k\,r)}a_{(j\,p)}a_{(i\,q)}\\
&+3\xi^2 \lambda_{qrr}\lambda_{iij}\chi_{qqr}(g_{(i\,p)})\chi_{i}(g_{(k\,q)})a_{(k\,q)}a_{(j\,r)}a_{(i\,p)}\\
&-3\xi^2\lambda_{iij} \lambda_{qrr}\Big(\chi_{(k\,q)}(g_{(i\,r)})\chi_{(j\,p)}(g_{r}) \\
&\qquad -2\chi_{(i\,q)}(g_{(j\,p)})\\
&\qquad \chi_{(i\,r)}(g_{(i\,p)})\chi_{(k\,q)}(g_{(i\,r)})\chi_{r}(g_{i})\Big)a_{(k\,q)}a_{(j\,p)}a_{(i\,r)} \\
&+3\xi^2 \lambda_{qrr}\lambda_{iij}\Big(\chi_{(k\,p)}(g_{ij})\chi_{r}(g_{(i\,q)})\\
&\qquad +\xi\chi_{i}(g_{(k\,p)})\chi_{r}(g_{(i\,p)})\Big)a_{(k\,p)}a_{(j\,r)}a_{(i\,q)}\\
&+3\xi^2 \lambda_{qrr}\lambda_{iij}\Big((1-\xi^2)\chi_{(k\,p)}(g_{ij})\\
&\qquad +\xi\chi_{r}(g_{(i\,p)})\chi_{i}(g_{(k\,p)})\chi_{(i\,q)}(g_{r})\\
 &\qquad +\xi\chi_{(k\,p)}(g_{(j\,r)})\chi_{(j\,q)}(g_{(i\,r)})\Big)a_{(k\,p)}a_{(j\,q)}a_{(i\,r)}.
\end{align*}
First, $\lambda_{qrr}(1+\xi^2\chi_{qrr}(g_{(i\,p)})+\xi)=\lambda_{qrr}(1+\xi^2+\xi)=0$. Next, observe that
\begin{align*}
\lambda_{qrr}&\lambda_{iij}\Big(\chi_{qqr}(g_{(i\,p)})\chi_{(k\,r)}(g_{(j\,q)})\chi_{(j\,p)} (g_{(i\,q)})+\xi\chi_{(i\,r)}(g_{(j\,p)})\chi_{r}(g_{(i\,p)})\Big)\\
&=\lambda_{qrr}\lambda_{iij}\xi^2\chi_{r}(g_{i})\chi_{i}(g_{(j\,p)})\Big(\chi_q(g_{iij})\chi_{(k\,p)} (g_{iij})+1\Big)\\
&=\lambda_{qrr}\lambda_{iij}\xi^2\chi_{r}(g_{i})\chi_{i}(g_{(j\,p)})(\xi^2+1)=-\lambda_{qrr}\lambda_{iij}\chi_{r}(g_{i})\chi_{i}(g_{(j\,p)}).
\end{align*}
Similarly,
\begin{align*}
\lambda_{iij} &\lambda_{qrr}\Big(\chi_{(k\,q)}(g_{(i\,r)})\chi_{(j\,p)}(g_{r}) + \chi_{(i\,r)}(g_{(i\,p)})\chi_{(k\,q)}(g_{(i\,r)})\chi_{r}(g_{i})\\
 &-2\chi_{(i\,q)}(g_{(j\,p)})\Big)=
\lambda_{iij} \lambda_{qrr}\chi_{(k\,q)}(g_{(i\,r)})\chi_{(j\,p)}(g_{r})\Big(1 + \xi-2\xi\Big)\\
&=
\lambda_{iij} \lambda_{qrr}\chi_{(k\,q)}(g_{(i\,r)})\chi_{(j\,p)}(g_{r})(1 - \xi).
\end{align*}
Also, we have
\begin{align*}
\lambda_{qrr}&\lambda_{iij}\Big(\chi_{(k\,p)}(g_{ij})\chi_{r}(g_{(i\,q)})
+\xi\chi_{i}(g_{(k\,p)})\chi_{r}(g_{(i\,p)})\Big)\\
&=2\lambda_{qrr}\lambda_{iij}\xi\chi_{r}(g_{(i\,p)})\chi_{i}(g_{(k\,p)}).
\end{align*}
Finally, 
\begin{align*}
\lambda_{qrr}&\lambda_{iij}\Big((1-\xi^2)\chi_{(k\,p)}(g_{ij})
 +\xi\chi_{r}(g_{(i\,p)})\chi_{i}(g_{(k\,p)})\chi_{(i\,q)}(g_{r})\\
 &+\xi\chi_{(k\,p)}(g_{(j\,r)})\chi_{(j\,q)}(g_{(i\,r)})\Big)=\lambda_{qrr}\lambda_{iij}\chi_{(k\,p)}(g_{ij})(1-2\xi^2)
\end{align*}
Thus, we have
\begin{align*}
\varsigma^{p}&(\bs\lambda,\bs\mu)=
-3\xi \lambda_{qrr}\lambda_{qqr}\chi_{(i\,r)}(g_{q})a_{(i\,p)}^3g_{qqr}\\
&-3\xi \lambda_{qrr}\lambda_{iij}\chi_{qqr}(g_{(i\,p)})\chi_{(i\,r)}(g_{(j\,q)})a_{(k\,r)}a_{(j\,q)} a_{(i\,p)}\\
&+3\xi \lambda_{qrr}\lambda_{iij}\chi_{r}(g_{i})\chi_{i}(g_{(j\,p)})a_{(k\,r)}a_{(j\,p)}a_{(i\,q)}\\
&-3(1 - \xi^2)\lambda_{qrr}\lambda_{iij}\chi_{qqr}(g_{(i\,p)})\chi_{i}(g_{(k\,q)})a_{(k\,q)}a_{(j\,r)}a_{(i\,p)}\\
&-3\xi^2\lambda_{iij} \lambda_{qrr}\chi_{(k\,q)}(g_{(i\,r)})\chi_{(j\,p)}(g_{r})a_{(k\,q)}a_{(j\,p)}a_{(i\,r)} \\
&+6 \lambda_{qrr}\lambda_{iij}\chi_{r}(g_{(i\,p)})\chi_{i}(g_{(k\,p)})a_{(k\,p)}a_{(j\,r)}a_{(i\,q)}\\
&-3\xi (2-\xi)\lambda_{qrr}\lambda_{iij}\chi_{(k\,p)}(g_{ij})a_{(k\,p)}a_{(j\,q)}a_{(i\,r)}.
\end{align*}
Hence the lemma follows by defining the scalars $h_{\sigma,i}$ appropriately.
\epf
\subsection{Technical identities}\label{sec:technical}

To compute the elements
$\varsigma_{(i\,l)}$ in \eqref{eqn:monomials-intro-1} in the PBW basis, we need a large series of technical identities involving commutators. 
This is the 
content of this section.
\begin{lem}\label{lem:tech1-lift}
The following identities hold in $\mtL$.
 \begin{enumerate}
\item $[a_{(1\,l)},a_{2}]_c=\begin{cases}
                   \lambda_{122}(1-\xi^2)\chi_2(g_{3})a_{3}-\lambda_{223}
(1-\xi^2)a_1g_{223}, & l=3;\\
          \lambda_{122}(1-\xi^2)\chi_2(g_{(3\,l)})a_{(3\,l)}, &l\geq 4.
                            \end{cases}$
\item $
[a_{(1\,l)},a_p]_c=0
$, $3\leq p<l-1$.
\item $[a_{(1\,l)},a_{(p\,k)}]_c=0
$, $3\leq p\leq k<l-1$.
\item  $
[a_{(1\,l)},a_{l-1}]_c=-\lambda_{l-1l-1l}(1-\xi^2)a_{(1\,l-2)}g_{l-1l-1l}.
$
\item $
[a_{(1\,l)},a_l]_c=-\lambda_{l-1ll}(1-\xi^2)a_{(1\,l-2)}g_{l-1ll}.
$
 \end{enumerate}
\end{lem}
\pf
(1) Case $l=3$ follows once again mimicking
\cite[Lemma 1.11]{AS2} as in Lemma \ref{lem:tech1}. The general case $l\geq 4$ follows as
in Lemma
\ref{lem:tech1}: in this situation, if $\lambda_{223}\neq0$, then $\lambda_{122}=0$ by
Lemma \ref{lem:l112} and
\begin{align*}
 [a_{(1\,l)},a_{2}]_c=-\lambda_{223}
(1-\xi^2)\chi_{(4\,l)}(g_{23})[a_1,a_{(4\,l)}]_cg_{223}=0,
\end{align*}
using q-Jacobi \eqref{eq:qjacob}.
For (2), first we have that
\begin{multline*}
[a_{(1\,l)},a_p]_c=[[a_{(1\,p+1)},a_{(p+2\,l)}]_c,a_p]_c=\chi_p(g_{(p+2\,l)})
[a_{(1\,p+1)
},a_p]_ca_{(p+2\,l)}\\-\chi_{(p+2\,l)}(g_{(1\,p+1)})a_{(p+2\,l)}[a_{(1\,p+1)
},a_p]_c.
\end{multline*}
Now, by (1),
$\lambda_{ppp+1}\chi_{ppp+1}=\lambda_{ppp+1}\eps$ and $[a_{(1\,p-2)},a_p]_c=0$,
\begin{align*}
 [a_{(1\,p+1)},a_p]_c&=[a_{(1\,p-2)},[a_{(p-1\,p+1)},a_p]_c]_c
\\
 &
=-\lambda_{ppp+1}
(1-\xi^2)[a_{(1\,p-2)},a_{p-1}g_{ppp+1}]_c
\\
&
=-\lambda_{ppp+1}(1-\xi^2)\Big(a_{(1\,p-2)}a_{p-1}g_{ppp+1}\\
&\qquad -\chi_{p-1}(g_{(1\,p-2)})a_{
p-1}g_{ppp+1}a_{(1\,p-2)}\Big)\\
&=-\lambda_{ppp+1}(1-\xi^2)a_{(1\,p-1)}g_{ppp+1},
\end{align*}
as $\lambda_{ppp+1}\chi_{(1\,p-2)}(g_{ppp+1})=\lambda_{ppp+1}$. In particular, this shows (4) for $p=l-1$. Now, if $p<l-1$ we get
\begin{align*}
[a_{(1\,l)},a_p]_c&=
 -\lambda_{ppp+1}(1-\xi^2)\chi_p(g_{(p+2\,l)})\big(a_{(1\,p-1)}g_{ppp+1}a_{(p+2\,l)}\\
&\qquad
-\chi_{(p+2\,l)}(g_p)\chi_{(p+2\,l)}(g_{(1\,p+1)})a_{(p+2\,l)}a_{(1\,p-1)}g_{ppp+1}\big)\\
&=-\lambda_{ppp+1}(1-\xi^2)\chi_p(g_{(p+2\,l)})\chi_{(p+2\,l)}(g_{ppp+1})\\
&\qquad\qquad \qquad \qquad \qquad\qquad \qquad \qquad  [a_{(1\,p-1)},
a_{(p+2\,l)}]_cg_{ppp+1}=0.
\end{align*}
(3) follows from (2) by induction. For (5), we get, as
$\lambda_{l-1ll}\chi_{l-1ll}=\lambda_{l-1ll}\eps$,
\begin{align*}
[a_{(1\,l)},a_l]_c&=[[a_{(1\,l-2)},a_{l-1l}]_c,a_l]_c
=\lambda_{l-1ll}(\chi_{(1\,l-2)}(g_{
l-1ll})-1)a_{1\,l-2}g_{l-1ll}
\end{align*}
and since
$$\lambda_{l-1ll}\chi_{1\,l-2}(g_{l-1ll})=\lambda_{l-1ll}\chi_{(1\,l-2)}(g_{l-1ll})
\chi_{l-1ll}(g_{(1\,l-2)})=\lambda_{l-1ll}\xi^2$$
the lemma follows.
\epf

\begin{rem}
As $
[a_{2},a_{(1\,3)}]_c=-\chi_{(1\,3)}(g_2)[a_{(1\,3)},a_{2}]_c$, we get
$$
[a_{2},a_{(1\,3)}]_c=\lambda_{122}
(1-\xi)a_{3}-\lambda_{223}\chi_1(g_2)(\xi^2-\xi)a_1g_{223}.
$$
\end{rem}

\begin{lem}\label{lem:corchete-1p,3l}
The following identities hold in $\mtL$.
\begin{enumerate}
\item $[a_{(1\,l)},a_{(3\,p)}]_c=[a_{(3\,p)},a_{(1\,l)}]_c=0$, $3\leq p<l-1$.
\item
$[a_{(1\,l)},a_{(3\,l)}]_c=-\lambda_{l-1ll}(1-\xi)a_{(1\,l-1)}g_{l-1ll}$\newline
\hspace*{4cm} $+3\lambda_{l-1ll}\chi_{l-1}(g_{(1\,l-2)})a_{l-1}a_{(1\,l-2)}g_{l-1ll}$.
\item $[a_{(1\,l)},a_{(3\,l-1)}]_c=\begin{cases}
                                    -\lambda_{334}(1-\xi^2)a_{12}g_{334}, & 
l=4,\\
 3\xi^2\lambda_{l-1l-1l}\chi_{(3\,l-2)}(g_{(1\,l)}) & \\
\qquad \qquad   a_{(3\,l-2)}a_{(1\,l-2)}g_{l-1l-1l}, & l\geq 5.
                                   \end{cases}
$
\end{enumerate}
\end{lem}
\pf
(1)  follows by induction on $p$ and using q-Jacobi \eqref{eq:qjacob}, case $p=3$
being Lemma \ref{lem:tech1-lift} (2).

(2) Using q-Jacobi \eqref{eq:qjacob} and Lemma \ref{lem:tech1-lift} (4-5), we have
\begin{align*}
 [a_{(1\,l)},&a_{(3\,l)}]_c=[[a_{(1\,l)},a_{l-1}]_c,a_l]_c+\chi_{l-1}(g_{(1\,l)})a_{l-1}[a_{(1\,l)},a_l]_c\\
& -\chi_l(g_{l-1})[a_{(1\,l)},a_l]_ca_{l-1}=-\lambda_{l-1l-1l}(1-\xi^2)[a_{(1\,l-2)},a_l]_cg_{l-1l-1l}\\
& -\lambda_{l-1ll}(1-\xi)a_{(1\,l-1)}g_{l-1ll}+3\lambda_{l-1ll}\chi_{l-1}(g_{(1\,l-2)})a_{l-1}a_{(1\,l-2)}g_{l-1ll}.
\end{align*}

(3) Case $l=4$ is Lemma \ref{lem:tech1-lift} (1), using q-Jacobi \eqref{eq:qjacob}.

Now, if $l\geq 5$, using q-Jacobi \eqref{eq:qjacob}, item (1) and Lemma \ref{lem:tech1-lift} (4) we get
\begin{align*}
[a_{(1\,l)},& a_{(3\,l-1)}]_c=-\lambda_{l-1l-1l}(1-\xi)\chi_{l-1}(g_{(3\,l-2)})[a_{(1\,l-2)},a_{(3\,l-2)}]_cg_{l-1l-1l}\\
&- \lambda_{l-1l-1l}(1-\xi^2)^2\chi_{(3\,l-2)}(g_{(1\,l)})a_{(3\,l-2)}a_{(1\,l-2)}g_{l-1l-1l}.
\end{align*}
Hence (3) follows using (2) and $\lambda_{l-1l-1l}\lambda_{l-3l-2l-2}=0$. Also, we use the fact that
$\lambda_{l-1l-1l}\chi_{l-1}(g_{(3\,l-2)})=\lambda_{l-1l-1l}\xi\chi_{(3\,l-2)}(g_{l-1l})$.
\epf
\begin{lem}\label{lem:corchete-1p,1l-parte1}
The following identities hold in $\mtL$.
\begin{enumerate}
\item $[a_{1},a_{12}]_c= \lambda_{112}(1-g_{112})$.
\item\label{it5-lift} $[a_{1},a_{(1\,l)}]_c=
\lambda_{112}(1-\xi^2)a_{(3\,l)}$, $l\geq 3$.
\item\label{it6a-lift} $[a_{12},a_{(1\,3)}]_c=
-3\xi^2\lambda_{112}
\chi_{(1\,3)}(g_2)a_3a_2+\lambda_{112}(1-\xi)a_{23}\newline \hspace*{7cm} -3\lambda_{223}
\xi\chi_1(g_2)a_1^2g_{223}$.
\item\label{it6b-lift} $[a_{(12)},a_{(1\,l)}]_c=-3\xi^2
\lambda_{112}\chi_{(1\,l)}(g_2)a_{(3\,l)}a_2+\lambda_{112}(1-\xi)a_{(2\,l)}$, $l\geq 4$.
\end{enumerate}
\end{lem}
\pf
(1) is by definition. For \eqref{it5-lift}, we have
\begin{multline*}
 [a_{1},a_{(1\,l)}]_c=
[a_{1},[a_{12},a_{(3\,l)}]_c]_c=\lambda_{112}(1-\xi^2)a_{(3\,l)}\\-\lambda_{112}
\left(g_{112}a_{(3\,l)}-\chi_{(3\,l)}(g_{112})a_{(3\,l)}g_{112}\right)
=\lambda_{112}(1-\xi^2)a_{(3\,l)}.
\end{multline*}
\eqref{it6a-lift} and \eqref{it6b-lift} follow as in Lemma \ref{lem:tech3}. In this case,
when $l=3$ and extra term
involving $a_1^2g_{223}$ arises, which gets killed for bigger $l$.
\epf

\begin{pro}\label{pro:corchete-1p,1l-parte2}
The following identities hold in $\mtL$.
\begin{enumerate}
\item\label{item:corchete-1p,1l-1}  For $3\leq p <l-1$:
\begin{align*}
 [a_{(1\,p)},a_{(1\,l)}]_c=-3\xi^2&\lambda_{112}\chi_{(1\,l)}(g_{(2\,p)})
a_{(3\,l)}a_{(2\,p)} +3\lambda_{112}\chi_{1}(g_{(3\,p)})a_{(3\,p)}a_{(2\,l)}.
\end{align*}
\item\label{item:corchete-1p,1l-2} For $l\geq 5$,
\begin{align*}
[a_{(1\,l-1)},&a_{(1\,l)}]_c= -3\xi^2\chi_{(1\,l)}(g_{l-1})\lambda_{l-1l-1l}a_{(1\,l-2)}^2g_{l-1l-1l}\\
&-3\xi^2\lambda_{112}\chi_{(1\,l)}(g_{(2\,l-1)})
a_{(3\,l)}a_{(2\,l-1)}+3\lambda_{112}\chi_{1}(g_{(3\,l-1)})a_{(3\,l-1)}a_{(2\,l)}.
\end{align*}
\end{enumerate}
\end{pro}
\pf
\eqref{item:corchete-1p,1l-1} We use q-Jacobi \eqref{eq:qjacob} and Lemma \ref{lem:corchete-1p,3l} (1) to
get
\begin{align*}
[a_{(1\,p)},a_{(1\,l)}]_c&=\chi_{(1\,l)}(g_{(3\,p)})[a_{(12)},a_{(1\,l)}]_ca_{(3\,p)}-\chi_
{(3\,p)}(g_{12})a_{(3\,p)}[a_{(12)},a_{(1\,l)}]_c.\\
[a_{(12)},a_{(1\,l)}]_c&a_{(3\,p)}=-3\xi^2
\lambda_{112}\chi_{(1\,l)}(g_2)a_{(3\,l)}a_2a_{(3\,p)}+\lambda_{112}(1-\xi)a_{(2\,l)}a_{
(3\,p)}\\
&=-3\xi^2
\lambda_{112}\chi_{(1\,l)}(g_2)\Big(a_{(3\,l)}a_{(2\,p)}+\chi_{(3\,p)}(g_2)a_{(3\,l)}a_{
(3\,p)}a_2\Big)\\
&\qquad +\lambda_{112}(1-\xi)\chi_{(3\,p)}(g_{(2\,l)})a_{(3\,p)}a_{(2\,l)},
 \end{align*}
as $\lambda_{112}\lambda_{223}=0$. We arrive to \eqref{item:corchete-1p,1l-1} using  $\lambda_{112}\lambda_{334}=0$:
\begin{align*}
a_{(3\,p)}[a_{(12)},a_{(1\,l)}]_c&=-3\xi
\lambda_{112}\chi_{(3\,l)}(g_{(2\,p)})a_{(3\,l)}a_{(3\,p)}a_2+\lambda_{
112 } (1-\xi)a_{(3\,p)}a_{ (2\,l)}.
\end{align*}

\noindent \eqref{item:corchete-1p,1l-2} We have, using q-Jacobi \eqref{eq:qjacob},
\begin{align*}
&[a_{(1\,l-1)},a_{(1\,l)}]_c=[[a_{(1\,l-2)},a_{l-1}]_c,a_{(1\,l)}]_c=[a_{(1\,l-2)},[a_{l-1},a_{(1\,l)}]_c]_c\\
&\,+\chi_{(1\,l)}(g_{l-1})\big([a_{(1\,l-2)},a_{(1\,l)}]_ca_{l-1}-\chi_{l-1}(g_{(1\,l)}g_{(1\,l-2)})a_{l-1}[a_{(1\,l-2)}
, a_ {(1\,l)}]_c\big).
\end{align*}
Now, by Lemma \ref{lem:tech1-lift} (4),
\begin{align*}
[a_{l-1},a_{(1\,l)}]_c&=-\chi_{(1\,l)}(g_{l-1})[a_{(1\,l)},a_{l-1}]_c\\
&=\lambda_{l-1l-1l}(1-\xi^2)\chi_{(1\,l)}(g_{l-1})a_{(1\,l-2)}g_{l-1l-1l}.
\end{align*}
Hence, $ [a_{(1\,l-2)},[a_{l-1},a_{(1\,l)}]_c]_c=-3\xi^2\chi_{(1\,l)}(g_{l-1})\lambda_{l-1l-1l}a_{(1\,l-2)}^2g_{l-1l-1l}.$

On the other hand, we have that, by item \eqref{item:corchete-1p,1l-1},
\begin{align*}
\chi_{l-1}(g_{(1\,l)}g_{(1\,l-2)})&a_{l-1}[a_{(1\,l-2)},a_{(1\,l)}]_c\\&=-3\xi^2\lambda_{112}\chi_{(1\,l)}(g_{(2\,l-2)
} )
\chi_{l-1}(g_{(1\,l)}g_{(1\,l-2)})a_{l-1}a_{(3\,l)}a_{(2\,l-2)}\\
&+3\lambda_{112}\chi_{1}(g_{(3\,l-2)})\chi_{l-1}(g_{(1\,l)}g_{(1\,l-2)})a_{l-1}a_{(3\,l-2)}a_{(2\,l)}.
\end{align*}
Now, by Lemma \ref{lem:tech1-lift}
\begin{align*}
 [a_{(1\,l-2)},&a_{(1\,l)}]_ca_{l-1}=-3\xi^2\lambda_{112}\chi_{(1\,l)}(g_{(2\,l-2)})
a_{(3\,l)}a_{(2\,l-2)}a_{l-1} \\
&+3\lambda_{112}\chi_{1}(g_{(3\,l-2)})a_{(3\,l-2)}a_{(2\,l)}a_{l-1}\\
&=-3\xi^2\lambda_{112}\chi_{(1\,l)}(g_{(2\,l-2)})
a_{(3\,l)}a_{(2\,l-1)} \\
&-3\xi^2\lambda_{112}\chi_{(1\,l)}(g_{(2\,l-2)})\chi_{l-1}(g_{(2\,l-2)})
a_{(3\,l)}a_{l-1}a_{(2\,l-2)} \\
&-3\lambda_{112}\lambda_{l-1l-1l}(1-\xi^2)\chi_{1}(g_{(3\,l-2)})a_{(3\,l-2)}a_{(2\,l-2)}g_{l-1l-1l}\\
&+3\lambda_{112}\chi_{1}(g_{(3\,l-2)})\chi_{l-1}(g_{(2\,l)})a_{(3\,l-2)}a_{l-1}a_{(2\,l)}\\
&=-3\xi^2\lambda_{112}\chi_{(1\,l)}(g_{(2\,l-2)})
a_{(3\,l)}a_{(2\,l-1)} \\
&+3\xi\lambda_{112}\lambda_{l-1l-1l}(1-\xi^2)\chi_{1}(g_{(2\,l-2)})
a_{(3\,l-2)}a_{(2\,l-2)}g_{l-1l-1l} \\
&-3\xi^2\lambda_{112}\chi_{(1\,l)}(g_{(2\,l-2)})\chi_{l-1}(g_{(2\,l-2)})
\chi_{l-1}(g_{(3\,l)})a_{l-1}a_{(3\,l)}a_{(2\,l-2)} \\
&-3\lambda_{112}\lambda_{l-1l-1l}(1-\xi^2)\chi_{1}(g_{(3\,l-2)})a_{(3\,l-2)}a_{(2\,l-2)}g_{l-1l-1l}\\
&+3\lambda_{112}\chi_{1}(g_{(3\,l-2)})\chi_{l-1}(g_{(2\,l)})a_{(3\,l-1)}a_{(2\,l)}\\
&+3\lambda_{112}\chi_{1}(g_{(3\,l-2)})\chi_{l-1}(g_{(2\,l)})\chi_{l-1}(g_{(3\,l-2)})
a_{l-1}a_{(3\,l-2)}a_{(2\,l)}
\end{align*}
Hence, using that $\chi_{l-1}(g_{112})=1$ and adding up the terms, we get (6).
\epf

Notice that for $0\leq p<q<l$, Lemmas \ref{lem:tech2b} \eqref{it8} and \ref{lem:tech1-gral} give:
\begin{multline}\label{eq:conm-lift-1}
 [y_{(p+1\,l)},y_{(q+1\,l)}]_c\\
=\begin{cases}
              0, & p+1<q< l;\\
-3\lambda_{p+1p+2p+2}
\chi_{p+2}(g_{(p+3\,l)})
y_{(p+3\,l)}^2, & p+1=q<l-1;\\
\lambda_{l-1ll}, & p+1=q=l-1.
             \end{cases}
\end{multline}

We fix $q=p+1$, $r=p+2$. Some of the identities computed in  Lemma 
\ref{lem:corchete-1p,1l-parte1} and Proposition \ref{pro:corchete-1p,1l-parte2} become 
simpler when multiplied by a factor $\lambda_{qrr}$, using Lemma \ref{lem:l112}. 
This will be of great importance in the computations, as 
$\varsigma_{(i\,l)}$ is a linear combination of the elements $\varsigma^{p}(\bs\lambda,\bs\mu)$ in \eqref{eqn:monomials-1} and each one of these terms is multiplied by $\lambda_{qrr}$.

We interpret these 
identities in the following corollary.

\begin{cor}\label{cor:p-p+1}
The following identities hold in $\mtL$.
\begin{enumerate}
\item If $p<3$, then
$$ \lambda_{qrr}[a_{(1\,p)},a_{(1\,p+1)}]_c= \lambda_{qrr}[a_{(1\,p)},a_{(1\,p+2)}]_c= \lambda_{qrr}[a_{(1\,p+1)},a_{(1\,p+2)}]_c=0.$$
 \item If $4\leq s=p+1,p+2$, then
\begin{multline*}
 \lambda_{qrr}[a_{(1\,p)},a_{(1\,s)}]_c=
-3\xi^2\lambda_{112} \lambda_{qrr}\chi_{(1\,s)}(g_{(2\,p)})
a_{(3\,s)}a_{(2\,p)}\\  +3\lambda_{112} \lambda_{qrr}\chi_{1}(g_{(3\,p)})a_{(3\,p)}a_{(2\,s)}.
\end{multline*}
\item
If $p\geq
3$, then
\begin{align*}
 \lambda_{qrr}[a_{(1\,q)},&a_{(1\,r)}]_c= -3 \lambda_{qrr}\xi^2\chi_{(1\,r)}(g_{q})\lambda_{qqr}a_{(1\,p)}^2g_{qqr}\\
&-3 \lambda_{qrr}\xi^2\lambda_{112}\chi_{(1\,r)}(g_{(2\,q)})
a_{(3\,r)}a_{(2\,q)} \\
&+3 \lambda_{qrr}\lambda_{112}\chi_{1}(g_{(3\,p)})\chi_{q}(g_{(2\,r)})\chi_{(1\,r)}(g_{q})a_{(3\,q)}a_{(2\,r)}.
\end{align*}
\end{enumerate}
\end{cor}

\pf
(1) follows using that $ \lambda_{233}\lambda_{112}= \lambda_{344}\lambda_{112}=0$.

(2) follows by Proposition \ref{pro:corchete-1p,1l-parte2} \eqref{item:corchete-1p,1l-1} using that $ \lambda_{qrr}\lambda_{ppp+1}=0$ by
Lemma \ref{lem:l112}. (3) is precisely Proposition \ref{pro:corchete-1p,1l-parte2} \eqref{item:corchete-1p,1l-2}.
\epf

In particular, Corollary \ref{cor:p-p+1} gives
\begin{cor}\label{cor:3q-3r}
Let $p\geq 3$. The following identities hold in $\mtL$.
\begin{enumerate}
 \item $
 \lambda_{qrr}\lambda_{112}[a_{(3\,q)},a_{(3\,r)}]_c=
-3 \lambda_{qrr}\lambda_{112}\lambda_{qqr}\xi^2\chi_{(3\,r)}(g_{q})a_{(3\,p)}^2g_{qqr}
.
$
\item $
 \lambda_{qrr}\lambda_{112}[a_{(3\,p)},a_{(3\,r)}]_c= \lambda_{qrr}\lambda_{112}[a_{(3\,p)},a_{(3\,q)}]_c=0.$
\end{enumerate}
 \end{cor}
\begin{cor}\label{cor:alpha-1-3} Let $p\geq 3$, $s=q,r$. The following identities hold in $\mtL$.
 \begin{enumerate}
  \item $\lambda_{112} \lambda_{qrr}[a_{(1\,q)},a_{(3\,p)}]_c=
\lambda_{112} \lambda_{qrr}[a_{(1\,r)}
,a_{(3\,p)}]_c=0$.
\item
$\lambda_{112} \lambda_{qrr}[a_{(1\,r)},a_{(3\,q)}]_c=
 3\lambda_{112} \lambda_{qrr}\xi^2\lambda_{qqr}\chi_{(3\,p)}(g_{(1\,r)})a_{(3\,p)}a_{(1\,p)}g_{qqr}$.
\item 
$\lambda_{112} \lambda_{qrr}[a_{(1\,q)},a_{(3\,r)}]_c=
\lambda_{112} \lambda_{qrr}\chi_{(3\,q)}(g_{(1\,q)})(1-\xi^2)a_{(3\,q)}a_{(1\,r)}$

 $\qquad \qquad \qquad \qquad -3\lambda_{112} \lambda_{qrr}\lambda_{qqr}\chi_{(3\,p)}(g_{(1\,q)})a_{(3\,p)}a_{(1\,p)}g_{qqr}$.
\item $\lambda_{112} \lambda_{qrr}[a_{(1\,p)},a_{(3\,s)}]_c=-\lambda_{112} \lambda_{qrr}\chi_{(3\,p)}(g_{12})(1-\xi)a_{(3\,p)}a_{(1\,r)}$.
 \end{enumerate}
\end{cor}
\pf
(1) and (2) follow from Lemma \ref{lem:corchete-1p,3l}. Also, q-Jacobi \eqref{eq:qjacob} and Lemma \ref{lem:corchete-1p,3l} together with Lemma 
\ref{lem:l112} give:
\begin{align*}
\lambda_{112}& \lambda_{qrr}[a_{(1\,q)},a_{(3\,r)}]_c=\lambda_{112} \lambda_{qrr}[a_{(1\,q)},[a_{(3\,q)},a_r]_c]_c\\
&=\lambda_{112} \lambda_{qrr}[[a_{(1\,q)},a_{(3\,q)}]_c
,a_r]_c
\\
&\quad -\lambda_{112} \lambda_{qrr}\chi_r(g_{(3\,q)})a_{(1\,r)}a_{(3\,q)}
+\lambda_{112} \lambda_{qrr}\chi_{(3\,q)}(g_{(1\,q)})a_{(3\,q)}a_{(1\,r)}\\
&=-\lambda_{112} \lambda_{qrr}\chi_r(g_{(3\,q)})[a_{(1\,r)},a_{(3\,q)}]_c\\
&\quad 
+\lambda_{112} \lambda_{qrr}\chi_{(3\,q)}(g_{(1\,q)})(1-\xi^2)a_{(3\,q)}a_{(1\,r)}.
\end{align*}
Hence now (3) follows using (2). (4) follows by Corollary \ref{cor:3q-3r} (2), using q-Jacobi.
\epf

\begin{cor}\label{cor:1-2}
Let $p\geq 3$. The following identities hold in $\mtL$.
\begin{enumerate}
 \item If $s=p,q,r$, then
\begin{align*}
 \lambda_{112} \lambda_{qrr}[a_{(1\,s)},a_{(2\,p)}]_c&=-3\xi^2\lambda_{112}\lambda_{122}
 \lambda_{qrr}\chi_2(g_{(3\,s)})a_{(3\,s)}a_{(3\,p)}.
\end{align*}
\item \begin{align*}
\lambda_{112} \lambda_{qrr}[a_{(1\,r)},a_{(2\,q)}]_c=
-3\xi^2&\lambda_{112}\lambda_{122} \lambda_{qrr}\chi_2(g_{(3\,r)})
a_{(3\,r)}a_{(3\,q)}\\
&
+3\lambda_{112} \lambda_{qrr}\xi^2\lambda_{qqr}\chi_{(2\,p)}(g_{(1\,r)})
a_{(2\,p)}a_{(1\,p)}g_{qqr}.
\end{align*}
\item
\begin{align*}
&\lambda_{112} \lambda_{qrr}[a_{(1\,q)},a_{(2\,r)}]_c=9\xi\lambda_{112}\lambda_{122}\lambda_{qqr} \lambda_{qrr}\chi_2(g_{(3\,q)})\chi_{(3\,r)}(g_{q})a_{(3\,p)}
^2g_{qqr}\\
&\qquad +3(1+\xi)\lambda_{112}\lambda_{122} \lambda_{qrr}\chi_{(2\,r)}(g_{(3\,q)})a_{(3\,r)}a_{(3\,q)}\\
&\qquad+\lambda_{112} \lambda_{qrr}\chi_{(2\,q)}(g_{(1\,q)})(1-\xi^2)a_{(2\,q)}a_{(1\,r)}\\
&\qquad -3\lambda_{112} \lambda_{qrr}\lambda_{qqr}\chi_{(2\,p)}(g_{(1\,q)})a_{(2\,p)}a_{(1\,p)}g_{qqr}.
\end{align*}
\item If $s=q,r$, then 
\begin{align*}
 \lambda_{112} \lambda_{qrr}[a_{(1\,p)},a_{(2\,s)}]_c=-3\xi&\lambda_{112}\lambda_{122}
 \lambda_{qrr}\chi_{2}(g_{(2\,p)})a_{(3\,s)}a_{(3\,p)}\\
& - \lambda_{112} \lambda_{qrr}\chi_{(2\,p)}(g_1)(1-\xi)a_{(2\,p)}a_{(1\,s)}.
\end{align*}
\end{enumerate}
\end{cor}
\pf
(1) follows by Lemmas \ref{lem:l112}, \ref{lem:tech1-lift} (1) and \ref{lem:corchete-1p,3l} and Corollary
\ref{cor:p-p+1} (1). We use
that $ \lambda_{qrr}\lambda_{p-1p-1p}= \lambda_{qrr}\lambda_{p-1pp}=0$.

(2) We have, using q-Jacobi \eqref{eq:qjacob}:
\begin{align*}
\lambda_{112}& \lambda_{qrr}[a_{(1\,r)},a_{(2\,q)}]_c=
\lambda_{112} \lambda_{qrr}[a_{(1\,r)},[a_2,a_{(3\,q)}]_c]_c\\
&=\lambda_{112} \lambda_{qrr}[[a_{(1\,r)},a_2]_c,a_{(3\,q)}]_c
+\lambda_{112} \lambda_{qrr}\chi_2(g_{(1\,r)})
\big(a_2[a_{(1\,r)},a_{(3\,q)}]_c\\
&\quad
-\chi_{(3\,q)}(g_2)\chi_{(1\,r)}(g_2)[a_{(1\,r)},a_{(3\,q)}]_ca_2\big).
\end{align*}
Now, by Lemma \ref{lem:tech1-lift} (1) and Corollary \ref{cor:3q-3r},
\begin{align*}
\lambda_{112}& \lambda_{qrr}[[a_{(1\,r)},a_2]_c,a_{(3\,q)}]_c\\
&=
\lambda_{112}\lambda_{122} \lambda_{qrr}(1-\xi^2)\chi_2(g_{(3\,r)})
\big(a_{(3\,r)}a_{(3\,q)}-\chi_{(3\,q)}(g_{(3\,r)}g_{122})
a_{(3\,q)}a_{(3\,r)}\big)\\
&=-3\xi^2\lambda_{112}\lambda_{122} \lambda_{qrr}\chi_2(g_{(3\,r)})
a_{(3\,r)}a_{(3\,q)}\\
&\quad -\lambda_{112}\lambda_{122} \lambda_{qrr}(1-\xi^2)\xi\chi_2(g_{(3\,r)})
\chi_{(3\,q)}(g_{(3\,r)})
[a_{(3\,q)},a_{(3\,r)}]_c\\
&=-3\xi^2\lambda_{112}\lambda_{122} \lambda_{qrr}\chi_2(g_{(3\,r)})
a_{(3\,r)}a_{(3\,q)}\\
&\quad +3\lambda_{112}\lambda_{122} \lambda_{qrr}\lambda_{qqr}(1-\xi^2)
\chi_{(2\,p)}(g_{(3\,r)})
a_{(3\,p)}^2g_{qqr}.
\end{align*}
Set $s=3\lambda_{112} \lambda_{qrr}\xi^2\lambda_{qqr}\chi_{(3\,p)}(g_{(1\,r)})$, by Corollary \ref{cor:alpha-1-3} (2):
\begin{align*}
\lambda_{112}& \lambda_{qrr}
\big(a_2[a_{(1\,r)},a_{(3\,q)}]_c
-\chi_{(3\,q)}(g_2)\chi_{(1\,r)}(g_2)[a_{(1\,r)},a_{(3\,q)}]_ca_2\big)\\
&=s
\big(a_2a_{(3\,p)}a_{(1\,p)}
-\chi_{(3\,p)}(g_2)\chi_{(1\,p)}(g_2)a_{(3\,p)}a_{(1\,p)}a_2\big)g_{qqr}\\
&=s
(a_{(2\,p)}a_{(1\,p)}-\lambda_{122}(\xi^2-\xi)\chi_{(1\,p)}(g_2)a_{(3\,p)}^2
)g_{qqr}\\
&\qquad +s
(\chi_{(3\,p)}(g_2)
-\chi_{(3\,p)}(g_2)\chi_{(1\,p)}(g_2)\chi_2(g_{(1\,p)}))a_{(3\,p)}a_2a_{(1\,p)}g_{qqr}\\
&=s
(a_{(2\,p)}a_{(1\,p)}+\lambda_{122}(1-\xi)\chi_{(3\,p)}(g_2)a_{(3\,p)}^2
)g_{qqr},
\end{align*}
using once again Lemma \ref{lem:tech1-lift} (1) and $1=\chi_{(1\,p)}(g_2)\chi_2(g_{(1\,p)})$. Adding up,
\begin{align*}
\lambda_{112} \lambda_{qrr}[a_{(1\,r)},a_{(2\,q)}]_c=&
-3\xi^2\lambda_{112}\lambda_{122} \lambda_{qrr}\chi_2(g_{(3\,r)})
a_{(3\,r)}a_{(3\,q)}\\
&\quad
+\chi_2(g_{(1\,r)})s
a_{(2\,p)}a_{(1\,p)}g_{qqr}.
\end{align*}
Here we have used that, as $\chi_{(3\,p)}(g_2)\chi_{(2\,p)}(g_{12})=\chi_{(3\,p)}(g_{122})\xi^2=1$:
\begin{multline*}
\lambda_{122}(1-\xi)\chi_2(g_{(1\,r)})\chi_{(3\,p)}(g_2)s\\
=
3\lambda_{112}\lambda_{122}\lambda_{qqr} \lambda_{qrr}\xi^2(1-\xi)\chi_2(g_{(1\,r)})\chi_{(3\,p)}(g_2)\chi_{(3\,p)}(g_{(1\,r)})\\
=-3\lambda_{112}\lambda_{122}\lambda_{qqr} \lambda_{qrr}(1-\xi^2)\chi_{(3\,p)}(g_2)\chi_{(2\,p)}(g_{12})\chi_{(2\,p)}(g_{(3\,r)})\\
=-3\lambda_{112}\lambda_{122}\lambda_{qqr} \lambda_{qrr}(1-\xi^2)\chi_{(2\,p)}(g_{(3\,r)}).
\end{multline*}
Hence the terms corresponding to $a_{(3\,p)}^2g_{qqr}$ cancel.

(3) We have, using q-Jacobi:
\begin{align*}
\lambda_{112} &\lambda_{qrr}[a_{(1\,q)},a_{(2\,r)}]_c=\lambda_{112} \lambda_{qrr}[[a_{(1\,q)},[a_2,a_{(3\,r)}]_c]_c\\
&=\lambda_{112} \lambda_{qrr}[[a_{(1\,q)},a_2]_c,a_{(3\,r)}]_c -\lambda_{112} \lambda_{qrr}\chi_{(3\,r)}(g_{2})[a_{(1\,q)},a_{(3\,r)}]_ca_2\\
&\qquad 
+\lambda_{112} \lambda_{qrr}\chi_{2}(g_{(1\,q)})a_2[a_{(1\,q)},a_{(3\,r)}]_c.
\end{align*}
Now, by Lemma \ref{lem:tech1-lift} and Corollary \ref{cor:3q-3r},
\begin{align*}
\lambda_{112}& \lambda_{qrr}[[a_{(1\,q)},a_2]_c,a_{(3\,r)}]_c=\lambda_{112}\lambda_{122} \lambda_{qrr}(1-\xi^2)\chi_2(g_{(3\,q)})[a_{(3\,q)},a_{(3\,r)}]_c\\
&+\lambda_{112}\lambda_{122} \lambda_{qrr}(1-\xi^2)\chi_2(g_{(3\,q)})\chi_{(3\,r)}(g_{(3\,q)})(1-\xi)a_{(3\,r)}a_{(3\,q)}\\
&=-3\xi^2(1-\xi^2)\lambda_{112}\lambda_{122}\lambda_{qqr} \lambda_{qrr}\chi_2(g_{(3\,q)})\chi_{(3\,r)}(g_{q})a_{(3\,p)}^2g_{qqr}\\
&\qquad +3\lambda_{112}\lambda_{122} \lambda_{qrr}\chi_{(2\,r)}(g_{(3\,q)})a_{(3\,r)}a_{(3\,q)}.
\end{align*}
On the other hand, by Corollary \ref{cor:alpha-1-3}, Lemma \ref{lem:tech1-lift} and Corollary \ref{cor:3q-3r},
\begin{align*}
\lambda_{112}& \lambda_{qrr}[a_{(1\,q)},a_{(3\,r)}]_ca_2=\\
&=\lambda_{112} \lambda_{qrr}\chi_{(3\,q)}(g_{(1\,q)})(1-\xi^2)a_{(3\,q)}a_{(1\,r)}a_2\\
&\quad -3\lambda_{112} \lambda_{qrr}\lambda_{qqr}\chi_{(3\,p)}(g_{(1\,q)})a_{(3\,p)}a_{(1\,p)}a_2g_{qqr}\\
&=\lambda_{112} \lambda_{qrr}\chi_{(3\,q)}(g_{(1\,q)})\chi_2(g_{(1\,r)})(1-\xi^2)a_{(3\,q)}a_2a_{(1\,r)}\\
&\quad +\lambda_{112}\lambda_{122} \lambda_{qrr}(1-\xi^2)^2\chi_{(3\,q)}(g_{(1\,q)})\chi_2(g_{(3\,r)})a_{(3\,q)}a_{(3\,r)}\\
&\quad  -3\lambda_{112} \lambda_{qrr}\lambda_{qqr}\chi_{(3\,p)}(g_{(1\,q)})\chi_2(g_{(1\,p)})
a_{(3\,p)}a_2a_{(1\,p)}g_{qqr}\\
&\quad  -3\lambda_{112}\lambda_{122} \lambda_{qrr}\lambda_{qqr}(1-\xi^2)\chi_{(3\,p)}
(g_{(1\,q)})\chi_2(g_{(3\,p)})a_{(3\,p)}^2g_{qqr}\\
&=\lambda_{112} \lambda_{qrr}\chi_{(3\,q)}(g_{(1\,q)})\chi_2(g_{(1\,r)})(1-\xi^2)a_{(3\,q)}a_2a_{(1\,r)}\\
&\quad -3\xi^2\lambda_{112}\lambda_{122} \lambda_{qrr}\chi_{(3\,q)}(g_{(1\,q)})\chi_2(g_{(3\,r)})\chi_{(3\,r)}(g_{(3\,q)})a_{(3\,r)}a_{(3\,q)}\\
&\quad -3\xi^2\lambda_{112}\lambda_{122}\lambda_{qqr}
 \lambda_{qrr}(1-\xi^2)^2\chi_{(3\,q)}(g_{(1\,q)})\chi_2(g_{(3\,r)})\chi_{(3\,r)}(g_q)a_{(3\,p)}^2g_{qqr}\\
&\quad  -3\lambda_{112} \lambda_{qrr}\lambda_{qqr}\chi_{(3\,p)}(g_{(1\,q)})\chi_2(g_{(1\,p)})
a_{(3\,p)}a_2a_{(1\,p)}g_{qqr}\\
&\quad  -3\lambda_{112}\lambda_{122} \lambda_{qrr}\lambda_{qqr}(1-\xi^2)\chi_{(3\,p)}
(g_{(1\,q)})\chi_2(g_{(3\,p)})a_{(3\,p)}^2g_{qqr}\\
&=-3\xi^2\lambda_{112}\lambda_{122} \lambda_{qrr}
\chi_{(3\,q)}(g_{(1\,q)})\chi_2(g_{(3\,r)})\chi_{(3\,r)}(g_{(3\,q)}) a_{(3\,r)}a_{(3\,q)}\\
&\quad  +(1-\xi^2)\lambda_{112} \lambda_{qrr}\chi_{(2\,q)}(g_{(1\,q)})\chi_2(g_{r})a_{(3\,q)}a_2a_{(1\,r)}\\
&\quad   +3(1-\xi)\lambda_{112}\lambda_{122}\lambda_{qqr} \lambda_{qrr}\chi_{(3\,p)}(g_{(1\,q)}) \chi_2(g_{(3\,p)})a_{(3\,p)}^2
 g_{qqr}\\
&\quad   -3\lambda_{112} \lambda_{qrr}\lambda_{qqr}\chi_{(3\,p)}(g_{(1\,q)})\chi_2(g_{(1\,p)})a_{(3\,p)}a_2a_{(1\,p)}g_{qqr}.
\end{align*}
Again, by  Corollary \ref{cor:alpha-1-3},
\begin{align*}
\lambda_{112}& \lambda_{qrr}a_2[a_{(1\,q)},a_{(3\,r)}]_c=
\lambda_{112} \lambda_{qrr}\chi_{(3\,q)}(g_{(1\,q)})(1-\xi^2)a_{(2\,q)}a_{(1\,r)}\\
&\quad+\lambda_{112} \lambda_{qrr}\chi_{(3\,q)}(g_{(1\,q)})\chi_{(3\,q)}(g_2)(1-\xi^2)a_{(3\,q)}a_2a_{(1\,r)}\\
&\quad -3\lambda_{112} \lambda_{qrr}\lambda_{qqr}\chi_{(3\,p)}(g_{(1\,q)})a_{(2\,p)}a_{(1\,p)}g_{qqr}\\
&\quad -3\lambda_{112} \lambda_{qrr}\lambda_{qqr}\chi_{(3\,p)}(g_{(1\,q)})\chi_{(3\,p)}(g_2)a_{(3\,p)}a_2a_{(1\,p)}g_{qqr}.
\end{align*}
Adding up,
\begin{align*}
&\lambda_{112} \lambda_{qrr}[a_{(1\,q)},a_{(2\,r)}]_c=9\xi\lambda_{112}\lambda_{122}\lambda_{qqr} \lambda_{qrr}\chi_2(g_{(3\,q)})\chi_{(3\,r)}(g_{q})a_{(3\,p)}
^2g_{qqr}\\
&\qquad +3(1+\xi)\lambda_{112}\lambda_{122} \lambda_{qrr}\chi_{(2\,r)}(g_{(3\,q)})a_{(3\,r)}a_{(3\,q)}\\
&\qquad+\lambda_{112} \lambda_{qrr}\chi_{(2\,q)}(g_{(1\,q)})(1-\xi^2)a_{(2\,q)}a_{(1\,r)}\\
&\qquad -3\lambda_{112} \lambda_{qrr}\lambda_{qqr}\chi_{(2\,p)}(g_{(1\,q)})a_{(2\,p)}a_{(1\,p)}g_{qqr}.
\end{align*}
(4) follows using (1) together with q-Jacobi and Lemma \ref{lem:l112}.
\epf 

Next, we order the elements \eqref{eqn:apqr} in terms of the PBW basis. We have:
\begin{pro}\label{pro:zpqr}
Assume $p\geq 3$. Then
\begin{align}
\label{eqn:zpqr} \lambda_{qrr}a_{p,q,r}&= \lambda_{qrr}\xi\chi_{q}(g_{(1\,p)})
a_{(1\,r)}a_{(1\,q)}a_{(1\,p)}\\
\notag &
-3\xi^2 \lambda_{qrr}\lambda_{qqr}\chi_{(1\,r)}(g_{q})
a_{(1\,p)}^3g_{qqr}\\
\notag &
-3\xi^2 \lambda_{qrr}\lambda_{112}\chi_{qqr}(g_{(1\,p)})
\chi_{(1\,r)}(g_{(2\,q)})
a_{(3\,r)}a_{(2\,q)} a_{(1\,p)}\\
\notag &
+3 \lambda_{qrr}\lambda_{112}\chi_{qqr}(g_{(1\,p)})
\chi_{1}(g_{(3\,q)})
a_{(3\,q)}a_{(2\,r)}a_{(1\,p)}\\
\notag & -3\xi^2 \lambda_{qrr}\lambda_{112}
\chi_{qqr}(g_{(1\,p)})\chi_{(3\,r)}(g_{(2\,q)})\chi_{(2\,p)} (g_{(1\,q)})
a_{(3\,r)}a_{(2\,p)} a_{(1\,q)}\\
\notag & +3 \lambda_{qrr}
\lambda_{112}
\chi_{(3\,p)}(g_{12})
\chi_{r}(g_{(1\,q)})
a_{(3\,p)}a_{(2\,r)} a_{(1\,q)}\\
\notag & +3 \lambda_{qrr}
\lambda_{112}(1-\xi^2)\chi_{(3\,p)}(g_{12})
a_{(3\,p)}a_{(2\,q)}a_{(1\,r)}\\
\notag &+6 \lambda_{qrr}\lambda_{112}\chi_{(1\,q)}(g_{(2\,p)})
a_{(3\,q)}a_{(2\,p)}a_{(1\,r)}.
\end{align}
\begin{align}
\label{eqn:zprq} \lambda_{qrr}a_{p,r,q}&= \lambda_{qrr}\xi^2\chi_{qqr} (g_{(1\,p)})
a_{(1\,r)}a_{(1\,q)}a_{(1\,p)}\\
\notag 
&-3\xi^2\lambda_{112} \lambda_{qrr}\chi_{(1\,q)}(g_{(1\,p)})\chi_{(3\,q)}(g_{(1\,r)}
)\chi_{r}(g_{1})
a_{(3\,q)}a_{(2\,p)}a_{(1\,r)} \\
\notag &+3\lambda_{112} \lambda_{qrr}\chi_{1}(g_{(3\,p)})\chi_{(1\,q)}(g_{r})
a_{(3\,p)}a_{(2\,q)}a_{(1\,r)}\\
\notag &-3\xi^2\lambda_{112} \lambda_{qrr}\chi_{(1\,r)}(g_{(2\,p)})
a_{(3\,r)}a_{(2\,p)}a_{(1\,q)}\\
\notag &+3\lambda_{112} \lambda_{qrr}\chi_{1}(g_{(3\,p)})
a_{(3\,p)}a_{(2\,r)}a_{(1\,q)}.
\end{align}
\begin{align}
\label{eqn:zrpq} \lambda_{qrr}a_{r,p,q}&=\xi\chi_{q}(g_{(1\,p)}) \lambda_{qrr}
a_{(1\,r)}a_{(1\,q)}a_{(1\,p)}\\
\notag&-3\xi^2\lambda_{112} \lambda_{qrr}\chi_{(3\,q)}(g_{(1\,r)})\chi_{(2\,p)}(g_{r})
a_{(3\,q)}a_{(2\,p)}a_{(1\,r)} \\
\notag 
&+3\lambda_{112} \lambda_{qrr}\chi_{(3\,p)}(g_{(2\,r)})\chi_{(2\,q)}(g_{(1\,r)})
a_{(3\,p)}a_{(2\,q)}a_{(1\,r)}.
\end{align}
\end{pro}
\pf
We will go through the description of $ \lambda_{qrr}a_{p,q,r}$ step by step, 
following the
identities in the lemmas. The other two summands are simpler and 
will be presented in their final form.  Every time there is a monomial that 
needs to be ordered, we shall highlight it on bold letters, for the reader to 
identify which is the bracket that needs to be computed for the next step. 
To do this, we shall use Corollaries \ref{cor:p-p+1}, \ref{cor:alpha-1-3}, 
\ref{cor:1-2} and \ref{cor:3q-3r}. 

We shall also reduce some of the scalars, for instance we consider:
$$ \lambda_{qrr}\chi_{(1\,q)}(g_{(1\,p)})\chi_{(1\,r)}(g_{(1\,p)}g_{(1\,q)})=
 \lambda_{qrr}\chi_{r}(g_{(1\,q)}).
$$
However, we leave a full reduction to the end.


We have, using Corollary \ref{cor:p-p+1}:
\begin{align*} \lambda_{qrr}&a_{p,q,r}= \lambda_{qrr}\chi_{(1\,q)}(g_{(1\,p)})\chi_{
(1\,r)}(g_{(1\,p)}g_{(1\,q)})a_{(1\,r)}a_{(1\,q)}a_{(1\,p)}\\
&
+\chi_{(1\,q)}(g_{(1\,p)})\chi_{(1\,r)}(g_{(1\,p)}) \lambda_{qrr}[a_{(1\,q)},a_{(1\,r)
}]_ca_{(1\,p)}\\
& +\chi_{(1\,q)}(g_{(1\,p)}) \lambda_{qrr}a_{(1\,q)}
[a_{(1\,p)},a_{(1\,r)}]_c+ \lambda_{qrr}[a_{(1\,p)},a_{(1\,q)}]_ca_{(1\,r)}\\
&= \lambda_{qrr}\xi\chi_{q}(g_{(1\,p)})
a_{(1\,r)}a_{(1\,q)}a_{(1\,p)}\\
&
-\chi_{(1\,q)}(g_{(1\,p)})\chi_{(1\,r)}(g_{(1\,p)})
3 \lambda_{qrr}\xi^2\chi_{(1\,r)}(g_{q})\lambda_{qqr}\chi_{(1\,p)}(g_{qqr})
a_{(1\,p)}^3g_{qqr}\\
&
-\chi_{(1\,q)}(g_{(1\,p)})\chi_{(1\,r)}(g_{(1\,p)})
3 \lambda_{qrr}\xi^2\lambda_{112}\chi_{(1\,r)}(g_{(2\,q)})
a_{(3\,r)}a_{(2\,q)} a_{(1\,p)}\\
&
+\chi_{(1\,q)}(g_{(1\,p)})\chi_{(1\,r)}(g_{(1\,p)})
3 \lambda_{qrr}\lambda_{112}\chi_{1}(g_{(3\,p)})\chi_{q}(g_{(2\,r)})\chi_{(1\,r)}(g_{q
})\\
&\hspace*{9cm}
a_{(3\,q)}a_{(2\,r)}a_{(1\,p)}\\
& -\chi_{(1\,q)}(g_{(1\,p)}) \lambda_{qrr}
3\xi^2\lambda_{112}\chi_{(1\,r)}(g_{(2\,p)})
{\bf a_{(1\,q)}a_{(3\,r)}}a_{(2\,p)}\\
& +\chi_{(1\,q)}(g_{(1\,p)}) \lambda_{qrr}
3\lambda_{112}\chi_{1}(g_{(3\,p)}){\bf a_{(1\,q)}a_{(3\,p)}}a_{(2\,r)}\\
&-3\xi^2\lambda_{112} \lambda_{qrr}\chi_{(1\,q)}(g_{(2\,p)})
a_{(3\,q)}a_{(2\,p)}a_{(1\,r)}.
\end{align*}
Using $\chi_{(1\,q)}(g_{(1\,p)})=\xi\chi_{q}(g_{(1\,p)})$, and
\begin{align*}
\lambda_{qqr}\chi_{(1\,p)}(g_{qqr})&=\lambda_{qqr}\xi,\\
\lambda_{qqr}\chi_{(1\,q)}(g_{(1\,p)})\chi_{(1\,r)}(g_{(1\,p)})&=\lambda_{qqr}
\xi^2\\
 \lambda_{qrr}\lambda_{112}\chi_{(1\,q)}(g_{(1\,p)})\chi_{(1\,r)}(g_{(1\,p)})
&=
 \lambda_{qrr}\lambda_{112}\chi_{qqr}(g_{(1\,p)}),\\
\chi_{1}(g_{(3\,p)})\chi_{q}(g_{(2\,r)})\chi_{(1\,r)}(g_{q})&=\chi_{1}(g_{(3\,q)
})\\
\chi_{(1\,q)}(g_{(1\,p)})\chi_{(1\,r)}(g_{(2\,p)})&=\xi\chi_{qqr}(g_{(1\,p)}
)\chi_1(g_{(1\,r)}),
\end{align*}
this becomes, applying Corollary \ref{cor:alpha-1-3}
\begin{align*}
 \lambda_{qrr}a_{p,q,r}&= \lambda_{qrr}\chi_{r}(g_{(1\,q)})
a_{(1\,r)}a_{(1\,q)}a_{(1\,p)}\\
&
-3\xi^2 \lambda_{qrr}\lambda_{qqr}\chi_{(1\,r)}(g_{q})
a_{(1\,p)}^3g_{qqr}\\
&
-3\xi^2 \lambda_{qrr}\lambda_{112}\chi_{qqr}(g_{(1\,p)})
\chi_{(1\,r)}(g_{(2\,q)})
a_{(3\,r)}a_{(2\,q)} a_{(1\,p)}\\
&
+3 \lambda_{qrr}\lambda_{112}\chi_{qqr}(g_{(1\,p)})
\chi_{1}(g_{(3\,q)})
a_{(3\,q)}a_{(2\,r)}a_{(1\,p)}\\
& -3 \lambda_{qrr}\lambda_{112}
\chi_{qqr}(g_{(1\,p)})\chi_1(g_{(1\,r)})\chi_{(3\,r)}(g_{(1\,q)})
a_{(3\,r)}{\bf a_{(1\,q)}a_{(2\,p)}}\\
& -3 \lambda_{qrr}\lambda_{112}
\chi_{qqr}(g_{(1\,p)})\chi_1(g_{(1\,r)})
[ a_{(1\,q)},a_{(3\,r)}]_ca_{(2\,p)}\\
& +3 \lambda_{qrr}
\lambda_{112}\chi_{(1\,q)}(g_{(1\,p)})\chi_{1}(g_{(3\,p)})\chi_{(3\,p)}(g_{(1\,
q)}) a_{(3\,p)}{\bf a_{(1\,q)}a_{(2\,r)}}\\
& +3 \lambda_{qrr}
\lambda_{112}\chi_{(1\,q)}(g_{(1\,p)})\chi_{1}(g_{(3\,p)})[ 
a_{(1\,q)},a_{(3\,p)}]_ca_{(2\,r)}\\
&-3\xi^2 \lambda_{qrr}\lambda_{112}\chi_{(1\,q)}(g_{(2\,p)})
a_{(3\,q)}a_{(2\,p)}a_{(1\,r)}\\
&= \lambda_{qrr}\chi_{r}(g_{(1\,q)})
a_{(1\,r)}a_{(1\,q)}a_{(1\,p)}\\
&
-3\xi^2 \lambda_{qrr}\lambda_{qqr}\chi_{(1\,r)}(g_{q})
a_{(1\,p)}^3g_{qqr}\\
&
-3\xi^2 \lambda_{qrr}\lambda_{112}\chi_{qqr}(g_{(1\,p)})
\chi_{(1\,r)}(g_{(2\,q)})
a_{(3\,r)}a_{(2\,q)} a_{(1\,p)}\\
&
+3 \lambda_{qrr}\lambda_{112}\chi_{qqr}(g_{(1\,p)})
\chi_{1}(g_{(3\,q)})
a_{(3\,q)}a_{(2\,r)}a_{(1\,p)}\\
& -3 \lambda_{qrr}\lambda_{112}
\chi_{qqr}(g_{(1\,p)})\chi_1(g_{(1\,r)})\chi_{(3\,r)}(g_{(1\,q)})
a_{(3\,r)}{\bf a_{(1\,q)}a_{(2\,p)}}\\
& -3(1-\xi^2) \lambda_{qrr}\lambda_{112}
\chi_{qqr}(g_{(1\,p)})\chi_1(g_{(1\,r)})
\chi_{(3\,q)}(g_{(1\,q)})a_{(3\,q)}{\bf a_{(1\,r)}a_{(2\,p)}}\\
& +9\xi \lambda_{qrr}\lambda_{112}\lambda_{qqr}
\chi_1(g_{(1\,r)})
\chi_{(3\,p)}(g_{(1\,q)})a_{(3\,p)}{\bf a_{(1\,p)}a_{(2\,p)}}g_{qqr}\\
& +3 \lambda_{qrr}
\lambda_{112}\chi_{(1\,q)}(g_{(1\,p)})\chi_{(3\,p)}(g_{(2\,q)}) a_{(3\,p)}{\bf 
a_{(1\,q)}a_{(2\,r)}}\\
&-3\xi^2 \lambda_{qrr}\lambda_{112}\chi_{(1\,q)}(g_{(2\,p)})
a_{(3\,q)}a_{(2\,p)}a_{(1\,r)}.
\end{align*}
We have also used 
$$\lambda_{qqr}\chi_{qqr}(g_{(1\,p)})\chi_{(2\,p)}(g_{qqr})=\lambda_{qqr}\xi, 
\chi_{1}(g_{(3\,p)})\chi_{(3\,p)}(g_{(1\,q)})=\chi_{(3\,p)}(g_{(2\,q)}).$$
We obtain: 
\begin{align*}
 \lambda_{qrr}&a_{p,q,r}= \lambda_{qrr}\chi_{r}(g_{(1\,q)})
a_{(1\,r)}a_{(1\,q)}a_{(1\,p)}\\
&
-3\xi^2 \lambda_{qrr}\lambda_{qqr}\chi_{(1\,r)}(g_{q})
a_{(1\,p)}^3g_{qqr}\\
&
-3\xi^2 \lambda_{qrr}\lambda_{112}\chi_{qqr}(g_{(1\,p)})
\chi_{(1\,r)}(g_{(2\,q)})
a_{(3\,r)}a_{(2\,q)} a_{(1\,p)}\\
&
+3 \lambda_{qrr}\lambda_{112}\chi_{qqr}(g_{(1\,p)})
\chi_{1}(g_{(3\,q)})
a_{(3\,q)}a_{(2\,r)}a_{(1\,p)}\\
& -3 \lambda_{qrr}\lambda_{112}
\chi_{qqr}(g_{(1\,p)})\chi_1(g_{(1\,r)})\chi_{(3\,r)}(g_{(1\,q)})
\chi_{(2\,p)} (g_{(1\,q)})
a_{(3\,r)}a_{(2\,p)} a_{(1\,q)}\\
& -3 \lambda_{qrr}\lambda_{112}
\chi_{qqr}(g_{(1\,p)})\chi_1(g_{(1\,r)})\chi_{(3\,r)}(g_{(1\,q)})
a_{(3\,r)}[ a_{(1\,q)},a_{(2\,p)}]_c\\
& -3(1-\xi^2) \lambda_{qrr}\lambda_{112}
\chi_{qqr}(g_{(1\,p)})\chi_{1} (g_{(1\,r)})
\chi_{(3\,q)}(g_{(1\,q)})\chi_{(2\,p)} (g_{(1\,r)})\\
&\hspace*{9cm}
a_{(3\,q)}a_{(2\,p)} a_{(1\,r)}\\
& -3(1-\xi^2) \lambda_{qrr}\lambda_{112}
\chi_{qqr}(g_{(1\,p)})\chi_1(g_{(1\,r)})
\chi_{(3\,q)}(g_{(1\,q)})
a_{(3\,q)}[ a_{(1\,r)},a_{(2\,p)}]_c\\
& +9\xi \lambda_{qrr}\lambda_{112}\lambda_{qqr}
\chi_1(g_{(1\,r)}) \chi_{(2\,p)}(g_{(1\,p)})
\chi_{(3\,p)}(g_{(1\,q)})a_{(3\,p)} a_{(2\,p)}a_{(1\,p)}g_{qqr}\\
& +9\xi \lambda_{qrr}\lambda_{112}\lambda_{qqr}
\chi_1(g_{(1\,r)})
\chi_{(3\,p)}(g_{(1\,q)})a_{(3\,p)}[ a_{(1\,p)},a_{(2\,p)}]_cg_{qqr}\\
& +3 \lambda_{qrr}
\lambda_{112}\chi_{(1\,q)}(g_{(1\,p)})\chi_{(3\,p)}(g_{(2\,q)}) 
\chi_{(2\,r)}(g_{(1\,q)})
a_{(3\,p)}a_{(2\,r)} a_{(1\,q)}\\
& +3 \lambda_{qrr}
\lambda_{112}\chi_{(1\,q)}(g_{(1\,p)})\chi_{(3\,p)}(g_{(2\,q)}) 
a_{(3\,p)}[a_{(1\,q)},a_{(2\,r)}]_c\\
&-3\xi^2 \lambda_{qrr}\lambda_{112}\chi_{(1\,q)}(g_{(2\,p)})
a_{(3\,q)}a_{(2\,p)}a_{(1\,r)}.
\end{align*}
Observe that $\chi_1(g_{(1\,r)})\chi_{(2\,p)} (g_{(1\,r)})=\chi_{(1\,p)} 
(g_{(1\,r)})$ and that we can add the two terms $a_{(3\,q)}a_{(2\,p)}a_{(1\,r)}$ 
and the corresponding scalar becomes:
\begin{align*}
&-3(1-\xi^2) \lambda_{qrr}\lambda_{112}
\chi_{qqr}(g_{(1\,p)})
\chi_{(3\,q)}(g_{(1\,q)})\chi_{(1\,p)} (g_{(1\,r)})\\
&\quad -3\xi^2 \lambda_{qrr}\lambda_{112}\chi_{(1\,q)}(g_{(2\,p)})
=-3 \lambda_{qrr}\lambda_{112}\chi_{(1\,q)}(g_{(2\,p)})\xi^2\times\\
&\quad \times\Big(
(1-\xi^2)\xi
\chi_{qqr}(g_{(1\,p)})
\chi_{(3\,q)}(g_{(1\,q)})\chi_{(1\,p)} (g_{(1\,r)})\chi_{(2\,p)}(g_{(1\,q)})
+1\Big)\\
&=-3 \lambda_{qrr}\lambda_{112}\chi_{(1\,q)}(g_{(2\,p)})\xi^2
\Big(
(1-\xi^2)\xi^2
+1\Big)=6 \lambda_{qrr}\lambda_{112}\chi_{(1\,q)}(g_{(2\,p)}).
\end{align*}

Now, we apply Corollary \ref{cor:1-2}:
\begin{align*}
 \lambda_{qrr}&a_{p,q,r}= \lambda_{qrr}\chi_{r}(g_{(1\,q)})
a_{(1\,r)}a_{(1\,q)}a_{(1\,p)}\\
&
-3\xi^2 \lambda_{qrr}\lambda_{qqr}\chi_{(1\,r)}(g_{q})
a_{(1\,p)}^3g_{qqr}\\
&
-3\xi^2 \lambda_{qrr}\lambda_{112}\chi_{qqr}(g_{(1\,p)})
\chi_{(1\,r)}(g_{(2\,q)})
a_{(3\,r)}a_{(2\,q)} a_{(1\,p)}\\
&
+3 \lambda_{qrr}\lambda_{112}\chi_{qqr}(g_{(1\,p)})
\chi_{1}(g_{(3\,q)})
a_{(3\,q)}a_{(2\,r)}a_{(1\,p)}\\
& -3 \lambda_{qrr}\lambda_{112}
\chi_{qqr}(g_{(1\,p)})\chi_1(g_{(1\,r)})\chi_{(3\,r)}(g_{(1\,q)})
\chi_{(2\,p)} (g_{(1\,q)})
a_{(3\,r)}a_{(2\,p)} a_{(1\,q)}\\
& +9\xi^2 \lambda_{qrr}\lambda_{112}\lambda_{122}
\chi_{qqr}(g_{(1\,p)})\chi_1(g_{(1\,r)})\chi_{(3\,r)}(g_{(1\,q)})\chi_2(g_{(3\,
q)})\\
&\hspace*{8cm}
a_{(3\,r)}a_{(3\,q)}a_{(3\,p)}\\
& +9\xi^2(1-\xi^2) \lambda_{qrr}\lambda_{112}\lambda_{122}
\chi_{qqr}(g_{(1\,p)})\chi_1(g_{(1\,r)})
\chi_{(3\,q)}(g_{(1\,q)})\chi_2(g_{(3\,r)})\\
&\hspace*{8cm}
{\bf a_{(3\,q)}a_{(3\,r)}}a_{(3\,p)}\\
& +9\xi \lambda_{qrr}\lambda_{112}\lambda_{qqr}
\chi_1(g_{(1\,r)}) \chi_{(2\,p)}(g_{(1\,p)})
\chi_{(3\,p)}(g_{(1\,q)})a_{(3\,p)} a_{(2\,p)}a_{(1\,p)}g_{qqr}\\
& -27 \lambda_{qrr}\lambda_{112}\lambda_{122}\lambda_{qqr}
\chi_1(g_{(1\,r)})
\chi_{(3\,p)}(g_{(1\,q)})\chi_2(g_{(3\,p)})
a_{(3\,p)}^3g_{qqr}\\
& +3 \lambda_{qrr}
\lambda_{112}\chi_{(1\,q)}(g_{(1\,p)})\chi_{(3\,p)}(g_{(2\,q)}) 
\chi_{(2\,r)}(g_{(1\,q)})
a_{(3\,p)}a_{(2\,r)} a_{(1\,q)}\\
& +27\xi \lambda_{qrr}
\lambda_{112}\lambda_{122}\lambda_{qqr}\chi_{(1\,q)}(g_{(1\,p)})\chi_{(3\,p)}(g_
{(2\,q)}) 
\chi_2(g_{(3\,q)})\chi_{(3\,r)}(g_{q})\\
&\hspace*{8cm}
a_{(3\,p)}^3g_{qqr}\\
& +9(1+\xi) \lambda_{qrr}
\lambda_{112}\lambda_{122}\chi_{(1\,q)}(g_{(1\,p)})\chi_{(3\,p)}(g_{(2\,q)}) 
\chi_{(2\,r)}(g_{(3\,q)})\\
&\hspace*{8cm}
{\bf a_{(3\,p)}a_{(3\,r)}}a_{(3\,q)}\\
& +3 \lambda_{qrr}
\lambda_{112}(1-\xi^2)\chi_{(1\,q)}(g_{(1\,p)})\chi_{(3\,p)}(g_{(2\,q)}) 
\chi_{(2\,q)}(g_{(1\,q)})
a_{(3\,p)}a_{(2\,q)}a_{(1\,r)}\\
& -9 \lambda_{qrr}
\lambda_{112}\lambda_{qqr}\chi_{(1\,q)}(g_{(1\,p)})\chi_{(3\,p)}(g_{(2\,q)}) 
\chi_{(2\,p)}(g_{(1\,q)})
a_{(3\,p)}a_{(2\,p)}a_{(1\,p)}g_{qqr}\\
&+6 \lambda_{qrr}\lambda_{112}\chi_{(1\,q)}(g_{(2\,p)})
a_{(3\,q)}a_{(2\,p)}a_{(1\,r)}.
\end{align*}
Observe that the terms $a_{(3\,p)}^3g_{qqr}$ have the scalar
\begin{align*}
27(1-\xi^2) \lambda_{qrr}\lambda_{112}\lambda_{122}\lambda_{qqr}
\chi_{(3\,p)}(g_{q}) \chi_{q}(g_{1})
\end{align*}
On the other hand, the terms $a_{(3\,p)}a_{(2\,p)}a_{(1\,p)}g_{qqr}$ cancel with 
each other.

We order the terms $a_{(3\,*)}a_{(3\,*)}a_{(3\,*)}$  using Corollary 
\ref{cor:3q-3r} and we get: 
\begin{align*}
 \lambda_{qrr}&a_{p,q,r}= \lambda_{qrr}\chi_{r}(g_{(1\,q)})
a_{(1\,r)}a_{(1\,q)}a_{(1\,p)}\\
&
-3\xi^2 \lambda_{qrr}\lambda_{qqr}\chi_{(1\,r)}(g_{q})
a_{(1\,p)}^3g_{qqr}\\
&
-3\xi^2 \lambda_{qrr}\lambda_{112}\chi_{qqr}(g_{(1\,p)})
\chi_{(1\,r)}(g_{(2\,q)})
a_{(3\,r)}a_{(2\,q)} a_{(1\,p)}\\
&
+3 \lambda_{qrr}\lambda_{112}\chi_{qqr}(g_{(1\,p)})
\chi_{1}(g_{(3\,q)})
a_{(3\,q)}a_{(2\,r)}a_{(1\,p)}\\
& -3 \lambda_{qrr}\lambda_{112}
\chi_{qqr}(g_{(1\,p)})\chi_1(g_{(1\,r)})\chi_{(3\,r)}(g_{(1\,q)})
\chi_{(2\,p)} (g_{(1\,q)})
a_{(3\,r)}a_{(2\,p)} a_{(1\,q)}\\
& +9\xi^2 \lambda_{qrr}\lambda_{112}\lambda_{122}
\chi_{qqr}(g_{(1\,p)})\chi_1(g_{(1\,r)})\chi_{(3\,r)}(g_{(1\,q)})\chi_2(g_{(3\,
q)})\\
&\hspace*{9cm}
a_{(3\,r)}a_{(3\,q)}a_{(3\,p)}\\
& +9\xi^2(1-\xi^2) \lambda_{qrr}\lambda_{112}\lambda_{122}
\chi_{qqr}(g_{(1\,p)})\chi_1(g_{(1\,r)})
\chi_{(3\,q)}(g_{(1\,q)})\chi_2(g_{(3\,r)})\\
&\hspace*{7cm}
\chi_{(3\,r)}(g_{(3\,q)})
 a_{(3\,r)}a_{(3\,q)}a_{(3\,p)}\\
& +9\xi^2(1-\xi^2) \lambda_{qrr}\lambda_{112}\lambda_{122}
\chi_{qqr}(g_{(1\,p)})\chi_1(g_{(1\,r)})
\chi_{(3\,q)}(g_{(1\,q)})\chi_2(g_{(3\,r)})\\
&\hspace*{8cm}
 [a_{(3\,q)},a_{(3\,r)}]_ca_{(3\,p)}\\
&+ 27(1-\xi^2) \lambda_{qrr}\lambda_{112}\lambda_{122}\lambda_{qqr}
\chi_{(3\,p)}(g_{q}) \chi_{q}(g_{1})
a_{(3\,p)}^3g_{qqr}\\
& +3 \lambda_{qrr}
\lambda_{112}\chi_{(1\,q)}(g_{(1\,p)})\chi_{(3\,p)}(g_{(2\,q)}) 
\chi_{(2\,r)}(g_{(1\,q)})
a_{(3\,p)}a_{(2\,r)} a_{(1\,q)}\\
& +9(1+\xi) \lambda_{qrr}
\lambda_{112}\lambda_{122}\chi_{(1\,q)}(g_{(1\,p)})\chi_{(3\,p)}(g_{(2\,q)}) 
\chi_{(2\,r)}(g_{(3\,q)})
\chi_{(3\,r)}(g_{(3\,p)})\\
&\hspace*{7cm}
\chi_{(3\,q)}(g_{(3\,p)})a_{(3\,r)} a_{(3\,q)}a_{(3\,p)}\\
& +9(1+\xi) \lambda_{qrr}
\lambda_{112}\lambda_{122}\chi_{(1\,q)}(g_{(1\,p)})\chi_{(3\,p)}(g_{(2\,q)}) 
\chi_{(2\,r)}(g_{(3\,q)})
\chi_{(3\,r)}(g_{(3\,p)})\\
&\hspace*{8cm}
a_{(3\,r)} [a_{(3\,p)},a_{(3\,q)}]_c\\
& +9(1+\xi) \lambda_{qrr}
\lambda_{112}\lambda_{122}\chi_{(1\,q)}(g_{(1\,p)})\chi_{(3\,p)}(g_{(2\,q)}) 
\chi_{(2\,r)}(g_{(3\,q)})\\
&\hspace*{8cm}
[ a_{(3\,p)},a_{(3\,r)}]_ca_{(3\,q)}\\
& +3 \lambda_{qrr}
\lambda_{112}(1-\xi^2)\chi_{(1\,q)}(g_{(1\,p)})\chi_{(3\,p)}(g_{(2\,q)}) 
\chi_{(2\,q)}(g_{(1\,q)})
a_{(3\,p)}a_{(2\,q)}a_{(1\,r)}\\
&+6 \lambda_{qrr}\lambda_{112}\chi_{(1\,q)}(g_{(2\,p)})
a_{(3\,q)}a_{(2\,p)}a_{(1\,r)}.
\end{align*}
That is,
\begin{align*}
 \lambda_{qrr}&a_{p,q,r}= \lambda_{qrr}\chi_{r}(g_{(1\,q)})
a_{(1\,r)}a_{(1\,q)}a_{(1\,p)}\\
&
-3\xi^2 \lambda_{qrr}\lambda_{qqr}\chi_{(1\,r)}(g_{q})
a_{(1\,p)}^3g_{qqr}\\
&
-3\xi^2 \lambda_{qrr}\lambda_{112}\chi_{qqr}(g_{(1\,p)})
\chi_{(1\,r)}(g_{(2\,q)})
a_{(3\,r)}a_{(2\,q)} a_{(1\,p)}\\
&
+3 \lambda_{qrr}\lambda_{112}\chi_{qqr}(g_{(1\,p)})
\chi_{1}(g_{(3\,q)})
a_{(3\,q)}a_{(2\,r)}a_{(1\,p)}\\
& -3 \lambda_{qrr}\lambda_{112}
\chi_{qqr}(g_{(1\,p)})\chi_1(g_{(1\,r)})\chi_{(3\,r)}(g_{(1\,q)})
\chi_{(2\,p)} (g_{(1\,q)})
a_{(3\,r)}a_{(2\,p)} a_{(1\,q)}\\
& +9\xi^2 \lambda_{qrr}\lambda_{112}\lambda_{122}
\chi_{qqr}(g_{(1\,p)})\chi_1(g_{(1\,r)})\chi_{(3\,r)}(g_{(1\,q)})\chi_2(g_{(3\,
q)})\\
&\hspace*{8cm}
a_{(3\,r)}a_{(3\,q)}a_{(3\,p)}\\
& +9\xi^2(1-\xi^2) \lambda_{qrr}\lambda_{112}\lambda_{122}
\chi_{qqr}(g_{(1\,p)})\chi_1(g_{(1\,r)})
\chi_{(3\,q)}(g_{(1\,q)})\chi_2(g_{(3\,r)})\\
&\hspace*{7cm}
\chi_{(3\,r)}(g_{(3\,q)})
 a_{(3\,r)}a_{(3\,q)}a_{(3\,p)}\\
& -27\xi^2(1-\xi^2)
 \lambda_{qrr}\lambda_{112}\lambda_{122}\lambda_{qqr}
\chi_1(g_{(1\,r)})
\chi_{(3\,q)}(g_{(1\,q)})\chi_2(g_{(3\,r)})\\
&\hspace*{8cm}
\chi_{(3\,r)}(g_{q})
a_{(3\,p)}^3g_{qqr}\\
&+ 27(1-\xi^2) \lambda_{qrr}\lambda_{112}\lambda_{122}\lambda_{qqr}
\chi_{(3\,p)}(g_{q}) \chi_{q}(g_{1})
a_{(3\,p)}^3g_{qqr}\\
& +3 \lambda_{qrr}
\lambda_{112}\chi_{(1\,q)}(g_{(1\,p)})\chi_{(3\,p)}(g_{(2\,q)}) 
\chi_{(2\,r)}(g_{(1\,q)})
a_{(3\,p)}a_{(2\,r)} a_{(1\,q)}\\
& +9(1+\xi) \lambda_{qrr}
\lambda_{112}\lambda_{122}\chi_{(1\,q)}(g_{(1\,p)})\chi_{(3\,p)}(g_{(2\,q)}) 
\chi_{(2\,r)}(g_{(3\,q)})
\chi_{(3\,r)}(g_{(3\,p)})\\
&\hspace*{7cm}\chi_{(3\,q)}(g_{(3\,p)})
a_{(3\,r)} a_{(3\,q)}a_{(3\,p)}\\
& +3 \lambda_{qrr}
\lambda_{112}(1-\xi^2)\chi_{(1\,q)}(g_{(1\,p)})\chi_{(3\,p)}(g_{(2\,q)}) 
\chi_{(2\,q)}(g_{(1\,q)})
a_{(3\,p)}a_{(2\,q)}a_{(1\,r)}\\
&+6 \lambda_{qrr}\lambda_{112}\chi_{(1\,q)}(g_{(2\,p)})
a_{(3\,q)}a_{(2\,p)}a_{(1\,r)}.
\end{align*}

On the one hand, the terms involving $a_{(3\,p)}^3g_{qqr}$ cancel with each 
other, and so do the ones involving $a_{(3\,r)} a_{(3\,q)}a_{(3\,p)}$.
We use
\begin{align*}
&\lambda_{112}\chi_1(g_{(1\,r)})\chi_{(3\,r)}(g_{(1\,q)})=\lambda_{112}
\xi^2\chi_{(3\,r)}(g_{(2\,q)})\\
&\lambda_{112}\chi_{(1\,q)}(g_{(1\,p)})\chi_{(3\,p)}(g_{(2\,q)}) 
\chi_{(2\,q)}(g_{(1\,q)})=
\lambda_{112}\chi_{(3\,p)}(g_{12})
\end{align*}
 to simplify the scalars and we end up with \eqref{eqn:zpqr}. 

Similar computations lead to \eqref{eqn:zprq} and \eqref{eqn:zrpq}.
\epf


\end{document}